# THE BROWNIAN NET

BY RONGFENG SUN AND JAN M. SWART[1]

*TU Berlin and ÚTIA Prague*

The (standard) *Brownian web* is a collection of coalescing one-dimensional Brownian motions, starting from each point in space and time. It arises as the diffusive scaling limit of a collection of coalescing random walks. We show that it is possible to obtain a nontrivial limiting object if the random walks in addition branch with a small probability. We call the limiting object the *Brownian net*, and study some of its elementary properties.

## Contents

Received October 2006; revised June 2007.

[1]Supported by the DFG and by GAČR grants 201/06/1323 and 201/07/0237.

## 1. Introduction and main results.

1.1. *Arrow configurations and branching-coalescing random walks.* The Brownian web originated from the work of Arratia [1, 2], and has since been studied by Tóth and Werner [17], and Fontes, Isopi, Newman and Ravishankar [7, 8, 9]. It arises as the diffusive scaling limit of a collection of coalescing random walks. In this paper we show that it is possible to obtain a nontrivial limiting object if the random walks in addition branch with a small probability.

Let $\mathbb{Z}^2_{\mathrm{even}} := \{(x,t) : x, t \in \mathbb{Z},\ x+t \text{ is even}\}$ be the even sublattice of $\mathbb{Z}^2$. We interpret the first coordinate $x$ as space and the second coordinate $t$ as time, which is plotted vertically in figures. Fix a *branching probability* $\beta \in [0,1]$. Independently for each $(x,t) \in \mathbb{Z}^2_{\mathrm{even}}$, with probability $\frac{1-\beta}{2}$, draw an arrow from $(x,t)$ to $(x-1,t+1)$, with probability $\frac{1-\beta}{2}$, draw an arrow from $(x,t)$ to $(x+1,t+1)$, and with the remaining probability $\beta$, draw two arrows starting at $(x,t)$, one ending at $(x-1,t+1)$ and the other at $(x+1,t+1)$. (See Figure 1.) We denote the random configuration of all arrows by

$$\aleph_\beta := \{(z,z') \in \mathbb{Z}^2_{\mathrm{even}} \times \mathbb{Z}^2_{\mathrm{even}} : \text{ there is an arrow from } z \text{ to } z'\}. \tag{1.1}$$

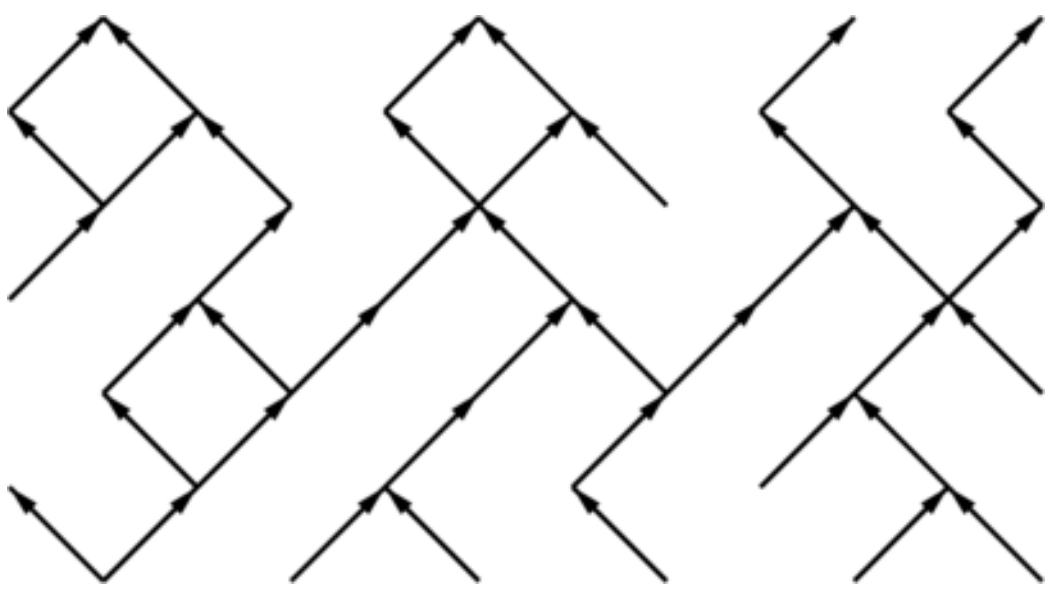

Fig. 1. *An arrow configuration.*

By definition, a *path along arrows in* $\aleph_\beta$, in short an $\aleph_\beta$*-path*, is the graph of a function $\pi:[\sigma_\pi,\infty]\to\mathbb{R}\cup\{*\}$, with $\sigma_\pi\in\mathbb{Z}\cup\{\pm\infty\}$, such that $((\pi(t),t),(\pi(t+1),t+1))\in\aleph_\beta$ and $\pi$ is linear on the interval $[t,t+1]$ for all $t\in[\sigma_\pi,\infty]\cap\mathbb{Z}$, while $\pi(\pm\infty)=*$ whenever $\pm\infty\in[\sigma_\pi,\infty]$. We call $\sigma_\pi$ the starting time, $\pi(\sigma_\pi)$ the starting position and $z_\pi:=(\pi(\sigma_\pi),\sigma_\pi)$ the starting point of the $\aleph_\beta$-path $\pi$.

For any $A\subset\mathbb{Z}^2_{\rm even}\cup\{(*,\pm\infty)\}$, we let $\mathcal{U}_\beta(A)$ denote the collection of all $\aleph_\beta$-paths with starting points in the set $A$, and we use the shorthands $\mathcal{U}_\beta(z):=\mathcal{U}_\beta(\{z\})$ and $\mathcal{U}_\beta:=\mathcal{U}_\beta(\mathbb{Z}^2_{\rm even}\cup\{(*,\pm\infty)\})$ for the collections of all $\aleph_\beta$-paths starting from a single point $z$ and from any point in space-time, respectively.

An arrow configuration $\aleph_\beta$ is in fact the graphical representation for a system of discrete time branching-coalescing random walks. Indeed, if we set

$$\eta^A_t:=\{\pi(t):\pi\in\mathcal{U}_\beta(A)\}\qquad(t\in\mathbb{Z},\ A\subset\mathbb{Z}^2_{\rm even}\cup\{(*,\pm\infty)\}),\tag{1.2}$$

and we interpret the points in $\eta^A_t$ as being occupied by a particle at time $t$, then $(\eta^A_t)_{t\in\mathbb{Z}}$ is a collection of random walks, which are introduced into the system at space-time points in $A$. At each time $t\in\mathbb{Z}$, independently each particle with probability $\frac{1-\beta}{2}$ jumps one step to the left (resp. right), and with probability $\beta$ branches into two particles, one jumping one step to the left and the other one step to the right. Two walks coalesce instantly when they jump to the same lattice site. Note that the case $\beta=0$ corresponds to coalescing random walks without branching.

We are interested in the limit of $\mathcal{U}_\beta$ under diffusive rescaling, letting at the same time $\beta\to0$. Thus, we rescale space by a factor $\varepsilon$, time by a factor $\varepsilon^2$, and let $\varepsilon\to0$ and $\beta\to0$ at the same time. For the case $\beta=0$, it has been shown in [8] that $\mathcal{U}_0$ diffusively rescaled converges weakly in law, with respect to an appropriate topology, to a random object $\mathcal{W}$, called the *Brownian web*. We will show that if $\beta/\varepsilon\to b$ for some $b\ge0$, then in (essentially) the same topology as in [8], $\mathcal{U}_\beta$ diffusively rescaled converges in law to a random object $\mathcal{N}_b$, which we call the *Brownian net with branching parameter* $b$. Here $\mathcal{N}_0$ is equal to $\mathcal{W}$ in distribution, while $\mathcal{N}_b$ with $b>0$ differ from $\mathcal{W}$, but are related to each other through scaling.

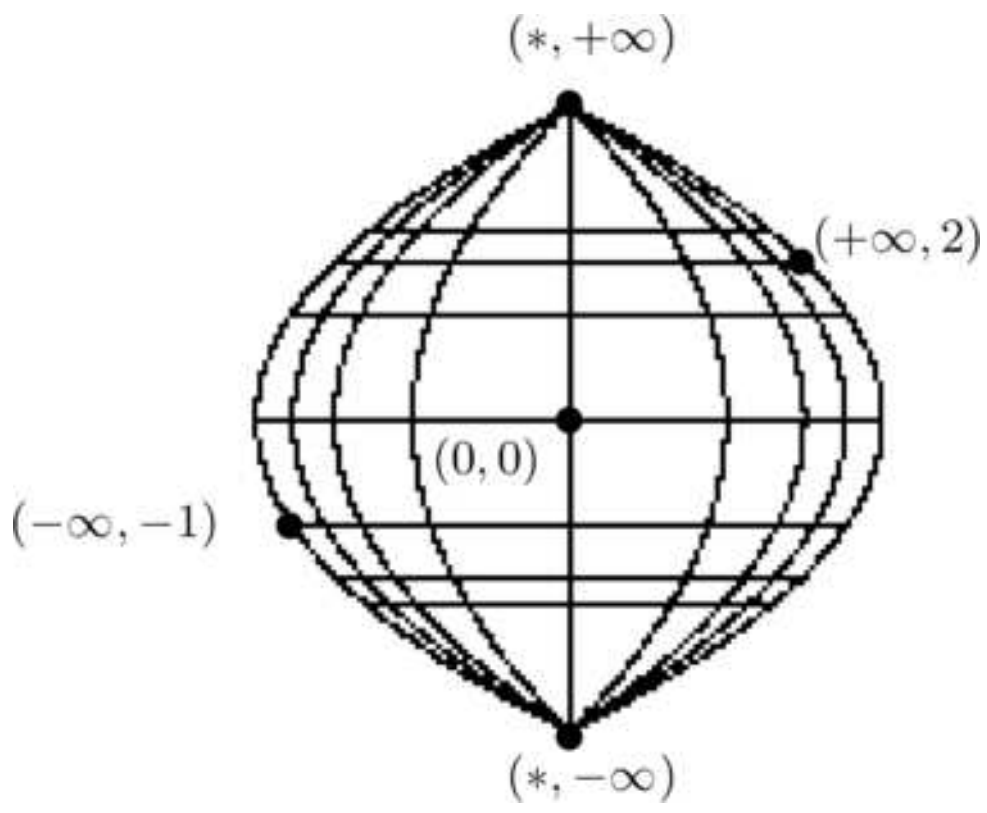

FIG. 2. *The compactification $R_c^2$ of $\mathbb{R}^2$.*

1.2. *Topology and convergence.* To formulate our main results, we first need to define the space in which our random variables take values and the topology with respect to which we will prove convergence. Our topology is essentially the same as the one used in [7, 8], except for a slight (and in most applications irrelevant) detail, as explained in the Appendix.

Let $R_c^2$ be the compactification of $\mathbb{R}^2$ obtained by equipping the set $R_c^2 := \mathbb{R}^2 \cup \{(\pm\infty, t) : t \in \mathbb{R}\} \cup \{(*, \pm\infty)\}$ with a topology such that $(x_n, t_n) \to (\pm\infty, t)$ if $x_n \to \pm\infty$ and $t_n \to t \in \mathbb{R}$, and $(x_n, t_n) \to (*, \pm\infty)$ if $t_n \to \pm\infty$ (regardless of the behavior of $x_n$). In [7, 8], such a compactification is achieved by taking the completion of $\mathbb{R}^2$ with respect to the metric

$$\rho((x_1,t_1),(x_2,t_2)) = |\Theta_1(x_1,t_1) - \Theta_1(x_2,t_2)| \vee |\Theta_2(t_1) - \Theta_2(t_2)|, \tag{1.3}$$

where the map $\Theta = (\Theta_1, \Theta_2)$ is defined by

$$\Theta(x,t) = (\Theta_1(x,t), \Theta_2(t)) := \left(\frac{\tanh(x)}{1+|t|}, \tanh(t)\right). \tag{1.4}$$

We can think of $R_c^2$ as the image of $[-\infty,\infty]^2$ under the map $\Theta$. Of course, $\rho$ and $\Theta$ are by no means the only choices that achieve the desired compactification. See Figure 2 for a picture of $R_c^2$ (for a somewhat different choice of $\Theta$).

By definition, a (continuous) *path* in $R_c^2$ is a function $\pi : [\sigma_\pi, \infty] \to [-\infty, \infty] \cup \{*\}$, with $\sigma_\pi \in [-\infty, \infty]$, such that $\pi : [\sigma_\pi, \infty] \cap \mathbb{R} \to [-\infty, \infty]$ is continuous, and $\pi(\pm\infty) = *$ whenever $\pm\infty \in [\sigma_\pi, \infty]$. Equivalently, if we identify $R_c^2$ with the image of $[-\infty,\infty]^2$ under the map $\Theta$, then $\pi$ is a continuous map from $[\Theta_2(\sigma_\pi), \Theta_2(\infty)]$ to $\mathbb{R}$ whose graph is contained in $\Theta([-\infty,\infty]^2)$. Throughout the paper *we identify a path* $\pi$ *with its graph* $\{(\pi(t),t) : t \in [\sigma_\pi,\infty]\} \subset R_c^2$. Thus, we often view paths as compact subsets of $R_c^2$. We stress that the starting time is part of the definition of a path, that is, paths that are defined by the same function but have different starting times are considered to be different. Note that both the function defining a path and its starting time can be read off from its graph.

We let $\Pi$ denote the space of all paths in $R_c^2$, equipped with the metric

$$\begin{aligned} d(\pi_1,\pi_2) := {} & |\Theta_2(\sigma_{\pi_1}) - \Theta_2(\sigma_{\pi_2})| \\ & \vee \sup_{t \geq \sigma_{\pi_1} \wedge \sigma_{\pi_2}} |\Theta_1(\pi_1(t \vee \sigma_{\pi_1}),t) - \Theta_1(\pi_2(t \vee \sigma_{\pi_2}),t)|. \end{aligned} \tag{1.5}$$

The space $(\Pi,d)$ is complete and separable. Note that paths converge in $(\Pi,d)$ if and only if their starting times converge and the functions converge locally uniformly on $\mathbb{R}$. If fact, one gets the same topology on $\Pi$ (though not the same uniform structure) if one views paths as compact subsets of $R_c^2$ and then equips $\Pi$ with the Hausdorff metric.

Recall that if $(E,d)$ is a metric space and $\mathcal{K}(E)$ is the space of all compact subsets of $E$, then the *Hausdorff metric* $d_{\rm H}$ on $\mathcal{K}(E)$ is defined by

$$d_{\rm H}(K_1,K_2) = \sup_{x_1 \in K_1} \inf_{x_2 \in K_2} d(x_1,x_2) \vee \sup_{x_2 \in K_2} \inf_{x_1 \in K_1} d(x_1,x_2). \tag{1.6}$$

If $(E,d)$ is complete and separable, then so is $(\mathcal{K}(E), d_{\rm H})$. For a given topology on $E$, the *Hausdorff topology* generated by $d_{\rm H}$ depends only on the topology on $E$ and not on the choice of the metric $d$.

The Brownian net $\mathcal{N}_b$ and web $\mathcal{W}$ are $\mathcal{K}(\Pi)$-valued random variables. We define scaling maps $S_\varepsilon : R_c^2 \to R_c^2$ by

$$S_\varepsilon(x,t) := (\varepsilon x, \varepsilon^2 t) \qquad ((x,t) \in R_c^2). \tag{1.7}$$

We adopt the convention that if $f : R_c^2 \to R_c^2$ and $A \subset R_c^2$, then $f(A) := \{f(x) : x \in A\}$ denotes the image of $A$ under $f$. Likewise, if $K$ is a set of subsets of $R_c^2$ (e.g., a set of paths), then $f(K) = \{f(A) : A \in K\}$ is the image of $K$ under the map $A \mapsto f(A)$. So, for example, $S_\varepsilon(\mathcal{U}_\beta)$ is the set of all $\aleph_\beta$-paths (viewed as subsets of $R_c^2$), diffusively rescaled with $\varepsilon$. This will later also apply to notation such as $-A := \{-x : x \in A\}$ and $A + y := \{x + y : x \in A\}$. We will sometimes also use the shorthand $f(A_1,\dots,A_n) := (f(A_1),\dots,f(A_n))$ when $f$ is a function defined on $R_c^2$ and $A_1,\dots,A_n$ are elements of, or subsets of, or sets of subsets of $R_c^2$.

Recall from Section 1.1 the definition of an arrow configuration $\aleph_\beta$ and the set $\mathcal{U}_\beta$ of all $\aleph_\beta$-paths. Note that $\mathcal{U}_\beta$ is a random subset of $\Pi$. In order to make $\mathcal{U}_\beta$ compact, *from now on, we modify our definition of* $\mathcal{U}_\beta$ by adding all trivial paths $\pi$ that satisfy $\sigma_\pi \in \{\pm\infty\} \cup \mathbb{Z}$ and $\pi(t) = -\infty$ or $\pi(t) = \infty$ for all $t \in [\sigma_\pi, \infty]$. The main result of this paper is the following convergence theorem.

THEOREM 1.1 (Convergence to the Brownian net). *There exist $\mathcal{K}(\Pi)$-valued random variables $\mathcal{N}_b$ $(b \geq 0)$ such that, if $\varepsilon_n, \beta_n \to 0$ and $\beta_n/\varepsilon_n \to b \geq 0$, then $S_{\varepsilon_n}(\mathcal{U}_{\beta_n})$ are $\mathcal{K}(\Pi)$-valued random variables, and*

$$\mathcal{L}(S_{\varepsilon_n}(\mathcal{U}_{\beta_n})) \underset{n\to\infty}{\Longrightarrow} \mathcal{L}(\mathcal{N}_b), \tag{1.8}$$

*where $\mathcal{L}(\,\cdot\,)$ denotes law, and $\Rightarrow$ denotes weak convergence. The random variables $(\mathcal{N}_b)_{b>0}$ satisfy the scaling relation*

$$\mathcal{L}(S_\varepsilon(\mathcal{N}_b)) = \mathcal{L}(\mathcal{N}_{b/\varepsilon}) \qquad (\varepsilon, b > 0). \tag{1.9}$$

*We have $\mathcal{L}(\mathcal{N}_0) = \mathcal{L}(\mathcal{W})$, where $\mathcal{W}$ is the Brownian web. However, the random variables $\mathcal{N}_b$ with $b > 0$ are different from $\mathcal{W}$.*

For $\beta_n = 0$, that is, the case without branching, Theorem 1.1 follows from [8], Theorem 6.1. In the next sections we will give three equivalent characterizations of the random variables $\mathcal{N}_b$ with $b > 0$. In view of the scaling relation (1.9), it suffices to consider the case $b = 1$. We call $\mathcal{N}_b$ the *Brownian net with branching parameter* $b$ and $\mathcal{N} := \mathcal{N}_1$ the (standard) *Brownian net*.

1.3. *The Brownian web.* In order to prepare for our first characterization of the Brownian net $\mathcal{N}$, we start by recalling from [8], Theorem 2.1, the characterization of the Brownian web $\mathcal{W}$. For any $K \in \mathcal{K}(\Pi)$ and $A \subset R^2_{\mathrm{c}}$, we let $K(A) := \{\pi \in K : z_\pi \in A\}$ denote the collection of paths in $K$ with starting points $z_\pi = (\pi(\sigma_\pi), \sigma_\pi)$ in $A$, and for $z \in R^2_{\mathrm{c}}$, we write $K(z) := K(\{z\})$.

THEOREM 1.2 (Characterization of the Brownian web). *There exists a $\mathcal{K}(\Pi)$-valued random variable $\mathcal{W}$, the so-called (standard) Brownian web, whose distribution is uniquely determined by the following properties:*

(i) *For each deterministic $z \in \mathbb{R}^2$, $\mathcal{W}(z)$ consists a.s. of a single path $\mathcal{W}(z) = \{\pi_z\}$.*

(ii) *For any finite deterministic set of points $z_1, \ldots, z_k \in \mathbb{R}^2$, $(\pi_{z_1}, \ldots, \pi_{z_k})$ is distributed as a system of coalescing Brownian motions starting at space-time points $z_1, \ldots, z_k$.*

(iii) *For any deterministic countable dense set $\mathcal{D} \subset \mathbb{R}^2$,*

$$\mathcal{W} = \overline{\mathcal{W}(\mathcal{D})} \qquad a.s., \tag{1.10}$$

*where $^{-}$ denotes closure in $(\Pi, d)$.*

Note that by properties (i) and (iii), for any deterministic countable dense set $\mathcal{D} \subset \mathbb{R}^2$, the Brownian web is almost surely determined by the countable system of paths $\mathcal{W}(\mathcal{D}) = \{\pi_z : z \in \mathcal{D}\}$, whose distribution is uniquely determined by property (ii). We call $\mathcal{W}(\mathcal{D})$ a *skeleton* of the Brownian web (relative to the countable dense set $\mathcal{D}$). Since skeletons may be constructed using

Kolmogorov's extension theorem, Theorem 1.2 allows a direct construction of the Brownian web.

Although $\mathcal{W}(z)$ consists of a single path for each deterministic $z \in \mathbb{R}^2$, as a result of the closure in (1.10), there exist random points $z$ where $\mathcal{W}(z)$ contains more than one path. These are points where the map $z \mapsto \pi_z$ is discontinuous, that is, the limit $\lim_{n\to\infty} \pi_{z_n}$ depends on the choice of the sequence $z_n \in \mathcal{D}$ with $z_n \to z$. These *special points* of the Brownian web are classified according to the number of disjoint incoming and distinct outgoing paths at $z$, and play an important role in understanding the Brownian web, and, later on, also the Brownian net. We recall the classification of the special points of the Brownian web in Section 3.2.

1.4. *Characterization of the Brownian net using hopping.* Our first characterization of the Brownian net will be similar to the characterization of the Brownian web in Theorem 1.2. A difficulty is that in the Brownian net $\mathcal{N}$, there is a multitude of paths starting at any site $z = (x,t) \in \mathbb{R}^2$. There is, however, a.s. a well-defined left-most path and right-most path in $\mathcal{N}(z)$, that is, there exist $l_z, r_z \in \mathcal{N}(z)$ such that $l_z(s) \le \pi(s) \le r_z(s)$ for any $s \ge t$ and $\pi \in \mathcal{N}(z)$. These left-most and right-most paths will play a key role in our characterization.

Our first task is to characterize the distribution of a finite number of left-most and right-most paths, started from deterministic starting points. Thus, for given deterministic $z_1, \ldots, z_k$, $z'_1, \ldots, z'_{k'} \in \mathbb{R}^2$, we need to characterize the joint law of $(l_{z_1}, \ldots, l_{z_k}, r_{z'_1}, \ldots, r_{z'_{k'}})$. It turns out that $(l_{z_1}, \ldots, l_{z_k})$ is a collection of coalescing Brownian motions with drift one to the left, while $(r_{z'_1}, \ldots, r_{z'_{k'}})$ is a collection of coalescing Brownian motions with drift one to the right. Moreover, paths evolve independently when they do not coincide. Therefore, in order to characterize the joint law of $(l_{z_1}, \ldots, l_{z_k}, r_{z'_1}, \ldots, r_{z'_{k'}})$, it suffices to characterize the interaction between one left-most path $l_z = l_{(x,s)}$ and one right-most path $r_{z'} = r_{(x',s')}$. The joint evolution of such a pair after time $s \vee s'$ can be characterized as the unique weak solution of the two-dimensional *left-right SDE*

$$
\begin{aligned}
dL_t &= 1_{\{L_t \neq R_t\}}\, dB^{\mathrm{l}}_t + 1_{\{L_t = R_t\}}\, dB^{\mathrm{s}}_t - dt,\\
dR_t &= 1_{\{L_t \neq R_t\}}\, dB^{\mathrm{r}}_t + 1_{\{L_t = R_t\}}\, dB^{\mathrm{s}}_t + dt,
\end{aligned}
\tag{1.11}
$$

where $B^{\mathrm{l}}_t, B^{\mathrm{r}}_t, B^{\mathrm{s}}_t$ are independent standard Brownian motions, and $L_t$ and $R_t$ are subject to the constraint that $L_t \le R_t$ for all $t \ge T := \inf\{u \ge s \vee s' : L_u \le R_u\}$. These rules uniquely determine the joint law of $(l_{z_1}, \ldots, l_{z_k}, r_{z'_1}, \ldots, r_{z'_{k'}})$. We call such a system a collection of *left-right coalescing Brownian motions*. See Figure 5 for a picture. We refer to Sections 2.1 and 2.2 for the proof that solutions to (1.11) are weakly unique, and a more careful definition of left-right coalescing Brownian motions.

Since we are not only interested in left-most and right-most paths, but in all paths in the Brownian net, we need a way to construct general paths from left-most and right-most paths. The method we choose in this section is based on *hopping*, that is, concatenating pieces of paths together at times when the two paths are at the same position.

We call $t$ an *intersection time* of two paths $\pi_1,\pi_2\in\Pi$ if $\sigma_{\pi_1}\vee\sigma_{\pi_2}<t<\infty$ and $\pi_1(t)=\pi_2(t)$. We say that a path $\pi_1$ crosses a path $\pi_2$ from left to right at time $t$ if there exist $\sigma_{\pi_1}\vee\sigma_{\pi_2}\leq t_-<t<t_+<\infty$ such that $\pi_1(t_-)<\pi_2(t_-)$, $\pi_2(t_+)<\pi_1(t_+)$, and $t=\inf\{s\in(t_-,t_+):\pi_2(s)<\pi_1(s)\}$. We say that $t\in\mathbb{R}$ is a *crossing time* of $\pi_1$ and $\pi_2$ if either $\pi_1$ crosses $\pi_2$ or $\pi_2$ crosses $\pi_1$ from left to right at time $t$.

For any collection of paths $\mathcal{A}\subset\Pi$, we let $\mathcal{H}_{\rm int}(\mathcal{A})$ denote the smallest set of paths containing $\mathcal{A}$ that is closed under hopping at intersection times, that is, $\mathcal{H}_{\rm int}(\mathcal{A})$ is the set of all paths $\pi\in\Pi$ of the form

$$\pi=\bigcup_{k=1}^{m}\{(\pi_k(t),t):t\in[t_{k-1},t_k]\}, \tag{1.12}$$

where $\pi_1,\ldots,\pi_m\in\mathcal{A}$, $\sigma_{\pi_1}=t_0<\cdots<t_m=\infty$, and $t_k$ is an intersection time of $\pi_k$ and $\pi_{k+1}$ for each $k=1,\ldots,m-1$. Likewise, we let $\mathcal{H}_{\rm cros}(\mathcal{A})$ denote the smallest set of paths containing $\mathcal{A}$ that is closed under hopping at crossing times.

THEOREM 1.3 (Characterization of the Brownian net using hopping). *There exists a $\mathcal{K}(\Pi)$-valued random variable $\mathcal{N}$, which we call the (standard) Brownian net, whose distribution is uniquely determined by the following properties:*

(i) *For each deterministic $z\in\mathbb{R}^2$, $\mathcal{N}(z)$ a.s. contains a unique left-most path $l_z$ and right-most path $r_z$.*

(ii) *For any finite deterministic set of points $z_1,\ldots,z_k,z'_1,\ldots,z'_{k'}\in\mathbb{R}^2$, the collection of paths $(l_{z_1},\ldots,l_{z_k},r_{z'_1},\ldots,r_{z'_{k'}})$ is distributed as a collection of left-right coalescing Brownian motions.*

(iii) *For any deterministic countable dense sets $\mathcal{D}^{\rm l},\mathcal{D}^{\rm r}\subset\mathbb{R}^2$,*

$$\mathcal{N}=\overline{\mathcal{H}_{\rm cros}(\{l_z:z\in\mathcal{D}^{\rm l}\}\cup\{r_z:z\in\mathcal{D}^{\rm r}\})}\qquad a.s. \tag{1.13}$$

Instead of hopping at crossing times, we could also have built our construction on hopping at intersection times. In fact, a much stronger statement is true.

PROPOSITION 1.4 (The Brownian net is closed under hopping). *We have $\mathcal{H}_{\rm int}(\mathcal{N})=\mathcal{N}$.*

We note, however, that as a result of the existence of special points in the Brownian web with one incoming and two outgoing paths, the Brownian net is *not closed* under hopping at times $t$ such that $\pi_1(t) = \pi_2(t)$ but $t = \sigma_{\pi_1} \vee \sigma_{\pi_2}(t)$. Thus, it is generally not allowed to hop onto paths at their starting times.

1.5. *The left-right Brownian web.* Given a Brownian net $\mathcal{N}$, if we take the closures of the sets of all left-most and right-most paths, started respectively from deterministic countable dense sets $\mathcal{D}^{\mathrm{l}}, \mathcal{D}^{\mathrm{r}} \subset \mathbb{R}^2$, then we obtain two Brownian webs, tilted respectively with drift $-1$ and $+1$, that are coupled in a special way. Our next theorem introduces this object in its own right.

THEOREM 1.5 (Characterization of the left-right Brownian web). *There exists a $\mathcal{K}(\Pi)^2$-valued random variable $(\mathcal{W}^{\mathrm{l}}, \mathcal{W}^{\mathrm{r}})$, which we call the (standard) left-right Brownian web, whose distribution is uniquely determined by the following properties:*

(i) *For each deterministic $z \in \mathbb{R}^2$, $\mathcal{W}^{\mathrm{l}}(z)$ and $\mathcal{W}^{\mathrm{r}}(z)$ a.s. each contain a single path $\mathcal{W}^{\mathrm{l}}(z) = \{l_z\}$ and $\mathcal{W}^{\mathrm{r}}(z) = \{r_z\}$.*

(ii) *For any finite deterministic set of points $z_1, \ldots, z_k, z'_1, \ldots, z'_{k'} \in \mathbb{R}^2$, the collection of paths $(l_{z_1}, \ldots, l_{z_k}; r_{z'_1}, \ldots, r_{z'_{k'}})$ is distributed as a collection of left-right coalescing Brownian motions.*

(iii) *For any deterministic countable dense sets $\mathcal{D}^{\mathrm{l}}, \mathcal{D}^{\mathrm{r}} \subset \mathbb{R}^2$,*

$$\mathcal{W}^{\mathrm{l}} = \overline{\{l_z : z \in \mathcal{D}^{\mathrm{l}}\}} \quad \text{and} \quad \mathcal{W}^{\mathrm{r}} = \overline{\{r_z : z \in \mathcal{D}^{\mathrm{r}}\}} \qquad \text{a.s.} \tag{1.14}$$

Note that if we define *titling maps* by $\mathrm{Tilt}^{\pm}(x,t) := (x \pm t, t)$, then $\mathrm{Tilt}^{+}(\mathcal{W}^{\mathrm{l}})$ and $\mathrm{Tilt}^{-}(\mathcal{W}^{\mathrm{r}})$ are distributed as the (standard) Brownian web. The following lemma, the proof of which can be found in Section 4, is an easy consequence of Theorem 1.3.

LEMMA 1.6 (Associated left-right Brownian web). *Let $\mathcal{N}$ be the Brownian net. Then $\mathcal{N}$ a.s. uniquely determines a left-right Brownian web $(\mathcal{W}^{\mathrm{l}}, \mathcal{W}^{\mathrm{r}})$ such that, for each deterministic $z \in \mathbb{R}^2$, $\mathcal{W}^{\mathrm{l}}(z) = \{l_z\}$ and $\mathcal{W}^{\mathrm{r}}(z) = \{r_z\}$, where $l_z$ and $r_z$ are respectively the left-most and right-most path in $\mathcal{N}(z)$.*

If $(\mathcal{W}^{\mathrm{l}}, \mathcal{W}^{\mathrm{r}})$ and $\mathcal{N}$ are coupled as in Lemma 1.6, then we say that $(\mathcal{W}^{\mathrm{l}}, \mathcal{W}^{\mathrm{r}})$ is the *left-right Brownian web associated with* the Brownian net $\mathcal{N}$. Theorem 1.3 shows that, conversely, a left-right Brownian web uniquely determines its associated Brownian net a.s.

In the next section we give another way to construct a Brownian net from its associated left-right Brownian web. Since the left-right Brownian web is characterized by Theorem 1.5, this yields another way to characterize the Brownian net.

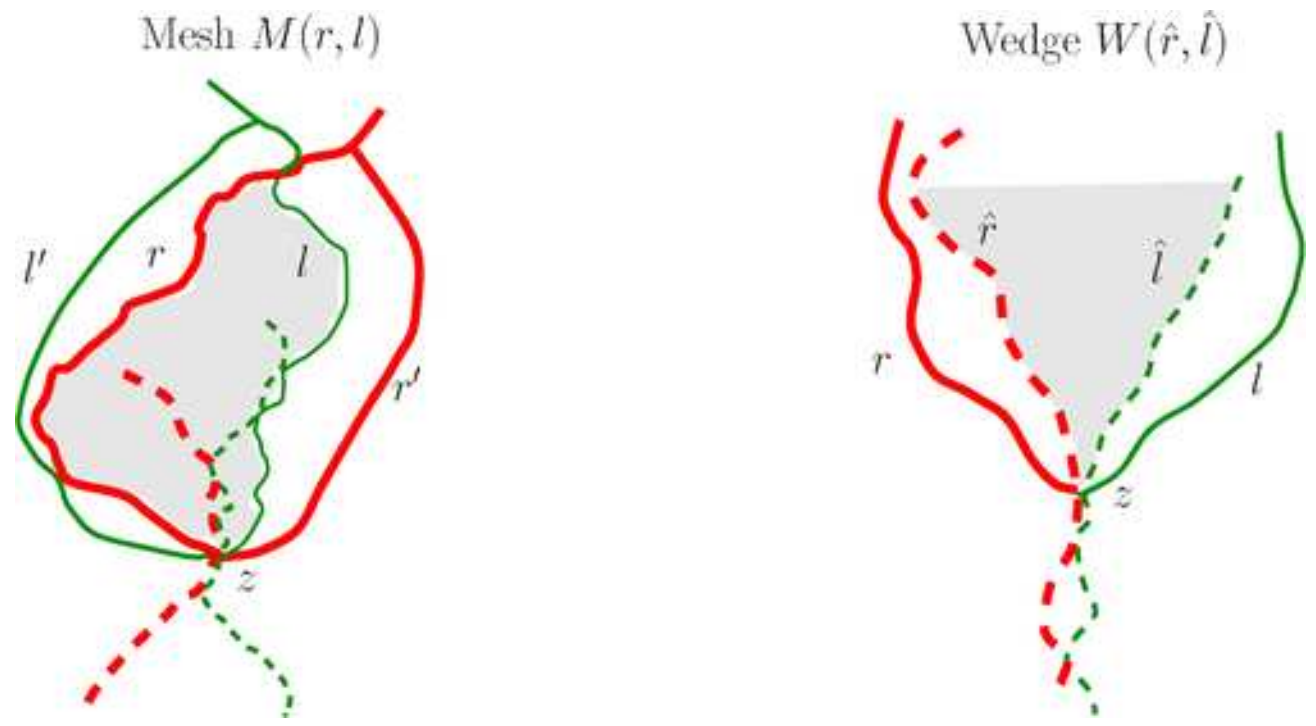

FIG. 3. *A mesh $M(r,l)$ with bottom point $z$ and a wedge $W(\hat{r},\hat{l})$ with bottom point $z$.*

1.6. *Characterization of the Brownian net using meshes.* If for some $z = (x,t) \in \mathbb{R}^2$, there exist $l \in \mathcal{W}^{\mathrm{l}}(z)$ and $r \in \mathcal{W}^{\mathrm{r}}(z)$ such that $r(s) < l(s)$ on $(t, t+\varepsilon)$ for some $\varepsilon > 0$, then denoting $T := \inf\{s > t : r(s) = l(s)\}$, we call the open set (see Figure 3)

$$M = M(r,l) := \{(y,s) \in \mathbb{R}^2 : t < s < T, r(s) < y < l(s)\} \tag{1.15}$$

the *mesh* with *bottom point* $z$, *top point* $(r(T), T)$, and left and right boundary $r$ and $l$, respectively. We call $x$ and $t$ the *bottom position* and *bottom time*, respectively, of the mesh $M$. We say that a path $\pi \in \Pi$ *enters* an open set $A \subset \mathbb{R}^2$ if there exist $\sigma_\pi < s < t$ such that $\pi(s) \notin A$ and $\pi(t) \in A$. Note the strict inequality in $s > \sigma_\pi$.

THEOREM 1.7 (Characterization of Brownian net with meshes). *Let $(\mathcal{W}^{\mathrm{l}}, \mathcal{W}^{\mathrm{r}})$ be the left-right Brownian web. Then almost surely,*

$$\begin{aligned} \mathcal{N} = \{\pi \in \Pi : \pi \text{ does not enter any mesh of} \\ (\mathcal{W}^{\mathrm{l}}, \mathcal{W}^{\mathrm{r}}) \text{ with bottom time } t > \sigma_\pi\} \end{aligned} \tag{1.16}$$

*is the Brownian net associated with $(\mathcal{W}^{\mathrm{l}}, \mathcal{W}^{\mathrm{r}})$.*

The next proposition implies that paths in the net $\mathcal{N}$ do not enter meshes of $(\mathcal{W}^{\mathrm{l}}, \mathcal{W}^{\mathrm{r}})$ at all (regardless of their bottom times), and hence, formula (1.16) stays true if one drops the restriction that the bottom time of the mesh should be larger than $\sigma_\pi$.

PROPOSITION 1.8 (Containment by left- and right-most paths). *Let $\mathcal{N}$ be the Brownian net and let $(\mathcal{W}^{\mathrm{l}}, \mathcal{W}^{\mathrm{r}})$ be its associated left-right Brownian web. Then, almost surely, there exist no $\pi \in \mathcal{N}$ and $l \in \mathcal{W}^{\mathrm{l}}$ such that $l(s) \leq \pi(s)$ and $\pi(t) < l(t)$ for some $\sigma_\pi \vee \sigma_l < s < t$. An analogue statement holds for right-most paths.*

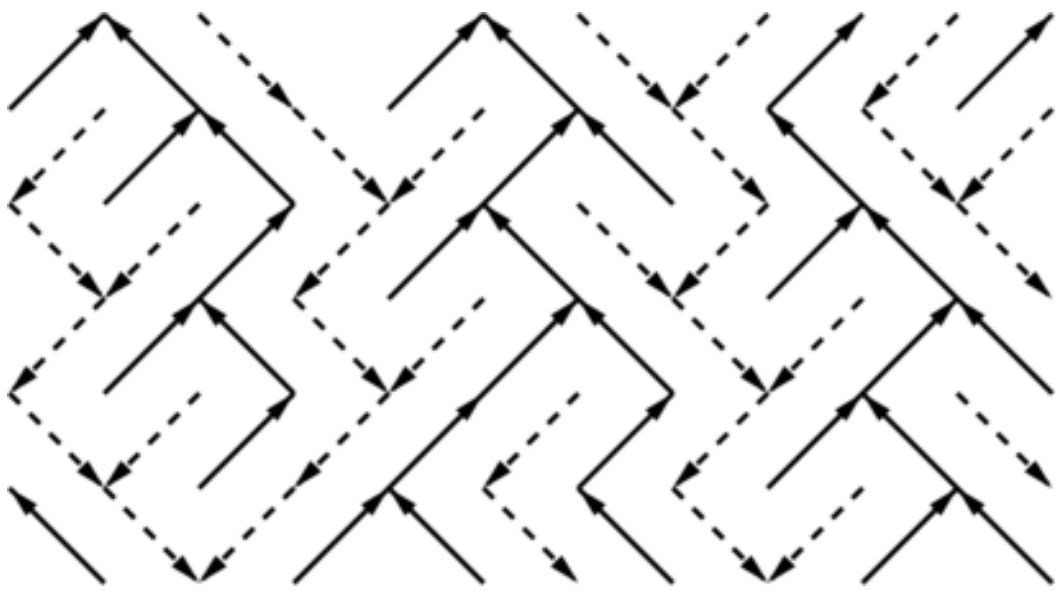

FIG. 4. *Dual arrow configuration with no branching.*

REMARK. Theorem 1.7 and Proposition 1.8 have analogues for the Brownian web. Indeed, generalizing our earlier definition, we can define a *left-right Brownian web* $(\mathcal{W}_b^{\mathrm{l}}, \mathcal{W}_b^{\mathrm{r}})$ *with drift* $b \geq 0$ by replacing the drift terms $+dt$ and $-dt$ in the left-right SDE (1.11) with $+b\,dt$ and $-b\,dt$, respectively. Then $\mathcal{W}_0^{\mathrm{l}} = \mathcal{N}_0 = \mathcal{W}_0^{\mathrm{r}}$ a.s. is distributed as the (standard) Brownian web, and Theorem 1.7 and Proposition 1.8 hold for any $b \geq 0$. The meshes of the Brownian web are called *bubbles* in [9].

1.7. *The dual Brownian web.* Arratia [1] observed that there is a natural dual for the arrow configuration $\aleph_0$, the graphical representation of discrete time coalescing simple random walks. More precisely, $\aleph_0$ uniquely determines a *dual arrow configuration* $\hat{\aleph}_0$ defined as follows (see Figure 4):

$$\begin{aligned} \hat{\aleph}_0 := \{((x, t+1), (x \pm 1, t)) \in \mathbb{Z}^2_{\mathrm{odd}} \times \mathbb{Z}^2_{\mathrm{odd}} : \\ ((x,t), (x \mp 1, t+1)) \in \aleph_0\}. \end{aligned} \tag{1.17}$$

Observe that directed edges in $\aleph_0$ and $\hat{\aleph}_0$ do not cross, and $\aleph_0$ and $\hat{\aleph}_0$ uniquely determine each other. A dual arrow configuration $\hat{\aleph}_0$ is the graphical representation of a system of coalescing simple random walks running backward in time, and $-\hat{\aleph}_0 + (0,1)$ is equally distributed with $\aleph_0$. In analogy with $\mathcal{U}_0$, let $\hat{\mathcal{U}}_0$ denote the set of backward paths along arrows in $\hat{\aleph}_0$. It follows from results in [8, 9] that

$$\mathcal{L}(S_\varepsilon(\mathcal{U}_0, \hat{\mathcal{U}}_0)) \underset{\varepsilon \to 0}{\Longrightarrow} \mathcal{L}(\mathcal{W}, \hat{\mathcal{W}}), \tag{1.18}$$

where $\mathcal{W}$ is the standard Brownian web, and $\hat{\mathcal{W}}$ is the so-called *dual Brownian web* associated with $\mathcal{W}$. One has

$$\mathcal{L}(-(\mathcal{W}, \hat{\mathcal{W}})) = \mathcal{L}(\hat{\mathcal{W}}, \mathcal{W}). \tag{1.19}$$

In particular, $\hat{\mathcal{W}}$ is equally distributed with $-\mathcal{W}$, the Brownian web rotated 180° around the origin. It was shown in [9, 15] that the interaction between paths in $\mathcal{W}$ and $\hat{\mathcal{W}}$ is that of Skorohod reflection.

A Brownian web $\mathcal{W}$ and its dual $\hat{\mathcal{W}}$ a.s. uniquely determine each other. There are several ways to construct $\mathcal{W}$ from $\hat{\mathcal{W}}$. We will describe one such way here, since this construction generalizes to the Brownian net. For any dual paths $\hat{\pi}_1, \hat{\pi}_2 \in \hat{\mathcal{W}}$ that are ordered as $\hat{\pi}_1(s) < \hat{\pi}_2(s)$ at the time $s := \hat{\sigma}_{\hat{\pi}_1} \wedge \hat{\sigma}_{\hat{\pi}_2}$, where $\hat{\sigma}_{\pi_i}$ denotes the starting time of $\hat{\pi}_i$ $(i = 1, 2)$, we let $T := \sup\{t < s : \hat{\pi}_1(t) = \hat{\pi}_2(t)\}$ denote the coalescence time of $\hat{\pi}_1$ and $\hat{\pi}_2$. We call the open set

(1.20) $$W = W(\hat{\pi}_1, \hat{\pi}_2) := \{(x, u) \in \mathbb{R}^2 : T < u < s, \hat{\pi}_1(u) < x < \hat{\pi}_2(u)\}$$

the *wedge* with left and right boundary $\hat{\pi}_1$ and $\hat{\pi}_2$. We say that a path $\pi \in \Pi$ *enters an open set* $A \subset \mathbb{R}^2$ *from outside* if there exist $\sigma_\pi < s < t$ such that $\pi(s) \notin \overline{A}$ and $\pi(t) \in A$.

THEOREM 1.9 (Construction of the Brownian web from its dual). *Let* $(\mathcal{W}, \hat{\mathcal{W}})$ *be a Brownian web and its dual. Then almost surely,*

(1.21) $$\mathcal{W} = \{\pi \in \Pi : \pi \textit{ does not enter any wedge of } \hat{\mathcal{W}} \textit{ from outside}\}.$$

The proof of Theorem 1.9 is contained in Section 4.2.

1.8. *Dual characterization of the Brownian net.* Let $(\mathcal{W}^{\mathrm{l}}, \mathcal{W}^{\mathrm{r}})$ be a left-right Brownian web. Then $\mathcal{W}^{\mathrm{l}}$ and $\mathcal{W}^{\mathrm{r}}$ each a.s. determine a dual web, which we denote respectively by $\hat{\mathcal{W}}^{\mathrm{l}}$ and $\hat{\mathcal{W}}^{\mathrm{r}}$. It will be proved in Section 5.2 below that

(1.22) $$\mathcal{L}(-(\mathcal{W}^{\mathrm{l}}, \mathcal{W}^{\mathrm{r}}, \hat{\mathcal{W}}^{\mathrm{l}}, \hat{\mathcal{W}}^{\mathrm{r}})) = \mathcal{L}(\hat{\mathcal{W}}^{\mathrm{l}}, \hat{\mathcal{W}}^{\mathrm{r}}, \mathcal{W}^{\mathrm{l}}, \mathcal{W}^{\mathrm{r}}).$$

In particular, the *dual left-right Brownian web* $(\hat{\mathcal{W}}^{\mathrm{l}}, \hat{\mathcal{W}}^{\mathrm{r}})$ is equally distributed with $-(\mathcal{W}^{\mathrm{l}}, \mathcal{W}^{\mathrm{r}})$, the left-right Brownian web rotated by 180° around the origin.

For any $\hat{r} \in \hat{\mathcal{W}}^{\mathrm{r}}$ and $\hat{l} \in \hat{\mathcal{W}}^{\mathrm{l}}$ that are ordered as $\hat{r}(s) < \hat{l}(s)$ at the time $s := \hat{\sigma}_{\hat{r}} \wedge \hat{\sigma}_{\hat{l}}$, we let $T := \sup\{t < s : \hat{r}(t) = \hat{l}(t)\}$ denote the first hitting time of $\hat{r}$ and $\hat{l}$, which may be $-\infty$. We call the open set (see Figure 3)

(1.23) $$W = W(\hat{r}, \hat{l}) := \{(x, u) \in \mathbb{R}^2 : T < u < s,\ \hat{r}(u) < x < \hat{l}(u)\}$$

the *wedge* with left and right boundary $\hat{r}$ and $\hat{l}$. The next theorem is analogous to Theorem 1.9.

THEOREM 1.10 (Dual characterization of the Brownian net). *Let* $(\mathcal{W}^{\mathrm{l}}, \mathcal{W}^{\mathrm{r}}, \hat{\mathcal{W}}^{\mathrm{l}}, \hat{\mathcal{W}}^{\mathrm{r}})$ *be a left-right Brownian web and its dual. Then, almost surely,*

(1.24) $$\mathcal{N} = \{\pi \in \Pi : \pi \textit{ does not enter any wedge of } (\hat{\mathcal{W}}^{\mathrm{l}}, \hat{\mathcal{W}}^{\mathrm{r}}) \textit{ from outside}\}$$

*is the Brownian net associated with* $(\mathcal{W}^{\mathrm{l}}, \mathcal{W}^{\mathrm{r}})$.

We note that there exist paths in $\mathcal{N}$ (even in $\mathcal{W}^{\rm l}$ and $\mathcal{W}^{\rm r}$) that enter wedges of $(\hat{\mathcal{W}}^{\rm l}, \hat{\mathcal{W}}^{\rm r})$ in the sense defined just before Theorem 1.7. Therefore, the condition in (1.24) that $\pi$ enters *from outside* cannot be relaxed.

1.9. *The branching-coalescing point set.* Just as the arrow configuration $\aleph_\beta$ is the graphical representation of a discrete system of branching-coalescing random walks, the Brownian net $\mathcal{N}$ is the graphical representation of a Markov process taking values in the space of compact subsets of $[-\infty, \infty]$, which we call the *branching-coalescing point set.* In analogy with (1.2), for any compact $A \subset R^2_{\rm c}$, we denote

$$\xi^A_t := \{\pi(t) : \pi \in \mathcal{N}(A)\} \qquad (t \in \mathbb{R},\ A \in \mathcal{K}(R^2_{\rm c})). \tag{1.25}$$

We set $\overline{\mathbb{R}} := [-\infty, \infty]$ and let $\mathcal{K}(\overline{\mathbb{R}})$ denote the space of compact subsets of $\overline{\mathbb{R}}$, equipped with the Hausdorff topology, under which $\mathcal{K}(\overline{\mathbb{R}})$ is itself a compact space. We recall that if $E$ is a compact metrizable space, then a *Feller process* in $E$ is a time-homogeneous Markov process in $E$, with cadlag sample paths, whose transition probabilities $P_t(x, dy)$ have the property that the map $(x, t) \mapsto P_t(x, \cdot\,)$ from $E \times [0, \infty)$ into the space of probability measures on $E$ is continuous with respect to the topology of weak convergence. Feller processes are strong Markov processes [6], Theorem 4.2.7.

THEOREM 1.11 (Branching-coalescing point set). *Let $\mathcal{N}$ be the Brownian net. Then for any $s \in \mathbb{R}$ and $K \in \mathcal{K}(\overline{\mathbb{R}})$,*

$$\xi_t := \xi^{K \times \{s\}}_t \qquad (s \le t < \infty) \tag{1.26}$$

*defines a Feller process $(\xi_t)_{t \ge s}$ in $\mathcal{K}(\overline{\mathbb{R}})$ with continuous sample paths, started from the initial state $K$ at time $s$. For each deterministic $t > s$, the set $\xi_t$ is a.s. locally finite in $\mathbb{R}$. If $K \in \mathcal{K}' := \{K \in \mathcal{K}(\overline{\mathbb{R}}) : K = \overline{K \cap \mathbb{R}}\}$, then*

$$\mathbb{P}[\xi_t \in \mathcal{K}' \ \forall t \ge s] = 1. \tag{1.27}$$

Note that $\mathcal{K}'$ excludes sets in which either $-\infty$ or $\infty$ is an isolated point, and hence, $\mathcal{K}'$ can in a natural way be identified with the space of all closed subsets of $\mathbb{R}$. Thus, property (1.27) says that we can view the branching-coalescing point set as a Markov process taking values in the space of closed subsets of $\mathbb{R}$.

The branching-coalescing point set $\xi_t$ arises as the scaling limit of the branching-coalescing random walks $\eta_t$ introduced in (1.2). The scaling regime considered in Theorem 1.1 allows us to interpret $\xi_t$ heuristically as a collection of Brownian particles which coalesce instantly when they meet but branch with an infinite rate. The infinite branching rate makes it difficult, however, to develop a good intuition from this simple picture. In particular, even for the process started at time zero from just one point, there is a dense collection of random times $t > 0$ such that $\xi_t$ is not locally finite. The proof of this fact is not difficult, but for lack of space, we defer it to a future paper.

For the branching-coalescing point set started from the whole extended real line $\overline{\mathbb{R}}$, we can explicitly calculate the expected density at any $t > 0$. Below, $|A|$ denotes the cardinality of a set and $\Phi(x) = \frac{1}{\sqrt{2\pi}} \int_{-\infty}^{x} e^{-y^2/2}\, dy$.

PROPOSITION 1.12 (Density of branching-coalescing point set). *We have*

$$\mathbb{E}[|\xi_t^{\overline{\mathbb{R}}\times\{0\}} \cap [a,b]|] = (b-a) \cdot \left( \frac{e^{-t}}{\sqrt{\pi t}} + 2\Phi(\sqrt{2t}) \right) \tag{1.28}$$

*for all* $[a,b] \subset \mathbb{R},\ t > 0$.

Note that the density of $\xi_t^{\overline{\mathbb{R}}\times\{0\}}$ is proportional to $t^{-1/2}$ as $t \downarrow 0$. This is consistent with the behavior of the Brownian web, but the decay is faster than is known for other coalescents such as Kingman's coalescent or the branching-coalescing particle systems in [3], Theorem 2(b). On the other hand, the density approaches the constant 2 as $t \to \infty$, in contrast to the Brownian web.

Our next proposition shows that it is possible to recover $\mathcal{N}(\overline{\mathbb{R}} \times \{0\})$ from $(\xi_t^{\overline{\mathbb{R}}\times\{0\}})_{t\geq 0}$. Below, for any $K \subset \mathcal{K}(R_c^2)$, we let

$$\cup K = \{z \in R_c^2 : \exists A \in K \text{ s.t. } z \in A\} \tag{1.29}$$

denote the union of all sets in $K$. We call $\cup K$ the *image set* of $K$. For $t \in [-\infty, \infty]$, let $\Pi_t := \{\pi \in \Pi : \sigma_\pi = t\}$ denote the space of all paths with starting time $t$. Note that $\cup(\mathcal{N} \cap \Pi_0) = \{(x,t) : t \geq 0, x \in \xi_t^{\overline{\mathbb{R}}\times\{0\}}\} \cup \{(*, \infty)\}$.

PROPOSITION 1.13 (Image set property). *Let $\mathcal{N}$ be the Brownian net. Then, almost surely for all $t \in [-\infty, \infty]$,*

$$\mathcal{N} \cap \Pi_t = \{\pi \in \Pi_t : \pi \subset \cup(\mathcal{N} \cap \Pi_t)\}. \tag{1.30}$$

1.10. *The backbone.* In this section we study $\mathcal{N}(*, -\infty)$, the set of paths in the Brownian net starting at time $-\infty$, and its discrete counterpart $\mathcal{U}_\beta(*, -\infty)$. These sets are relevant in the study of ergodic properties of the branching-coalescing point set and the branching-coalescing random walks. Borrowing terminology from branching theory, we call $\mathcal{N}(*, -\infty)$ and $\mathcal{U}(*, -\infty)$ respectively the *backbone* of the Brownian net and the *backbone* of an arrow configuration.

PROPOSITION 1.14 (Backbone of an arrow configuration). *For $\beta \geq 0$, the set of $\aleph_\beta$-paths, $\mathcal{U}_\beta$, satisfies the following properties:*

(i) $\{\pi(0): \pi \in \mathcal{U}_\beta(*,-\infty)\}$ *is a Bernoulli random field on* $\mathbb{Z}_{\text{even}}$ *with intensity* $\rho := \frac{4\beta}{(1+\beta)^2}$.

(ii) $\mathcal{U}_\beta(*,-\infty)$ *and* $-\mathcal{U}_\beta(*,-\infty)$ *are equal in law.*

(iii) *Almost surely,* $\mathcal{U}_\beta(x_n,t_n) \underset{n\to\infty}{\longrightarrow} \mathcal{U}_\beta(*,-\infty)$ *in* $\mathcal{K}(\Pi)$ *for any sequence* $(x_n,t_n) \in \mathbb{Z}^2_{\text{even}}$ *satisfying* $t_n \to -\infty$ *and* $\limsup_{n\to\infty} \frac{|x_n|}{|t_n|} < \beta$.

Note that [recall (1.2)]

$$\eta_t^{(*,-\infty)} = \{\pi(t): \pi \in \mathcal{U}_\beta(*,-\infty)\} \qquad (t \in \mathbb{Z}) \tag{1.31}$$

defines, modulo parity, a stationary system of branching-coalescing random walks $(\eta_t^{(*,-\infty)})_{t\in\mathbb{Z}}$. Thus, property (i) implies that, modulo parity, Bernoulli product measure with intensity $\frac{4\beta}{(1+\beta)^2}$ is an invariant measure for the branching-coalescing random walks with branching probability $\beta$. This is perhaps surprising, unless one is familiar with other branching-coalescing particle systems such as Schlögl models (see, e.g., [3, 5, 13]). Property (ii) says that this invariant law is, moreover, reversible in a rather strong sense. Note that an arrow configuration $\aleph_\beta$ is *not* symmetric with respect to time reversal, so this statement is not as obvious as it may seem. Property (iii) implies that the branching-coalescing random walks $(\eta_t)_{t\geq 0}$ exhibit complete convergence, that is, for any nonempty initial state $\eta_0 \subset \mathbb{Z}_{\text{even}}$, as $t\to\infty$, $\eta_{2t}$ (resp. $\eta_{2t+1}$) converges in law to a Bernoulli product measure on $\mathbb{Z}_{\text{even}}$ (resp. $\mathbb{Z}_{\text{odd}}$) with intensity $\rho = \frac{4\beta}{(1+\beta)^2}$.

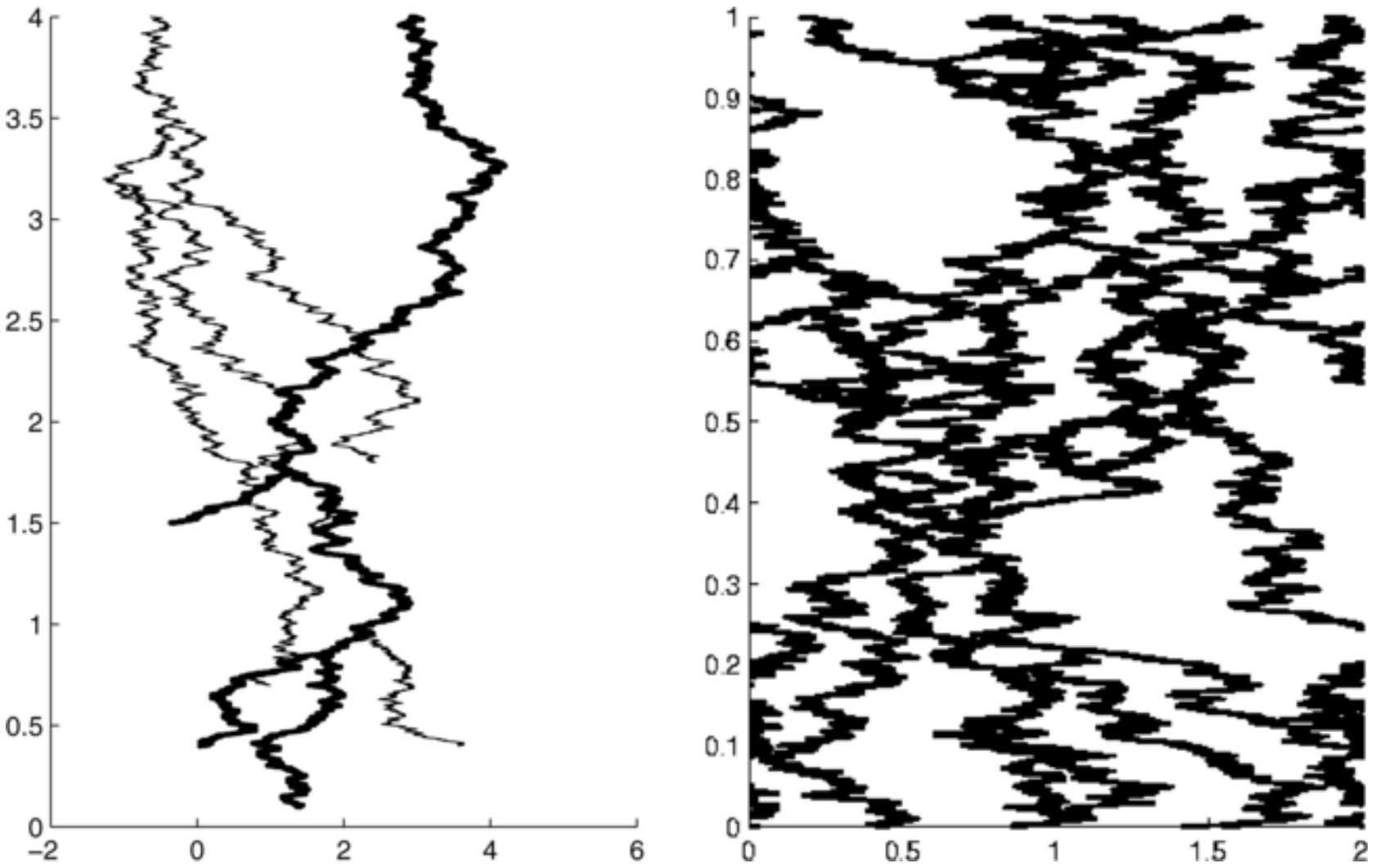

FIG. 5. *Left-right coalescing Brownian motions and the backbone of the Brownian net.*

For the Brownian net, we have the following analogue of Proposition 1.14.

PROPOSITION 1.15 (Backbone of the Brownian net). *The Brownian net $\mathcal{N}$ satisfies the following properties:*

(i) $\{\pi(0) : \pi \in \mathcal{N}(*, -\infty)\} \backslash \{\pm\infty\}$ *is a Poisson point process on* $\mathbb{R}$ *with intensity* 2.

(ii) $\mathcal{N}(*, -\infty)$ *and* $-\mathcal{N}(*, -\infty)$ *are equal in law.*

(iii) *Almost surely,* $\mathcal{N}(x_n, t_n) \underset{n\to\infty}{\longrightarrow} \mathcal{N}(*, -\infty)$ *in* $\mathcal{K}(\Pi)$ *for any sequence* $(x_n, t_n) \in \mathbb{R}^2$ *satisfying* $t_n \to -\infty$ *and* $\limsup_{n\to\infty} \frac{|x_n|}{|t_n|} < 1$.

In analogy with the branching-coalescing random walks, it follows that the law of a Poisson point set on $\mathbb{R}$ with intensity 2 is an invariant law for the branching-coalescing point set, that the latter exhibits complete convergence, and hence, this is its unique nontrivial invariant law. See Figure 5 for a picture of the backbone, or rather its image set $\cup\mathcal{N}(*, -\infty)$. Note that by Proposition 1.13, any path starting at time $-\infty$ that stays in $\cup\mathcal{N}(*, -\infty)$ is a path in $\mathcal{N}(*, -\infty)$.

1.11. *Discussion, applications and open problems.* This article began with the question of whether it is possible to add a small branching probability to the arrow configuration $\aleph_0$, which scales to the Brownian web, in such a way that one still obtains a nontrivial limit. At first sight, this may not seem possible because of the instantaneous coalescing of paths in the Brownian web. At second thought, for arrow configurations $\aleph_\beta$ with branching probability $\beta$, if we rescale space and time by $\varepsilon$ and $\varepsilon^2$ and let $\varepsilon \to 0$, then for the left-most and right-most $\aleph_\beta$-path starting from the origin to have a nontrivial limit, we need $\beta/\varepsilon \to b$ for some $b > 0$. It seems a coincidence that exactly the same scaling of $\beta$ and $\varepsilon$ is needed for the invariant measures of the branching-coalescing random walks from Proposition 1.14(i) to have a nontrivial limit. It was the observation of this coincidence that started off the present article.

Arratia's [1, 2] original motivation for studying the Brownian web came from one-dimensional voter models. In fact, coalescing simple random walks are dual to the one-dimensional nearest-neighbor voter model in two ways. They represent the genealogy lines of the voter model, and they also characterize the evolution of boundaries between domains of different types in an infinite type voter model. Voter models are used in population genetics to study the spread of genes in the absence of selection and mutation. They can also be viewed as the stochastic dynamics of zero-temperature one-dimensional Potts models. These points of view suggest several extensions of the Brownian web.

In [9] the *marked Brownian web* was introduced for the study of one-dimensional Potts models at small positive temperature. There, with small probability, a site may change its type, giving rise to a "nucleation event." In the biological context, such an event may be interpreted as a mutation. For the dual system of coalescing random walks, this results in a small death rate. The diffusive scaling limit of such a system is characterized by a Poisson marking of paths in the dual Brownian web, according to their length measure, where marks indicate deaths of particles.

There are at least two motivations for studying the Brownian net. First, in the biological interpretation, if instead of mutation, one adds a small selection rate, then one ends up with a biased voter model, which is dual to branching-coalescing random walks (compare [3]). Near the completion of this article, we learned that Newman, Ravishankar and Schertzer have been studying a differently motivated model that also leads to the Brownian net. Their model is a one-dimensional infinite-type Potts model, where, in contrast to the model in [9], nucleation events can only occur at the boundaries between different types. These boundaries then evolve as a system of continuous-time branching-coalescing random walks, which leads to the Brownian net. Rather than starting from the left-right Brownian web, their construction of the Brownian net is based on allowing hopping in the (standard) Brownian web at points that are chosen according to a Poisson marking of the set of intersection points between paths in the Brownian web $\mathcal{W}$ and its dual $\hat{\mathcal{W}}$. This construction will be published in [12].

There are a number of questions about the Brownian net that are worth investigating. First, we would like to give a complete classification of all special points in the Brownian net, in analogy with the classification of special points in the Brownian web. We have some results in that direction and will present them in a forthcoming paper [14]. An important ingredient in that work is to characterize the interaction between forward left-most and dual right-most paths, which can be used to give an alternative characterization of the left-right Brownian web not discussed in the present paper.

The second question regards the universality of the Brownian net and the branching-coalescing point set. Our convergence result is for the simplest system of branching-coalescing random walks. It is plausible that the same result holds for more general branching-coalescing systems, such as Schlögl models or the biased annihilating branching process from [16]. Related to this is the question of how to characterize the branching-coalescing point set by means of a generator or well-posed martingale problem.

The third question is to study the *marked Brownian net*, which can be defined by a Poisson marking of paths in the Brownian net in the same spirit as the marked Brownian web introduced in [9]. In the biological context, such a model describes genealogies in the presence of small selection and rare mutations. It can be shown that the resulting branching-coalescing point set with deaths undergoes a phase transition of contact-process type as the death rate is increased. This model might therefore be of some relevance in the study of the one-dimensional contact process.

Finally, it needs to be investigated how the Brownian net relates to certain other objects that have been introduced in the literature. In particular, it seems that a subclass of the stochastic flows of kernels introduced by Howitt and Warren in [10] is supported on the Brownian net. Also, it would be interesting to know if the branching-coalescing point set is related to some field theory used in theoretical physics. The physicist's way of viewing this process would probably be to say that these are coalescing Brownian motions with an infinite branching rate, but, due to the coalescence, most of this branching is not effective, so at macroscopic space scales one only observes the 'renormalized' branching rate, which is finite.

1.12. *Outline.* The rest of the paper is organized as follows. In Section 2 we construct and characterize the left-right Brownian web (Theorem 1.5) by first characterizing the left-right SDE and left-right coalescing Brownian motions described in Section 1.4. In Section 3 we establish some basic properties for the left-right SDE, recall some properties of the Brownian web and its dual, and prove some basic properties for the left-right Brownian web and its dual.

In Section 4 we prove the equivalence of the hopping construction (Theorem 1.3) and the dual construction (Theorem 1.10) of the Brownian net. In Section 5 we prove Theorem 1.1, our main convergence result. In fact, we prove something more: denoting the collections of all left-most and right-most paths in an arrow configuration $\aleph_\beta$ by $\mathcal{U}^{\rm l}_\beta$ and $\mathcal{U}^{\rm r}_\beta$, respectively, we show that $S_\varepsilon(\mathcal{U}^{\rm l}_\beta,\mathcal{U}^{\rm r}_\beta,\mathcal{U}_\beta)$ converges to $(\mathcal{W}^{\rm l},\mathcal{W}^{\rm r},\mathcal{N})$, where $(\mathcal{W}^{\rm l},\mathcal{W}^{\rm r})$ is a left-right Brownian web and $\mathcal{N}$ is the associated Brownian net. Here the hopping and dual characterizations of the Brownian net serve respectively as a stochastic lower and upper bound on limit points of $S_\varepsilon(\mathcal{U}_\beta)$.

In Section 6 we carry out two density calculations. The first of these yields Proposition 1.12, while the second estimates the density of the set of times when the left-most path starting at the origin first meets some path in the Brownian net starting at time 0 to the left of the origin. This second calculation is used in Section 7 to establish the characterization of the Brownian net using meshes (Theorem 1.7) and Proposition 1.8. These two results then in turn imply Propositions 1.4 and 1.13.

Finally, in Section 8 we prove Theorem 1.11 on the branching-coalescing point set, and in Section 9 we prove Propositions 1.14 and 1.15 on the backbones of arrow configurations and the Brownian net.

**2. The left-right Brownian web.** In Section 2.1 we characterize the *left-right SDE* described in Section 1.4 as the unique weak solution of the SDE (1.11). In Section 2.2 we give a rigorous definition of a collection of *left-right coalescing Brownian motions* described in Section 1.4. Finally, in Section 2.3 we construct the left-right Brownian web and prove Theorem 1.5.

2.1. *The left-right SDE.* Recall that a Markov transition probability kernel $P_t(x,dy)$ on a compact metrizable space has the Feller property if the map $(x,t)\mapsto P_t(x,\cdot)$ from $E\times[0,\infty)$ into the space of probability measures on $E$ is continuous with respect to the topology of weak convergence. Each Feller transition probability kernel gives rise to a strong Markov process with cadlag sample paths [6], Theorem 4.2.7. If $E$ is not compact, but locally compact, then let $E_\infty=E\cup\{\infty\}$ denote the one-point compactification of $E$. In this case, one says that a Markov transition probability kernel $P_t(x,dy)$ on $E$ has the Feller property if the extension of $P_t(x,dy)$ to $E_\infty$ defined by putting $P_t(\infty,\cdot):=\delta_\infty$ $(t\geq 0)$ has the Feller property. The corresponding Markov process is called a Feller process.

PROPOSITION 2.1 (Well-posedness and stickiness of the left-right SDE). *For each initial state $(L_0,R_0)\in\mathbb{R}^2$, there exists a unique weak solution to the SDE (1.11) subject to the constraint that $L_t\leq R_t$ for all $t\geq T:=\inf\{s\geq 0: L_s=R_s\}$. The family of solutions $\{(L_t,R_t)_{t\geq 0}\}_{(L_0,R_0)\in\mathbb{R}^2}$ defines a Feller process. The law of the total time that $L_t$ and $R_t$ are equal is given by*

$$\mathcal{L}\left(\int_0^\infty 1_{\{L_t=R_t\}}\,dt\right)=\mathcal{L}\left(\sup_{t\geq 0}\left(\frac{B_t}{\sqrt{2}}-t+\frac{(L_0-R_0)\wedge 0}{2}\right)\right), \tag{2.1}$$

*where $B_t$ is a standard Brownian motion (started at zero).*

Denote $R^2_\leq:=\{(x,y)\in\mathbb{R}^2: x\leq y\}$. A weak $R^2_\leq$-valued solution to (1.11) is a quintuple $(L,R,B^{\mathrm{l}},B^{\mathrm{r}},B^{\mathrm{s}})$, where $B^{\mathrm{l}},B^{\mathrm{r}},B^{\mathrm{s}}$ are independent Brownian motions and $(L,R)$ is a continuous, adapted $R^2_\leq$-valued process such that (1.11) holds in integral form (where the stochastic integrals are of Itô-type).

We rewrite the SDE (1.11) into a different equation, which has a pathwise unique solution. (In contrast, we believe that solutions to (1.11) are not pathwise unique; see [18] and the references therein for a similar equation where this has been proved.) Consider the following equation:

$$\begin{aligned}
&\text{(i)} && dL_t=d\tilde{B}^{\mathrm{l}}_{T_t}+d\tilde{B}^{\mathrm{s}}_{S_t}-dt,\\
&\text{(ii)} && dR_t=d\tilde{B}^{\mathrm{r}}_{T_t}+d\tilde{B}^{\mathrm{s}}_{S_t}+dt,\\
&\text{(iii)} && T_t+S_t=t,\\
&\text{(iv)} && \int_0^t 1_{\{L_s<R_s\}}\,dS_s=0.
\end{aligned} \tag{2.2}$$

Note that (2.2)(iv) says that $S_t$ increases only when $L_t = R_t$. By definition, by a weak $R^2_{\leq}$-valued solution to (2.2), we will mean a 7-tuple $(L, R, S, T, \tilde{B}^{\rm l}, \tilde{B}^{\rm r}, \tilde{B}^{\rm s})$, where $\tilde{B}^{\rm l}, \tilde{B}^{\rm r}, \tilde{B}^{\rm s}$ are independent Brownian motions, $S, T$ are nonnegative, nondecreasing, continuous, adapted processes such that (2.2)(iii) and (iv) hold, and $(L, R)$ is a continuous, adapted $R^2_{\leq}$-valued process such that (2.2)(i) and (ii) hold in integral form.

Proposition 2.1 follows from the following lemma.

Lemma 2.2 (Space-time SDE). (a) *There is a one-to-one correspondence in law between weak $R^2_{\leq}$-valued solutions of* (1.11) *and weak $R^2_{\leq}$-valued solutions of* (2.2).

(b) *For each initial state $(L_0, R_0) \in \mathbb{R}^2_{\leq}$, equation* (2.2) *has a pathwise unique solution.*

(c) *Solutions to* (2.2) *satisfy $S_t := \int_0^t 1_{\{L_s = R_s\}}\, ds$,*

$$S_t = \sup_{0 \leq s \leq T_t} \left( \tfrac{1}{2}(L_0 + \tilde{B}^{\rm l}_s - R_0 - \tilde{B}^{\rm r}_s) - s \right) \qquad a.s., \tag{2.3}$$

*and $\lim_{t \to \infty} T_t = \infty$.*

Proof of Proposition 2.1. Since $L_t$ and $R_t$ evolve independently until they meet, it suffices to consider the case $L_0 \leq R_0$. The existence and uniqueness of weak solutions to (1.11) under the given constraint follow from Lemma 2.2(a) and (b), while (2.1) follows from Lemma 2.2(c). To prove the Feller property, by the continuity of sample paths, it suffices to show that the law on path space of solutions to (1.11) depends continuously on the initial state. Since the first meeting time and position of two independent Brownian motions depend continuously on their initial states, it suffices to show continuity of $\mathbb{R}^2_{\leq}$-valued solutions to (2.2) in the initial state. Fix Brownian motions $\tilde{B}^{\rm l}, \tilde{B}^{\rm r}$ and $\tilde{B}^{\rm s}$, and let $(L^n, R^n, S^n, T^n)$ be a sequence of solutions to (2.2) with initial states $(L^n_0, R^n_0) = (l_n, r_n) \in \mathbb{R}^2_{\leq}$, such that $(l_n, r_n) \to (l, r) \in \mathbb{R}^2_{\leq}$. Since $L^n$ and $R^n$ are Brownian motions and $S^n, T^n$ increase with slope at most 1, the sequence $(L^n, R^n, S^n, T^n)$ is tight. It is not hard to see that any subsequential limit solves (2.2) (compare the proof of Proposition 5.1 in Section 5.1), and therefore, $(L^n, R^n)$ converges to the pathwise unique solution of (2.2) with initial state $(l, r)$. □

Proof of Lemma 2.2. We start with the proofs of parts (b) and (c). Our approach is to transform an equation with a sticky boundary into a SDE with immediate reflection, which is a standard technique to deal with sticky reflection. Given a solution to (2.2), put

$$\begin{aligned} D_t &:= R_t - L_t, \\ W_t &:= R_0 + \tilde{B}^{\rm r}_t - L_0 - \tilde{B}^{\rm l}_t. \end{aligned} \tag{2.4}$$

Then

$$dD_t = dW_{T_t} + 2\,dt. \tag{2.5}$$

It is easy to see from (2.2) that $D_t$ leaves 0 immediately, that is, there exist no $s<t$ such that $D_u = 0$ for all $u \in (s,t)$. Hence, by (2.2)(iii) and (iv), $T_t$ is strictly increasing and continuous in $t$. Making the random time change $\tau = T_t$, denoting the inverse of $T$ by $\tau \mapsto T_\tau^{-1}$, and writing $dt = dT_t + dS_t$, we can transform the equation for $D_t$ into

$$dD_{T_\tau^{-1}} = dW_\tau + 2\,d\tau + 2\,dS_{T_\tau^{-1}}, \tag{2.6}$$

where $D_{T_\tau^{-1}}$ is constrained to be nonnegative for all $\tau > 0$, and $2S_{T_\tau^{-1}}$ is a nondecreasing process that increases only when $D_{T_\tau^{-1}} = 0$. Equation (2.6) is an SDE with instant reflection, known as the Skorohod equation (see, e.g., Section 3.6.C of [11]). It can be solved (pathwise) uniquely for $2S_{T_\tau^{-1}}$, yielding

$$2S_{T_\tau^{-1}} = -\inf_{0\le s\le\tau}(W_s + 2s). \tag{2.7}$$

Time changing back and remembering the definition of $W$, we arrive at (2.3). By the fact that $S_t + T_t = t$, we find that

$$t = T_t + \sup_{0\le s\le T_t}\left(\tfrac{1}{2}(L_0 + \tilde{B}^{\mathrm{l}}_s - R_0 - \tilde{B}^{\mathrm{r}}_s) - s\right). \tag{2.8}$$

Since the function

$$\tau \mapsto \tau + \sup_{0\le s\le \tau}\left(\tfrac{1}{2}(L_0 + \tilde{B}^{\mathrm{l}}_s - R_0 - \tilde{B}^{\mathrm{r}}_s) - s\right) \tag{2.9}$$

is strictly increasing and continuous, it has a unique inverse, which is $t \mapsto T_t$. This proves that $S$ and $T$ are pathwise unique, and therefore, by (2.2)(i) and (ii), also $L$ and $R$ are pathwise unique.

Since the solution $D_{T_\tau^{-1}}$ of (2.6) spends zero Lebesgue time at 0, time-changing $\tau = T_s$, we see that

$$0 = \int_0^{T_t} 1_{\{D_{T_\tau^{-1}}=0\}}\,d\tau = \int_0^t 1_{\{D_s=0\}}\,dT_s. \tag{2.10}$$

By (2.2)(iii) and (iv), we conclude that $S_t = \int_0^t 1_{\{L_s=R_s\}}\,ds$ and $T_t = \int_0^t 1_{\{L_s<R_s\}}\,ds$. Finally, since $L$ and $R$ are Brownian motions with drift $-1$ and $+1$, respectively, there is a last time that $L$ and $R$ are equal, and therefore, $\lim_{t\to\infty} T_t = \infty$. This completes the proofs of parts (b) and (c).

To prove part (a), note that we have just proved that any solution to (2.2) solves the following equations:

$$\text{(i)}\quad dL_t = d\tilde{B}^{\mathrm{l}}_{T_t} + d\tilde{B}^{\mathrm{s}}_{S_t} - dt,$$

$$
(2.11)\qquad
\begin{aligned}
&\text{(ii)} && dR_t = d\tilde{B}^{\mathrm{r}}_{T_t} + d\tilde{B}^{\mathrm{s}}_{S_t} + dt,\\
&\text{(iii)} && T_t = \int_0^t 1_{\{L_s<R_s\}}\,ds,\\
&\text{(iv)} && S_t = \int_0^t 1_{\{L_s=R_s\}}\,ds.
\end{aligned}
$$

Conversely, solutions to (2.11) obviously solve (2.2).

Given a $\mathbb{R}^2_{\leq}$-valued solution to (2.2), setting

$$
(2.12)\qquad
\begin{aligned}
B^{\mathrm{l}}_t &:= \tilde{B}^{\mathrm{l}}_{T_t} + \int_0^t 1_{\{L_s=R_s\}}\,d\hat{B}^{\mathrm{l}}_s,\\
B^{\mathrm{r}}_t &:= \tilde{B}^{\mathrm{r}}_{T_t} + \int_0^t 1_{\{L_s=R_s\}}\,d\hat{B}^{\mathrm{r}}_s,\\
B^{\mathrm{s}}_t &:= \tilde{B}^{\mathrm{s}}_{S_t} + \int_0^t 1_{\{L_s<R_s\}}\,d\hat{B}^{\mathrm{s}}_s,
\end{aligned}
$$

where $\hat{B}^{\mathrm{l}}, \hat{B}^{\mathrm{r}}$ and $\hat{B}^{\mathrm{s}}$ are Brownian motions independent of each other and of $\tilde{B}^{\mathrm{l}}, \tilde{B}^{\mathrm{r}}$ and $\tilde{B}^{\mathrm{s}}$, yields a weak $\mathbb{R}^2_{\leq}$-valued solution to (1.11). Conversely, given a weak $\mathbb{R}^2_{\leq}$-valued solution to (1.11), let $S_t := \int_0^t 1_{\{L_s=R_s\}}\,ds$, $T_t := \int_0^t 1_{\{L_s<R_s\}}\,ds$, and

$$
(2.13)\qquad
\begin{aligned}
\tilde{B}^{\mathrm{l}}_{T_t} &:= \int_0^t 1_{\{L_s<R_s\}}\,dB^{\mathrm{l}}_t,\\
\tilde{B}^{\mathrm{r}}_{T_t} &:= \int_0^t 1_{\{L_s<R_s\}}\,dB^{\mathrm{r}}_t,\\
\tilde{B}^{\mathrm{s}}_{S_t} &:= \int_0^t 1_{\{L_s=R_s\}}\,dB^{\mathrm{s}}_t.
\end{aligned}
$$

Then $(\tilde{B}^{\mathrm{l}}_t)_{t\in[0,T_\infty)}$, $(\tilde{B}^{\mathrm{r}}_t)_{t\in[0,T_\infty)}$, and $(\tilde{B}^{\mathrm{s}}_t)_{t\in[0,S_\infty)}$ can be extended to independent Brownian motions defined for all $t \geq 0$, yielding a solution to (2.11). This completes the proof of part (a). $\square$

2.2. *Left-right coalescing Brownian motions.* In this section we give a rigorous definition of a collection $l_{z_1}, \ldots, l_{z_k}, r_{z'_1}, \ldots, r_{z'_{k'}}$ of paths of left-right coalescing Brownian motions, started at points $z_1, \ldots, z_k, z'_1, \ldots, z'_{k'} \in \mathbb{R}^2$. Write $z_i = (x_i, t_i)$ and $z'_i = (x'_i, t'_i)$. The times $t_1, \ldots, t_k, t'_1, \ldots, t'_{k'}$ divide $\mathbb{R}$ into a finite number of intervals. It suffices to define a Markov process that specifies the time evolution of the left-right coalescing Brownian motions during each such interval.

Thus, we need to construct a Markov process $(L_{1,t}, \ldots, L_{k,t}; R_{1,t}, \ldots, R_{k',t})_{t\geq 0}$ in $\mathbb{R}^{k+k'}$ such that $(L_{1,t}, \ldots, L_{k,t})$ and $(R_{1,t}, \ldots, R_{k',t})$ are each distributed as coalescing Brownian motions with drift $-1$ and $+1$ respectively, and the interaction between paths in $(L_{1,t}, \ldots, L_{k,t})$ and $(R_{1,t}, \ldots, R_{k',t})$ is that of the left-right SDE (1.11). Instead of characterizing the joint process $(L_{1,t}, \ldots, L_{k,t}; R_{1,t}, \ldots, R_{k',t})$ as the unique weak solution of a system of SDEs, which is rather laborious, we give an inductive construction using the distribution of $(L_t, R_t)$.

We first construct the system up to the first time two left Brownian motions coalesce, or two right Brownian motions coalesce, or a right Brownian motion hits a left Brownian motion from the left. In the last case, the right Brownian motion has to continue on the right of the left Brownian motion, so we call this a crossing. If our left and right coalescing Brownian motions are initially ordered as $LRRLRLRLLLRRLR$, say, then we partition them as $\{LR\}\{R\}\{LR\}\{LR\}\{L\}\{L\}\{LR\}\{R\}\{LR\}$, letting all pairs of a left Brownian motion followed by a right Brownian motion constitute a partition element with two members, and putting all remaining Brownian motions into partition elements with one member. We let the partition elements evolve independently until the first coalescing or crossing time. Here partition elements containing two members evolve according to the left-right SDE (1.11), while partition elements containing one member are just Brownian motions with drift $+1$ or $-1$. At the first coalescing or crossing time, we respectively coalesce or cross the motions that have hit each other, repartition the remaining Brownian motions and continue the process. Note that there can be at most $k + k'$ coalescence events and at most $kk'$ crossings, so this procedure is iterated at most finitely often and eventually leads to a single pair $(L, R)$.

The above construction uniquely defines the system of left-right coalescing Brownian motions $l_{z_1}, \ldots, l_{z_k}, r_{z'_1}, \ldots, r_{z'_{k'}}$. By the Feller property of coalescing Brownian motions and solutions to the left-right SDE, it is clear that the law of $(l_{z_1}, \ldots, l_{z_k}, r_{z'_1}, \ldots, r_{z'_{k'}})$ depends continuously on the starting points $z_1, \ldots, z_k, z'_1, \ldots, z'_{k'}$, and the marginal distribution of a subset of paths in $\{l_{z_1}, \ldots, l_{z_k}, r_{z'_1}, \ldots, r_{z'_{k'}}\}$ is also a system of left-right coalescing Brownian motions. This consistency property allows the definition of a countable system of left-right coalescing Brownian motions.

2.3. *The left-right Brownian web.* We now construct the left-right Brownian web and prove Theorem 1.5.

PROOF OF THEOREM 1.5. We first show uniqueness. Fix countable dense sets $\mathcal{D}^{\mathrm{l}}, \mathcal{D}^{\mathrm{r}} \subset \mathbb{R}^2$. Suppose that there exists a $\mathcal{K}(\Pi) \times \mathcal{K}(\Pi)$-valued random variable $(\mathcal{W}^{\mathrm{l}}, \mathcal{W}^{\mathrm{r}})$ satisfying properties (i)–(iii) in Theorem 1.5. By property (i), let $l_z$, $z \in \mathcal{D}^{\mathrm{l}}$, denote the almost sure unique element in $\mathcal{W}^{\mathrm{l}}$ starting from $z$, and let $r_z$, $z \in \mathcal{D}^{\mathrm{r}}$, denote the almost sure unique element in $\mathcal{W}^{\rm r}$ starting from $z$. Then by property (ii), $((l_z)_{z\in\mathcal{D}^{\rm l}},(r_z)_{z\in\mathcal{D}^{\rm r}})$ is a $\Pi^{\mathcal{D}^{\rm l}}\times\Pi^{\mathcal{D}^{\rm r}}$-valued random variable whose finite-dimensional distributions are that of left-right coalescing Brownian motions. Hence, the law of $((l_z)_{z\in\mathcal{D}^{\rm l}},(r_z)_{z\in\mathcal{D}^{\rm r}})$ is uniquely determined, and by property (iii), so is the law of $(\mathcal{W}^{\rm l},\mathcal{W}^{\rm r})$.

We now construct a $\mathcal{K}(\Pi)\times\mathcal{K}(\Pi)$-valued random variable $(\mathcal{W}^{\rm l},\mathcal{W}^{\rm r})$ satisfying properties (i)–(iii) in Theorem 1.5. By our construction of left-right coalescing Brownian motions in Section 2.2 and their consistency property when more paths are added, we can invoke Kolmogorov's extension theorem to assert that there exists a $\Pi^{\mathcal{D}^{\rm l}}\times\Pi^{\mathcal{D}^{\rm r}}$-valued random variable $((l_z)_{z\in\mathcal{D}^{\rm l}},(r_z)_{z\in\mathcal{D}^{\rm r}})$ whose finite-dimensional distributions are that of left-right coalescing Brownian motions. Define

$$\mathcal{W}^{\rm l}:=\overline{\{l_z:z\in\mathcal{D}^{\rm l}\}},\qquad \mathcal{W}^{\rm r}:=\overline{\{r_{z'}:z'\in\mathcal{D}^{\rm r}\}}. \tag{2.14}$$

Note that $\mathcal{W}^{\rm l}$ and $\mathcal{W}^{\rm r}$ are each distributed as a standard Brownian web with drift $-1$ and $+1$ respectively. Properties (i) and (iii) then follow from the analogous properties for the standard Brownian web. It only remains to show (ii). Let $\{u_1,\dots,u_k\}$ and $\{u'_1,\dots,u'_{k'}\}$ be deterministic finite subsets of $\mathbb{R}^2$. By (i), almost surely, a unique path $l_{u_i}\in\mathcal{W}^{\rm l}$ starts from each $u_i$, $1\le i\le k$, and a unique path $r_{u'_j}\in\mathcal{W}^{\rm r}$ starts from each $u'_j$, $1\le j\le k'$. Choose $z_{n,i}\in\mathcal{D}^{\rm l}$, $z'_{n,j}\in\mathcal{D}^{\rm r}$ such that $z_{n,i}\to u_i$ and $z'_{n,j}\to u'_j$ as $n\to\infty$. Since the Brownian web is a.s. continuous at deterministic points (see Proposition 3.2), we have $l_{z_{n,i}}\to l_{u_i}$ and $r_{z'_{n,j}}\to r_{u'_j}$ in $\Pi$, and hence,

$$\mathcal{L}(l_{z_{n,1}},\dots,l_{z_{n,k}},r_{z'_{n,1}},\dots,r_{z'_{n,k'}})\underset{n\to\infty}{\Longrightarrow}\mathcal{L}(l_{u_1},\dots,l_{u_k},r_{u'_1},\dots,r_{u'_{k'}}). \tag{2.15}$$

By the continuity of left-right coalescing Brownian motions in its starting points, it follows that $(l_{u_1},\dots,l_{u_k},r_{u'_1},\dots,r_{u'_{k'}})$ is distributed as a system of left-right coalescing Brownian motions starting from $u_1,\dots,u_k,u'_1,\dots,u'_{k'}$, verifying property (ii). $\square$

## 3. Properties of the left-right Brownian web.

In Sections 3.1–3.3 below we collect some properties of solutions to the left-right SDE, the Brownian web and its dual, and the left-right Brownian web and its dual, respectively.

### 3.1. *Properties of the left-right SDE.*

Recall that a set $X$ is *perfect* if $X$ is closed and $x\in\overline{X\backslash\{x\}}$ for all $x\in X$, that is, $X$ has no isolated points.

PROPOSITION 3.1 (Properties of the left-right SDE). *Let $(L_t,R_t)_{t\ge0}$ be the unique weak solution of the SDE (1.11) with initial condition $(L_0,R_0)\in\mathbb{R}^2$, subject to the constraint that $L_t\le R_t$ for all $t\ge T:=\inf\{s\ge0:L_s=R_s\}$. Let $I:=\{t\ge0:L_t=R_t\}$ and let $\mu_I$ be the measure on $\mathbb{R}$ defined by $\mu_I(A):=\ell(I\cap A)$, where $\ell$ denotes Lebesgue measure. Then:*

(a) *Almost surely, $I$ is a nowhere dense perfect set.*
(b) *Almost surely, $I$ is the support of $\mu_I$.*

PROOF. If $T=\infty$, the lemma is vacuous. Since $(L_t,R_t)_{t\geq 0}$ is a strong Markov process and $T$ is a stopping time, we may assume without loss of generality that $T=0$, that is, $L_0=R_0$. Define $W$ as in (2.4), put $\tilde{W}_\tau := W_\tau + 2\tau$ $(\tau\geq 0)$, and

$$X_\tau := \tilde{W}_\tau + R_\tau \qquad \text{where } R_\tau := -\inf_{0\leq s\leq\tau}\tilde{W}_s \ (\tau\geq 0). \tag{3.1}$$

Then $X$ is a Brownian motion with diffusion constant 2 and drift 2, instantaneously reflected at zero. It is well known (and not hard to prove) that $\{\tau\geq 0 : X_\tau = 0\}$ is the support of $dR$.

Setting $D_t := R_t - L_t$ $(t\geq 0)$, we see by (2.6), (2.7) and (2.11)(iii) that

$$D_t = X_{T_t} \qquad \text{where } T_t := \int_0^t 1_{\{D_s>0\}}\,ds \qquad (t\geq 0). \tag{3.2}$$

It follows that $I=\{t\geq 0: D_t=0\}$ is the image of $\{\tau\geq 0 : X_\tau=0\}$ under the map $\tau\mapsto T_\tau^{-1}$. Since by (2.7) and (2.11)(iv),

$$S_t = \int_0^t 1_{\{D_s=0\}}\,ds = \tfrac{1}{2}R_{T_t} \qquad (t\geq 0), \tag{3.3}$$

the measure $\mu_I$ is the image of the measure $\frac{1}{2}\,dR$ under the map $T^{-1}$. Since $T^{-1}$ is a continuous open map, it follows that $\mathrm{supp}(\mu_I) = T^{-1}(\mathrm{supp}(dR)) = T^{-1}(\{\tau\geq 0 : X_\tau = 0\}) = I$. This proves part (b). It follows that $I$ has no isolated points, that is, is perfect. To see that $I$ is nowhere dense, by the Markov property, it suffices to show that $D_t$ leaves the origin immediately. Indeed, setting $\sigma := \inf\{t\geq 0 : D_t > 0\}$ and using (2.2), we see that $0 = D_\sigma = \int_0^\sigma 2\,dt = 2\sigma$ a.s. This proves part (a). $\square$

3.2. *Properties of the Brownian web.* In this section we recall some properties of the standard Brownian web $\mathcal{W}$ and its dual $\hat{\mathcal{W}}$, which can all be found in [8, 9, 15, 17]. Recall that $\hat{\sigma}_{\hat{\pi}}$ denotes the starting time of a dual path $\hat{\pi}$. Thus, a dual path is a map $\hat{\pi} : [-\infty, \hat{\sigma}_{\hat{\pi}}] \to [-\infty,\infty] \cup \{*\}$ such that $\hat{\pi} : [-\infty,\hat{\sigma}_{\hat{\pi}}] \cap \mathbb{R} \to [-\infty,\infty]$ is continuous, and $\hat{\pi}(\pm\infty) := *$ whenever $\pm\infty \in [-\infty, \hat{\sigma}_{\hat{\pi}}]$.

PROPOSITION 3.2 (Properties of the Brownian web). *Let $\mathcal{W}$ be the Brownian web and $\hat{\mathcal{W}}$ its dual. Then:*

(a) *$(\mathcal{W},\hat{\mathcal{W}})$ is equally distributed with $-(\hat{\mathcal{W}},\mathcal{W})$.*

(b) *Almost surely, paths in $\mathcal{W}$ coalesce when they meet, that is, for each $\pi,\pi'\in\mathcal{W}$ and $t>\sigma_\pi\vee\sigma_{\pi'}$ such that $\pi(t)=\pi'(t)$, one has $\pi(s)=\pi'(s)$ for all $s\geq t$.*

(c) *Almost surely, paths and dual paths do not cross, that is, there exist no* $\pi \in \mathcal{W}$, $\hat{\pi} \in \hat{\mathcal{W}}$, *and* $s, t \in [\sigma_\pi, \hat{\sigma}_{\hat{\pi}}]$ *such that* $(\pi(s) - \hat{\pi}(s)) \cdot (\pi(t) - \hat{\pi}(t)) < 0$.

(d) *Almost surely, paths and dual paths spend zero Lebesgue time together, that is, we have* $\int_{\sigma_\pi}^{\hat{\sigma}_{\hat{\pi}}} 1_{\{\pi(t)=\hat{\pi}(t)\}}\,dt = 0$ *for all* $\pi \in \mathcal{W}$ *and* $\hat{\pi} \in \hat{\mathcal{W}}$.

(e) *Almost surely, for each point* $z = (x,t) \in \mathbb{R}^2$, $x_n^- \uparrow x$, $x_n^+ \downarrow x$, $\pi_n^- \in \mathcal{W}(x_n^-, t)$, *and* $\pi_n^+ \in \mathcal{W}(x_n^+, t)$, *the limits* $\pi_{z-} := \lim_{n\to\infty} \pi_n^-$ *and* $\pi_{z+} := \lim_{n\to\infty} \pi_n^+$ *exist and do not depend on the choice of* $\pi_n^- \in \mathcal{W}(x_n^-, t)$ *and* $\pi_n^+ \in \mathcal{W}(x_n^+, t)$.

Points $z \in \mathbb{R}^2$ in the Brownian web are classified according to the number of disjoint incoming and distinct outgoing paths at $z$. By definition, an *incoming path* at $z = (x,t)$ is a path $\pi \in \mathcal{W}$ such that $\sigma_\pi < t$ and $\pi(t) = x$. Two incoming paths $\pi, \pi'$ at $z$ are *equivalent* if $\pi = \pi'$ on $[s, \infty]$, for some $\sigma_\pi \vee \sigma_{\pi'} \le s < t$. Let $m_{\rm in}(z)$ denote the number of equivalence classes of incoming paths in $\mathcal{W}$ at $z$, and let $m_{\rm out}(z)$ denote the cardinality of $\mathcal{W}(z)$. Then $(m_{\rm in}(z), m_{\rm out}(z))$ is the *type* of the point $z$ in $\mathcal{W}$. Points of type $(1,2)$ are distinguished into points of type $(1,2)_{\rm l}$ and $(1,2)_{\rm r}$, according to whether the incoming path continues along the left or right of the two outgoing paths. We let $(\hat{m}_{\rm in}(z), \hat{m}_{\rm out}(z))$ denote the type of a point $z$ in $\hat{\mathcal{W}}$, which is defined to be the type of $-z$ in $-\hat{\mathcal{W}}$, the rotation of $\hat{\mathcal{W}}$ by 180° around the origin. We denote the joint type of $z$ with respect to $(\mathcal{W}, \hat{\mathcal{W}})$ by $(m_{\rm in}(z), m_{\rm out}(z))/(\hat{m}_{\rm in}(z), \hat{m}_{\rm out}(z))$. The next lemma, which was first established in [17] (see also [9], Theorems 3.11 and 3.14), classifies all points in $\mathbb{R}^2$ according to their types in $(\mathcal{W}, \hat{\mathcal{W}})$. Note the relations $\hat{m}_{\rm out} = m_{\rm in} + 1$ and $m_{\rm out} = \hat{m}_{\rm in} + 1$.

LEMMA 3.3 (Classification of points in the Brownian web).

(a) *Almost surely, all* $z \in \mathbb{R}^2$ *are in* $(\mathcal{W}, \hat{\mathcal{W}})$ *of one of the types* $(0,1)/(0,1)$, $(0,2)/(1,1)$, $(0,3)/(2,1)$, $(1,1)/(0,2)$, $(1,2)_{\rm l}/(1,2)_{\rm l}$, $(1,2)_{\rm r}/(1,2)_{\rm r}$ *and* $(2,1)/(0,3)$. *See Figure* 6.

(b) *For each deterministic* $t \in \mathbb{R}$, *almost surely each point on* $\mathbb{R} \times \{t\}$ *is of either type* $(0,1)$, $(0,2)$ *or* $(1,1)$ *in* $\mathcal{W}$.

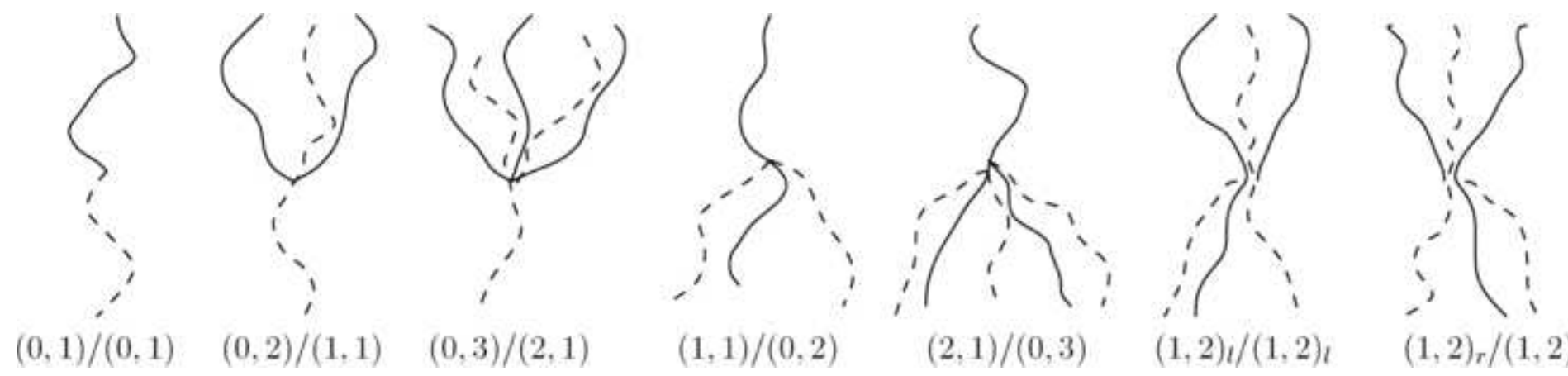

FIG. 6. *Types of points in the Brownian web and its dual* $(\mathcal{W}, \hat{\mathcal{W}})$.

(c) *Each deterministic point* $z \in \mathbb{R}^2$ *is almost surely of type* $(0,1)$ *in* $\mathcal{W}$.

The next lemma shows that convergent sequences of paths in $\mathcal{W}$ converge in a rather strong sense.

LEMMA 3.4 (Convergence of paths). *Let* $\mathcal{W}$ *be the standard Brownian web. Then:*

(a) *Almost surely, for any* $\{\pi_n\}_{n\in\mathbb{N}}, \pi \in \mathcal{W}$ *such that* $\pi_n \to \pi$, *one has* $\sigma_{\pi_n} \to \sigma_\pi$ *and* $\sup\{t \geq \sigma_{\pi_n} \vee \sigma_\pi : \pi_n(t) \neq \pi(t)\} \underset{n\to\infty}{\longrightarrow} \sigma_\pi$.

(b) *Let* $\mathcal{D}$ *be a deterministic countable dense subset of* $\mathbb{R}^2$. *Let* $\{\pi_z\}_{z\in\mathcal{D}}$ *be the skeleton of* $\mathcal{W}$ *relative to the starting set* $\mathcal{D}$. *Then almost surely, for all* $\pi \in \mathcal{W}$ *and* $\varepsilon > 0$, *there exists* $z = (x,t) \in \mathcal{D}$ *such that* $t \in (\sigma_\pi - \varepsilon, \sigma_\pi + \varepsilon)$ *and* $\pi_z(s) = \pi(s)$ *for all* $s \geq \sigma_\pi + \varepsilon$.

PROOF. By [8], Proposition 4.1, $\mathcal{W}^{t,\delta} := \{\gamma(t) : \gamma \in \mathcal{W}, \sigma_\gamma \leq t - \delta\}$ is a.s. locally finite for each $t, \delta \in \mathbb{Q}$ with $\delta > 0$. Therefore, $\pi_n \to \pi$ implies that, for each $\sigma_\pi < t \in \mathbb{Q}$, $\pi_n(t)$ eventually equals $\pi(t)$, and hence, $\pi_n = \pi$ on $[t, \infty)$, which proves part (a). Part (b) is a trivial consequence of part (a) and Theorem 1.2 (see also Proposition 2.2 of [17] and Proposition 4.2 of [8]). □

In applications of Lemma 3.4, one mostly needs part (b). Typically, a property is proved first for skeletal paths, and then extended to all paths in the web by Lemma 3.4(b).

We say that a path $\pi_1$ *crosses a path* $\pi_2$ *from left to right* if there exist $\sigma_{\pi_1} \vee \sigma_{\pi_2} \leq s < t$ such that $\pi_1(s) < \pi_2(s)$ and $\pi_2(t) < \pi_1(t)$. Likewise, we say that a path $\pi_1$ *crosses a dual path* $\hat{\pi}_2$ *from left to right* if there exist $\sigma_{\pi_1} \leq s < t \leq \hat{\sigma}_{\hat{\pi}_2}$ such that $\pi_1(s) < \hat{\pi}_2(s)$ and $\hat{\pi}_2(t) < \pi_1(t)$. The next lemma will be useful in what follows.

LEMMA 3.5 (Equivalence of crossing). *Let* $(\mathcal{W}, \hat{\mathcal{W}})$ *be the Brownian web and its dual. A path* $\gamma \in \Pi$ *crosses some* $\pi \in \mathcal{W}$ *from left to right if and only if it also crosses some* $\hat{\pi} \in \hat{\mathcal{W}}$ *from left to right. The same is true if we interchange left and right.*

PROOF. Assume $\gamma \in \Pi$ crosses $\pi \in \mathcal{W}$ from left to right, that is, $\gamma(s) < \pi(s)$ and $\gamma(t) > \pi(t)$ for some $\sigma_\gamma \vee \sigma_\pi \leq s < t$. Then by the noncrossing property of paths in $\mathcal{W}$ and $\hat{\mathcal{W}}$, for any $\hat{\pi} \in \hat{\mathcal{W}}(x,t)$ with $x \in (\pi(t), \gamma(t))$, we have $\gamma(s) < \pi(s) \leq \hat{\pi}(s)$. Hence, $\gamma$ crosses $\hat{\pi}$ from left to right. The proof of the converse implication is similar. By symmetry, the same statements hold for crossings from right to left. □

3.3. *Properties of the left-right Brownian web.* In this section we collect some basic properties of the left-right Brownian web $(\mathcal{W}^{\rm l}, \mathcal{W}^{\rm r})$ and its dual $(\hat{\mathcal{W}}^{\rm l}, \hat{\mathcal{W}}^{\rm r})$. Recall the definitions of intersection times and crossing times from Section 1.4. For any $\pi_1, \pi_2 \in \Pi$, we let

$$I(\pi_1, \pi_2) := \{t \in (\sigma_{\pi_1} \vee \sigma_{\pi_2}, \infty) : \pi_1(t) = \pi_2(t)\} \tag{3.4}$$

denote the set of intersection times of $\pi_1$ and $\pi_2$.

PROPOSITION 3.6 (Properties of the left-right Brownian web). *Let $(\mathcal{W}^{\rm l}, \mathcal{W}^{\rm r}, \hat{\mathcal{W}}^{\rm l}, \hat{\mathcal{W}}^{\rm r})$ be the standard left-right Brownian web and its dual. Then, almost surely, the following statements hold:*

(a) *For each $l \in \mathcal{W}^{\rm l}$ and $r \in \mathcal{W}^{\rm r}$ such that $\sigma_l \vee \sigma_r < \infty$, one has $T_{\rm cros} := \inf\{t > \sigma_l \vee \sigma_r : l(t) < r(t)\} = \inf\{t > \sigma_l \vee \sigma_r : l(t) \leq r(t)\} < \infty$, and $l(t) \leq r(t)$ for all $t \geq T_{\rm cros}$.*

(b) *For each $l \in \mathcal{W}^{\rm l}$ and $r \in \mathcal{W}^{\rm r}$, $\overline{I(l,r)}$ is a (possibly empty) nowhere dense perfect set.*

(c) *For each $l \in \mathcal{W}^{\rm l}$ and $r \in \mathcal{W}^{\rm l}$ such that $\sigma_l \vee \sigma_r < \infty$, the set $I(l,r)$ is the support of the measure $\mu_I$ on $(\sigma_l \vee \sigma_r, \infty)$ defined by $\mu_{I(l,r)}(A) := \ell(I(l,r) \cap A)$, where $\ell$ denotes Lebesgue measure.*

(d) *Paths in $\mathcal{W}^{\rm l}$ cannot cross paths in $\hat{\mathcal{W}}^{\rm r}$ from left to right, that is, there exist no $l \in \mathcal{W}^{\rm l}$, $\hat{r} \in \hat{\mathcal{W}}^{\rm r}$, and $\sigma_l \leq s < t \leq \hat{\sigma}_{\hat{r}}$ such that $l(s) < \hat{r}(s)$ and $\hat{r}(t) < l(t)$. Similarly, paths in $\mathcal{W}^{\rm r}$ cannot cross paths in $\hat{\mathcal{W}}^{\rm l}$ from right to left.*

PROOF. Let $\mathcal{D}^{\rm l}$ and $\mathcal{D}^{\rm r}$ be deterministic countable dense subsets of $\mathbb{R}^2$, and let $\{l_z\}_{z \in \mathcal{D}^{\rm l}}$ and $\{r_z\}_{z \in \mathcal{D}^{\rm r}}$ be the corresponding skeletons of $\mathcal{W}^{\rm l}$ and $\mathcal{W}^{\rm r}$. By Theorem 1.5 and Lemma 3.4(b), it suffices to prove parts (a)–(c) for skeletal paths, and hence for deterministic pairs $(l_z, r_{z'})$ where $z \in \mathcal{D}^{\rm l}$ and $z' \in \mathcal{D}^{\rm r}$. Since such deterministic pairs satisfy the SDE (1.11) by Theorem 1.5, parts (a)–(c) follow readily from Proposition 3.1(a) and (b). Property (d) is a consequence of (a) and Lemma 3.5. □

**4. The Brownian net.** Let $(\mathcal{W}^{\rm l}, \mathcal{W}^{\rm r}, \hat{\mathcal{W}}^{\rm l}, \hat{\mathcal{W}}^{\rm r})$ be a left-right Brownian web and its dual, and set

$$\mathcal{N}_{\rm hop} := \overline{\mathcal{H}_{\rm cros}(\mathcal{W}^{\rm l} \cup \mathcal{W}^{\rm r})}. \tag{4.1}$$

Note that if $\mathcal{D}^{\rm l}, \mathcal{D}^{\rm r} \subset \mathbb{R}^2$ are deterministic countable dense sets, then by Lemma 3.4(b), we also have $\mathcal{N}_{\rm hop} = \overline{\mathcal{H}_{\rm cros}(\mathcal{W}^{\rm l}(\mathcal{D}^{\rm l}) \cup \mathcal{W}^{\rm r}(\mathcal{D}^{\rm r}))}$. Define $\mathcal{N}_{\rm mesh}$ and $\mathcal{N}_{\rm wedge}$ by formulas (1.16) and (1.24), respectively. In Sections 4.1 and 4.2 we prove the inclusions $\mathcal{N}_{\rm hop} \subset \mathcal{N}_{\rm wedge}$ and $\mathcal{N}_{\rm wedge} \subset \mathcal{N}_{\rm hop}$, respectively. As an application, in Section 4.3 we establish Theorems 1.3 and 1.10, as well as Lemma 1.6. In addition, we prove Theorem 1.9 in Section 4.2, and, as a preparation for the characterization of the Brownian net using meshes, we prove the inclusion $\mathcal{N}_{\rm hop} \subset \mathcal{N}_{\rm mesh}$ in Section 4.1. The proof of the other inclusion is more difficult, and will be postponed to Section 7.

### 4.1. $\mathcal{N}_{\rm hop} \subset \mathcal{N}_{\rm wedge}$.

Set

$$\begin{aligned}\mathcal{P}_{\rm noncros} := \{\pi \in \Pi : \pi \text{ does not cross paths in } \mathcal{W}^{\rm l} \text{ from right}\\ \text{to left or paths in } \mathcal{W}^{\rm r} \text{ from left to right}\}.\end{aligned} \tag{4.2}$$

LEMMA 4.1 (Closedness of constructions). *The sets $\mathcal{N}_{\rm wedge}$, $\mathcal{N}_{\rm mesh}$ and $\mathcal{P}_{\rm noncros}$ are closed.*

PROOF. Note that if a path $\pi \in \Pi$ enters a mesh with bottom time $t > \sigma_\pi$, then it must enter from outside. Likewise, if $\pi$ crosses a dual path $\hat{l} \in \hat{\mathcal{W}}^{\rm l}$ from right to left, then it enters the open set $\{(x,t) \in \mathbb{R}^2 : t < \hat{\sigma}_{\hat{l}},\ x < \hat{l}(t)\}$ from outside. Thus, taking into account Lemma 3.5, all statements follow from the fact that if $\pi_n, \pi \in \Pi$ satisfy $\pi_n \to \pi$, and $\pi$ enters an open set $A$ from outside, then for $n$ sufficiently large, $\pi_n$ also enters $A$ from outside. □

LEMMA 4.2 (Noncrossing property). *We have $\mathcal{N}_{\rm hop} \subset \mathcal{P}_{\rm noncros}$ a.s.*

PROOF. It suffices to show that no path $\pi \in \mathcal{N}_{\rm hop}$ crosses paths in $\mathcal{W}^{\rm l}$ from right to left. By Lemma 4.1, it suffices to verify the statement for paths in $\mathcal{H}_{\rm cros}(\mathcal{W}^{\rm l} \cup \mathcal{W}^{\rm r})$. By Propositions 3.2(b) and 3.6(a), paths $\pi \in \mathcal{W}^{\rm l} \cup \mathcal{W}^{\rm r}$ have the stronger property that there exist no $\sigma_\pi < s < t$ and $l \in \mathcal{W}^{\rm l}$ such that $l(s) \le \pi(s)$ and $\pi(t) < l(t)$. It is easy to see that this stronger property is preserved under hopping. □

Let $A$ be either a mesh or wedge with (finite) bottom point $z = (x,t)$. We say that a path $\pi \in \Pi$ *enters $A$ through $z$* if $\sigma_\pi < t$ and there exists $s > t$ such that $(\pi(s), s) \in A$ and $(\pi(u), u) \in \overline{A}$ for all $u \in [t,s]$. Note that if a path enters a mesh (wedge) from outside, then it must either cross a left-most or right-most (dual) path in the wrong direction, or enter the mesh (wedge) through its bottom point.

LEMMA 4.3 (Finite wedges contained in meshes). *For every wedge $W$ with bottom point $z$, there exists a mesh $M$ with bottom point $z$ such that $W \subset M$.*

PROOF. Write $z = (x,t)$ and let $\hat{r}, \hat{l}$ be the left and right boundary of $W$. By Lemma 3.3, there exist $r \in \mathcal{W}^{\rm r}(z)$ and $l \in \mathcal{W}^{\rm l}(z)$ such that $r(s) \le \hat{r}(s)$ for all $s \in (t, \hat{\sigma}_{\hat{r}})$ and $\hat{l}(s) \le l(s)$ for all $s \in (t, \hat{\sigma}_{\hat{l}})$. It follows that $r$ and $l$ are the left and right boundary of a mesh containing $W$ (see Figure 3). □

LEMMA 4.4 (Hopping contained in mesh construction). *We have* $\mathcal{N}_{\rm hop} \subset \mathcal{N}_{\rm mesh}$ *a.s.*

PROOF. By Lemma 4.1, it suffices to show that $\mathcal{H}_{\rm cros}(\mathcal{W}^{\rm l} \cup \mathcal{W}^{\rm r}) \subset \mathcal{N}_{\rm mesh}$. We will show that, even stronger, paths in $\mathcal{H}_{\rm cros}(\mathcal{W}^{\rm l} \cup \mathcal{W}^{\rm r})$ do not enter meshes regardless of their bottom times. It is easy to see that this stronger property is preserved under hopping, so it suffices to show that paths in $\mathcal{W}^{\rm l} \cup \mathcal{W}^{\rm r}$ do not enter meshes. By symmetry, it suffices to show this for paths in $\mathcal{W}^{\rm l}$. By Propositions 3.2(b) and 3.6(a), it suffices to show that paths in $\mathcal{W}^{\rm l}$ cannot enter meshes through their bottom point. Let $M = M(r,l)$ be a mesh with left and right boundary $r$ and $l$ and bottom point $z = (x,t)$. Let $l' := l_{z-}$ and $r' := r_{z+}$ be the left-most path in $\mathcal{W}^{\rm l}(z)$ and the right-most path in $\mathcal{W}^{\rm r}(z)$, respectively, in the sense of Proposition 3.2(e). Then, by Proposition 3.6(a), $l'(s) \leq r(s)$ and $l(s) \leq r'(s)$ for all $s \geq t$ (see Figure 3). If some $l'' \in \mathcal{W}^{\rm l}$ enters $M$ through $z$, then by Lemma 3.3, $z$ must be of the type $(1,2)_{\rm l}$ or $(1,2)_{\rm r}$ in $\mathcal{W}^{\rm l}$, and therefore, $l''$ continues along either $l$ or $l'$. In either case, $l''$ does not enter $M$. □

LEMMA 4.5 (Hopping contained in wedge construction). *We have* $\mathcal{N}_{\rm hop} \subset \mathcal{N}_{\rm wedge}$ *a.s.*

PROOF. By Lemma 4.1, it suffices to show that $\mathcal{H}_{\rm cros}(\mathcal{W}^{\rm l} \cup \mathcal{W}^{\rm r}) \subset \mathcal{N}_{\rm wedge}$. Thus, we must show that paths in $\mathcal{H}_{\rm cros}(\mathcal{W}^{\rm l} \cup \mathcal{W}^{\rm r})$ do not cross paths in $\hat{\mathcal{W}}^{\rm l}, \hat{\mathcal{W}}^{\rm r}$ in the wrong direction or enter wedges through their bottom points. The first assertion follows from Lemmas 3.5 and 4.2, while the second assertion is a result of Lemmas 4.3 and 4.4. □

### 4.2. $\mathcal{N}_{\rm wedge} \subset \mathcal{N}_{\rm hop}$.

In this section we prove that $\mathcal{N}_{\rm wedge} \subset \mathcal{N}_{\rm hop}$. We start with a preparatory lemma.

LEMMA 4.6 (Compactness of $\mathcal{N}_{\rm hop}$). $\mathcal{N}_{\rm hop} \in \mathcal{K}(\Pi)$ *a.s.*

PROOF. Recall $(\Theta_1, \Theta_2)$ from (1.4). From the definition of the topology on $\Pi$ introduced in Section 1.2, by Arzela–Ascoli, we note that a set $K \subset \Pi$ is precompact if and only if the set of functions defined by the images of the graphs of $\pi \in K$ under the map $(\Theta_1, \Theta_2)$ is equicontinuous, that is, the modulus of continuity of $K$,

$$m_K(\delta) := \sup\{|\Theta_1(\pi(t),t) - \Theta_1(\pi(s),s)| : \pi \in K,\ s,t \geq \sigma_\pi,\ |\Theta_2(s) - \Theta_2(t)| \leq \delta\} \tag{4.3}$$

satisfies $m_K(\delta) \downarrow 0$ as $\delta \downarrow 0$.

Lemma 4.2 implies that, for each $\pi \in \mathcal{N}_{\rm hop}$ and $s \geq \sigma_\pi$, we have $l \leq \pi \leq r$ on $[s,\infty)$, where $l := l_{(\pi(s),s)-}$ and $r := r_{(\pi(s),s)+}$ denote respectively the left-most and the right-most paths in $\mathcal{W}^{\rm l}(\pi(s),s)$ and $\mathcal{W}^{\rm r}(\pi(s),s)$, in the sense of Proposition 3.2(e). It follows that, for any $t > s$,

$$\begin{aligned}&|\Theta_1(\pi(t),t) - \Theta_1(\pi(s),s)|\\ &\qquad \leq |\Theta_1(l(t),t) - \Theta_1(l(s),s)| \vee |\Theta_1(r(t),t) - \Theta_1(r(s),s)|.\end{aligned} \tag{4.4}$$

Taking the supremum over all $\pi \in \mathcal{N}_{\rm hop}$ and $\sigma_\pi \leq s < t$ such that $|\Theta_2(s) - \Theta_2(t)| \leq \delta$, we see that $m_{\mathcal{N}_{\rm hop}}(\delta) \leq m_{\mathcal{W}^{\rm l}\cup\mathcal{W}^{\rm r}}(\delta)$ (in fact, equality holds since $\mathcal{W}^{\rm l} \cup \mathcal{W}^{\rm r} \subset \mathcal{N}_{\rm hop}$), hence, the compactness of $\mathcal{N}_{\rm hop}$ follows from the compactness of $\mathcal{W}^{\rm l} \cup \mathcal{W}^{\rm r}$ a.s. $\square$

The next lemma is the main result of this section. This lemma and Proposition 1.8, which will be proved in Section 7, are the key technical results of this paper.

LEMMA 4.7 (Hopping contains wedge construction). *We have* $\mathcal{N}_{\rm wedge} \subset \mathcal{N}_{\rm hop}$ *a.s.*

PROOF. We must show that any path $\pi \in \mathcal{N}_{\rm wedge}$ can be approximated by a sequence of paths $\pi_n \in \mathcal{H}_{\rm cros}(\mathcal{W}^{\rm l} \cup \mathcal{W}^{\rm r})$. By the compactness of $\mathcal{N}_{\rm hop}$ (Lemma 4.6), it suffices to show that, for any $\pi \in \mathcal{N}_{\rm wedge}$, $\varepsilon > 0$, and $\sigma_\pi < t_1 < \cdots < t_n < \infty$, we can find $\pi^\varepsilon \in \mathcal{H}_{\rm cros}(\mathcal{W}^{\rm l} \cup \mathcal{W}^{\rm r})$ such that $\sigma_{\pi^\varepsilon} \in (\sigma_\pi, t_1)$ and $|\pi^\varepsilon(t_i) - \pi(t_i)| \leq \varepsilon$ for all $i = 1,\ldots,n$.

Our strategy is to first introduce piecewise continuous functions $\hat r$ and $\hat l$ on $[t_1,t_n]$, such that $\hat r(s) \leq \pi(s) \leq \hat l(s)$ for $s \in (t_1,t_n]$ and $|\hat r(t_i) - \pi(t_i)| \vee |\hat l(t_i) - \pi(t_i)| \leq \varepsilon$ for $i = 2,\ldots,n$. These functions will be constructed by piecing together paths in $\hat{\mathcal{W}}^{\rm r}$ and $\hat{\mathcal{W}}^{\rm l}$. We then construct $\pi^\varepsilon$ by *steering* a hopping path between $\hat r$ and $\hat l$.

We inductively choose $n = n_1 > \cdots > n_m > 1$ and $\hat r_1,\ldots,\hat r_m$ such that

$$\begin{aligned}&\hat r_k \in \hat{\mathcal{W}}^{\rm r}(\pi(t_{n_k}) - \varepsilon, t_{n_k}) \quad \text{and}\\ &n_{k+1} := \sup\{i : n_k > i > 1,\ \hat r_k(t_i) < \pi(t_i) - \varepsilon\}.\end{aligned} \tag{4.5}$$

This process terminates if $\hat r_k(t_i) \geq \pi(t_i) - \varepsilon$ for all $n_k > i > 1$. In this case we set $m := k$. We define $\hat r := \hat r_j$ on $(t_{n_{j+1}}, t_{n_j}]$ $(j = 1,\ldots,m-1)$ and $\hat r := \hat r_m$ on $[t_1, t_{n_m}]$. By left-right symmetry, we define $n = n'_1 > \cdots > n'_{m'} > 1$, $\hat l_1,\ldots,\hat l_{m'}$, and $\hat l$ analoguously. We make the following claims:

(1) $\hat r \leq \pi \leq \hat l$ on $[t_1,t_n]$.
(2) $\varepsilon' := \inf_{s\in[t_1,t_n]}(\hat l(s) - \hat r(s)) > 0$.
(3) $|\hat r(t_i) - \pi(t_i)| \vee |\hat l(t_i) - \pi(t_i)| \leq \varepsilon$ for $i = 2,\ldots,n$.

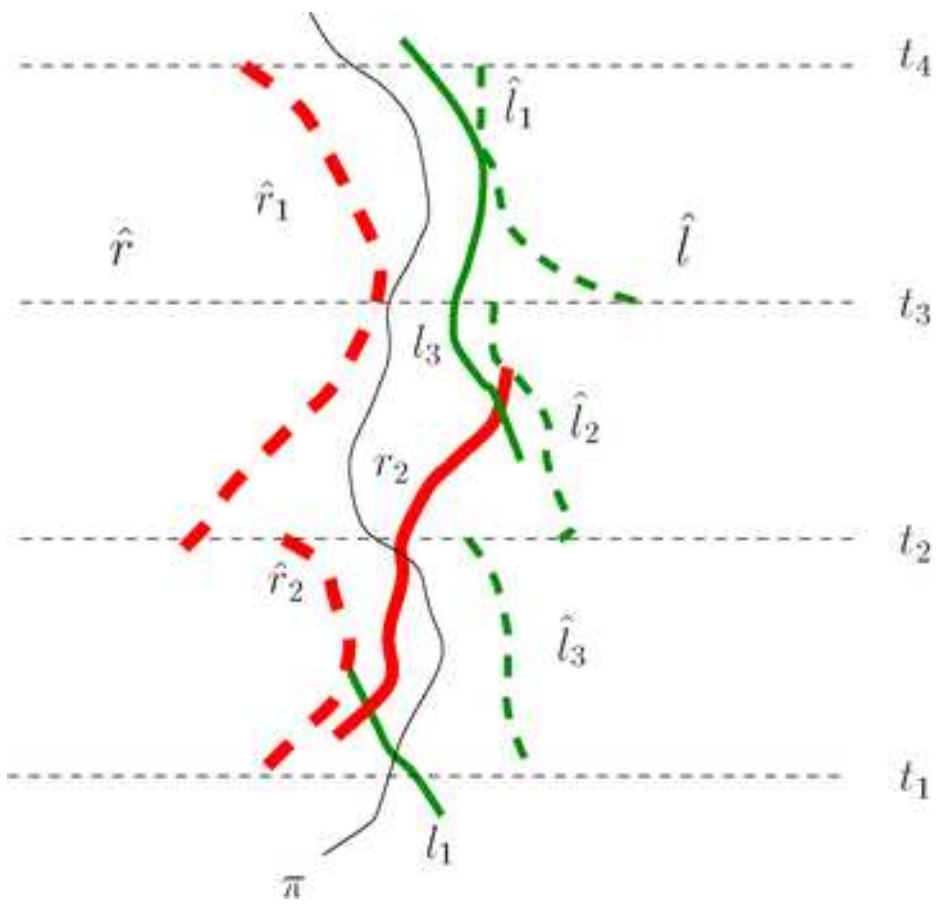

Fig. 7. *Steering a hopping path in the "fish-trap" $(\hat{r}, \hat{l})$.*

(4) $\lim_{t\downarrow t_i} \hat{r}(t) \leq \hat{r}(t_i)$ and $\lim_{t\downarrow t_i} \hat{l}(t) \geq \hat{l}(t_i)$ for $i = 2, \ldots, n-1$, which are the only possible discontinuities of $\hat{r}$ and $\hat{l}$.

Properties (1) and (2) follow from our assumption that $\pi$ does not enter wedges whose left and right boundaries are any of the dual paths $\hat{r}_1, \ldots, \hat{r}_m$ and $\hat{l}_1, \ldots, \hat{l}_{m'}$. Properties (3) and (4) are now obvious from our construction. The pair $(\hat{r}, \hat{l})$ resembles a fish-trap (see Figure 7).

We now construct a path $\pi^\varepsilon \in \mathcal{H}_{\mathrm{cros}}(\mathcal{W}^{\mathrm{l}} \cup \mathcal{W}^{\mathrm{r}})$ such that $\sigma_{\pi^\varepsilon} \in (\sigma_\pi, t_1)$, $|\pi^\varepsilon(t_1) - \pi(t_1)| \leq \varepsilon$, and $\hat{r}(s) \leq \pi^\varepsilon(s) \leq \hat{l}(s)$ for all $s \in [t_1, t_n]$. To this aim, we inductively choose $l_1, l_3, l_5, \ldots \in \mathcal{W}^{\mathrm{l}}$, $r_2, r_4, r_6, \ldots \in \mathcal{W}^{\mathrm{r}}$, and $\tau_1, \tau_2, \ldots$ such that $\tau_i$ is a crossing time of $l_i$ and $r_{i+1}$ if $i$ is odd and a crossing time of $r_i$ and $l_{i+1}$ if $i$ is even, in the following way. First, we choose $l_1$ such that $\sigma_{l_1} \in (\sigma_\pi, t_1)$ and $l_1(t_1) \in (\hat{r}(t_1), \hat{l}(t_1)) \cap [\pi(t_1) - \varepsilon, \pi(t_1) + \varepsilon]$. Assuming that we have already chosen $l_1, \ldots, l_i$ and $r_2, \ldots, r_{i-1}$, we proceed as follows. If $\hat{r}(s) < l_i(s) \leq \hat{l}(s)$ for all $s \in [\tau_{i-1}, t_n]$ (where $\tau_0 := t_1$), the process terminates. Otherwise, since paths cannot cross dual paths [Proposition 3.2(c)], $l_i$ must hit $\hat{r}$ before time $t_n$. In this case, we set $\sigma_i := \inf\{s \in [\tau_{i-1}, t_n] : l_i(s) = \hat{r}(s)\}$. Using Propositions 3.2(c) and 3.6(a), we can choose $\delta > 0$ sufficiently small and $r_{i+1} \in \mathcal{W}^{\mathrm{r}}$ started in $\{(x, s) : \sigma_i - \delta < s < \sigma_i,\ \hat{r}(s) < x < l_i(s)\}$, such that $r_{i+1}$ crosses $l_i$ at a time $\tau_i \in (\sigma_i - \delta, \sigma_i)$ and $r_{i+1}(\tau_i) - \hat{r}(\tau_i) \leq \frac{1}{3}\varepsilon'$. In case the last path we have chosen is a right-most path, by left-right symmetry, we proceed analogously. This process must terminate after a finite number of steps, for if this were not the case, then $\tau_i \uparrow \tau_\infty$ for some $\tau_\infty \leq t_n$. By the piecewise continuity of $\hat{l}$ and $\hat{r}$, we have $|r_i(\tau_i) - r_i(\tau_{i-1})| \geq \frac{1}{4}\varepsilon'$ for all sufficiently large even $i$, which contradicts the local equicontinuity, and hence compactness of $\mathcal{W}^{\mathrm{r}}$.

Defining $\pi^\varepsilon \in \mathcal{H}_{\rm cros}(\mathcal{W}^{\rm l} \cup \mathcal{W}^{\rm r})$ by hopping between the paths $l_1, l_3, \ldots$ and $r_2, r_4, \ldots$ at the times $\tau_1, \tau_2, \ldots$, we have found the desired approximation of $\pi$ by hopping paths. □

Since it is very similar to the proof of Lemma 4.7, we include here the proof of Theorem 1.9.

Proof of Theorem 1.9. Let $\mathcal{W}_{\rm wedge}$ be defined by the right-hand side of (1.21). Since paths in $\mathcal{W}$ cannot cross paths in $\hat{\mathcal{W}}$, to show that $\mathcal{W} \subset \mathcal{W}_{\rm wedge}$, it suffices that paths in $\mathcal{W}$ cannot enter wedges of $\hat{\mathcal{W}}$ through their bottom points. This can be proved by mimicking the proofs of Lemmas 4.3 and 4.4.

The inclusion $\mathcal{W}_{\rm wedge} \subset \mathcal{W}$ can be proved in the same way as the proof of Lemma 4.7. Since $\mathcal{W}$ is compact, it suffices to show that path that does not enter wedges from outside can be approximated by paths in $\mathcal{W}$. We can define a "fish-trap" whose left and right boundary are constructed by piecing dual paths together. In this case, any path in $\mathcal{W}$ entering the fish-trap from below must stay between its left and right boundary, so no hopping is necessary. □

4.3. *Characterizations with hopping and wedges.*

Proofs of Theorem 1.3, Lemma 1.6 and Theorem 1.10. Consider a left-right Brownian web and its dual $(\mathcal{W}^{\rm l}, \mathcal{W}^{\rm r}, \hat{\mathcal{W}}^{\rm l}, \hat{\mathcal{W}}^{\rm r})$, and let $\mathcal{N}_{\rm wedge}$ be defined as in (1.24) and $\mathcal{N}_{\rm hop}$ be defined as in (4.1). By Lemmas 4.5 and 4.7, $\mathcal{N}_{\rm hop} = \mathcal{N}_{\rm wedge}$. It follows from Lemma 4.2 that, for every $z = (x,t) \in \mathbb{R}^2$, we have $l_{z-}(s) \le \pi(s) \le r_{z+}(s)$ for all $\pi \in \mathcal{N}_{\rm hop}(z)$ and $s \ge t$, where $l_{z-}, r_{z+}$ are defined for $\mathcal{W}^{\rm l}, \mathcal{W}^{\rm r}$ as in Proposition 3.2(e). In particular, for deterministic $z$, the a.s. unique paths $l_z \in \mathcal{W}^{\rm l}(z)$ and $r_z \in \mathcal{W}^{\rm r}(z)$ are respectively the left-most and right-most paths in $\mathcal{N}_{\rm hop}(z)$. Setting $\mathcal{N} := \mathcal{N}_{\rm hop} = \mathcal{N}_{\rm wedge}$, we have found a $\mathcal{K}(\Pi)$-valued (by Lemma 4.6) random variable that satisfies conditions (i)–(ii) of Theorem 1.3. To see that condition (iii) is also satisfied, note that by Lemma 3.4(b), $\mathcal{N}_{\rm hop} = \overline{\mathcal{H}_{\rm cros}(\mathcal{W}^{\rm l}(\mathcal{D}^{\rm l}) \cup \mathcal{W}^{\rm r}(\mathcal{D}^{\rm r}))}$ for any deterministic countable dense sets $\mathcal{D}^{\rm l}, \mathcal{D}^{\rm r} \subset \mathbb{R}^2$. Since a random variable satisfying the conditions of Theorem 1.3 is obviously unique in distribution, the proof of Theorem 1.3 is complete.

Since for each deterministic $z$, the a.s. unique paths $l_z \in \mathcal{W}^{\rm l}(z)$ and $r_z \in \mathcal{W}^{\rm r}(z)$ are the left-most and right-most paths in $\mathcal{N}$, this also shows that to each Brownian net, there exists an associated left-right Brownian web, which is obviously unique by properties (i) and (ii) of Theorem 1.3. This proves Lemma 1.6.

Finally, since $\mathcal{N} = \mathcal{N}_{\rm wedge}$, we have also proved Theorem 1.10. □

**5. Convergence.** In this section we prove Theorem 1.1. In fact, we prove something more: we prove the joint convergence under diffusive scaling of the collections of all left-most and right-most paths (and their dual) in the arrow configuration $\aleph_\beta$ to the left-right Brownian web (and its dual), and of the collection of all $\aleph_\beta$-paths to the associated Brownian net. Throughout this section, $\mathcal{N}$ denotes the (standard) Brownian net, defined by the hopping or dual characterization (Theorem 1.3 or 1.10), which have been shown to be equivalent. We will not use the mesh characterization of the Brownian net (Theorem 1.7, yet to be proved) in this section.

In Section 5.1 we prove the convergence of a single pair of left-most and right-most paths in the arrow configuration $\aleph_\beta$ to a solution of the left-right SDE (1.11). In Section 5.2 we prove the convergence of all left-most and right-most paths and their dual to the left-right Brownian web and its dual. Finally, in Section 5.3, we prove Theorem 1.1.

5.1. *Convergence to the left-right SDE.* Recall the definition of $\aleph_\beta$ and $\mathcal{U}_\beta$ from Section 1.1. Let $\mathcal{U}^{\rm l}_\beta$ (resp. $\mathcal{U}^{\rm r}_\beta$) denote the set of left-most (resp. right-most) paths in $\mathcal{U}_\beta$, that is, $\aleph_\beta$-paths which follow arrows to the left (resp. right) at branching points. We have the following convergence result for a single pair of paths in $(\mathcal{U}^{\rm l}_{\beta_n}, \mathcal{U}^{\rm r}_{\beta_n})$. Below, $\mathcal{C}_{\mathbb{R}^n}[0,\infty)$ denotes the space of continuous functions from $[0,\infty)$ to $\mathbb{R}^n$, equipped with the topology of uniform convergence on compacta.

PROPOSITION 5.1 (Convergence of a pair of left and right paths). *Let $\beta_n, \varepsilon_n \to 0$ with $\beta_n/\varepsilon_n \to 1$. Let $x^{(n)}, y^{(n)} \in \mathbb{Z}_{\rm even}$ be points such that $(\varepsilon_n x^{(n)}, \varepsilon_n y^{(n)}) \to (x,y)$ for some $(x,y) \in \mathbb{R}^2$. Let $(L^{(n)}_t)_{t\geq 0}$ denote the path in $\mathcal{U}^{\rm l}_{\beta_n}$ starting at $(x^{(n)}, 0)$, and $(R^{(n)}_t)_{t\geq 0}$ the path in $\mathcal{U}^{\rm r}_{\beta_n}$ starting at $(y^{(n)}, 0)$. Then*

$$\mathcal{L}((\varepsilon_n L^{(n)}_{t/\varepsilon_n^2}, \varepsilon_n R^{(n)}_{t/\varepsilon_n^2})_{t\geq 0}) \underset{n\to\infty}{\Longrightarrow} \mathcal{L}((L_t, R_t)_{t\geq 0}), \tag{5.1}$$

*where $\Rightarrow$ denotes weak convergence of probability laws on $\mathcal{C}_{\mathbb{R}^2}[0,\infty)$, and $(L_t, R_t)_{t\geq 0}$ is the unique weak solution of (1.11) with initial state $(L_0, R_0) = (x,y)$, subject to the constraint that $L_t \leq R_t$ for all $t \geq T := \inf\{s \geq 0 : L_s = R_s\}$.*

PROOF. Set $T_n := \inf\{s \geq 0 : L^{(n)}_s = R^{(n)}_s\}$. Since up to time $T_n$, $L^{(n)}$ and $R^{(n)}$ are independent random walks with drift $-\beta_n$ and $+\beta_n$ respectively, it follows from Donsker's invariance principle and the almost sure continuity of the first intersection time between two independent Brownian motions with drift $\pm 1$ that

$$\mathcal{L}((\varepsilon_n L^{(n)}_{t/\varepsilon_n^2 \wedge T_n}, \varepsilon_n R^{(n)}_{t/\varepsilon_n^2 \wedge T_n})_{t\geq 0}) \underset{n\to\infty}{\Longrightarrow} \mathcal{L}((L_{t\wedge T}, R_{t\wedge T})_{t\geq 0}). \tag{5.2}$$

Therefore, it suffices to prove Proposition 5.1 for the case $x^{(n)} = y^{(n)}$. By translation invariance, we may take $x^{(n)} = y^{(n)} = 0$.

Note that $(\varepsilon_n L^{(n)}_{t/\varepsilon_n^2})_{t\ge 0}$ and $(\varepsilon_n R^{(n)}_{t/\varepsilon_n^2})_{t\ge 0}$ individually converge weakly to a Brownian motion with drift $-1$, respectively, $+1$. This implies tightness for the family of joint processes $\{(L^{(n)}, R^{(n)})\}_{n\in\mathbb{N}}$. Our strategy is to represent $(L^{(n)}_t, R^{(n)}_t)_{t\ge 0}$ as the solution of a difference equation, which in the limit yields an SDE with a unique solution. Since the discontinuous coefficients of the SDE (1.11) are problematic, we prefer to work with (2.2), which behaves better under limits.

Let $(V^{\mathrm{l}}_t)_{t\in\mathbb{N}_0}$, $(V^{\mathrm{r}}_t)_{t\in\mathbb{N}_0}$ and $(V^{\mathrm{s}}_t)_{t\in\mathbb{N}_0}$ be independent discrete-time simple symmetric random walks starting at the origin at time zero. For $\alpha = l, r, s$, let $(D^{(n),\alpha,-}_t)_{t\in\mathbb{N}_0}$ be a process such that whenever $V^\alpha_t$ jumps one step to the right, $D^{(n),\alpha,-}_t$ with probability $\beta_n$ jumps two steps to the left. Likewise, let $(D^{(n),\alpha,+}_t)_{t\in\mathbb{N}_0}$ be the process that with probability $\beta_n$ jumps two steps to the right whenever $V^\alpha_t$ jumps one step to the left. As a result, $V^\alpha_t + D^{(n),\alpha,-}_t$ is a random walk with drift $-\beta_n$, and $V^\alpha_t + D^{(n),\alpha,+}_t$ is a random walk with drift $+\beta_n$.

The unscaled process $(L^{(n)}_t, R^{(n)}_t)$ at integer times can be constructed as the solution of

$$
\begin{aligned}
L^{(n)}_t &= V^{\mathrm{l}}_{T^{(n)}_t} + D^{(n),l,-}_{T^{(n)}_t} + V^{\mathrm{s}}_{S^{(n)}_t} + D^{(n),s,-}_{S^{(n)}_t},\\
R^{(n)}_t &= V^{\mathrm{r}}_{T^{(n)}_t} + D^{(n),r,+}_{T^{(n)}_t} + V^{\mathrm{s}}_{S^{(n)}_t} + D^{(n),s,+}_{S^{(n)}_t},\\
T^{(n)}_t &= \sum_{s=0}^{t-1} 1_{\{L^{(n)}_s < R^{(n)}_s\}},\\
S^{(n)}_t &= \sum_{s=0}^{t-1} 1_{\{L^{(n)}_s = R^{(n)}_s\}}
\end{aligned}
\tag{5.3}
$$

[compare with (2.11)]. We define $L^{(n)}_t, R^{(n)}_t, V^\alpha_t, D^{(n),\alpha,\pm}_t, T^{(n)}_t$ and $S^{(n)}_t$ at noninteger times by linear interpolation. Note that $dT^{(n)}_t = 1_{\{L^{(n)}_{\lfloor t\rfloor} < R^{(n)}_{\lfloor t\rfloor}\}}\,dt$. The rescaled process then satisfies [compare with (2.2)] the following equations:

$$
\begin{aligned}
&\text{(i)} && \varepsilon_n L^{(n)}_{t/\varepsilon_n^2} = \varepsilon_n (V^{\mathrm{l}} + D^{(n),l,-})_{T^{(n)}_{t/\varepsilon_n^2}} + \varepsilon_n (V^{\mathrm{s}} + D^{(n),s,-})_{S^{(n)}_{t/\varepsilon_n^2}},\\
&\text{(ii)} && \varepsilon_n R^{(n)}_{t/\varepsilon_n^2} = \varepsilon_n (V^{\mathrm{r}} + D^{(n),r,-})_{T^{(n)}_{t/\varepsilon_n^2}} + \varepsilon_n (V^{\mathrm{s}} + D^{(n),s,-})_{S^{(n)}_{t/\varepsilon_n^2}},\\
&\text{(iii)} && \varepsilon_n^2 (T^{(n)} + S^{(n)})_{t/\varepsilon_n^2} = t,
\end{aligned}
\tag{5.4}
$$

(iv) $\displaystyle\int_0^t 1_{\{\varepsilon_n R^{(n)}_{s/\varepsilon_n^2} - \varepsilon_n L^{(n)}_{s/\varepsilon_n^2} > \varepsilon_n\}}\, d(\varepsilon_n^2 S^{(n)}_{s/\varepsilon_n^2}) = 0,$

where in the indicator event in (iv), we impose the lower bound of $\varepsilon_n$ instead of 0 for $\varepsilon_n R^{(n)}_{s/\varepsilon_n^2} - \varepsilon_n L^{(n)}_{s/\varepsilon_n^2}$ to compensate the effect of linearly interpolating $S^{(n)}$ between integer times.

Clearly,

$$(\varepsilon_n V^{\mathrm{l}}_{t/\varepsilon_n^2}, \varepsilon_n V^{\mathrm{r}}_{t/\varepsilon_n^2}, \varepsilon_n V^{\mathrm{s}}_{t/\varepsilon_n^2}, -\varepsilon_n D^{(n),l,-}_{t/\varepsilon_n^2},$$
$$\varepsilon_n D^{(n),r,+}_{t/\varepsilon_n^2}, -\varepsilon_n D^{(n),s,-}_{t/\varepsilon_n^2}, \varepsilon_n D^{(n),s,+}_{t/\varepsilon_n^2})_{t\geq 0} \tag{5.5}$$

converge weakly in law on $\mathcal{C}_{\mathbb{R}^7}[0,\infty)$ to

$$(\tilde B^{\mathrm{l}}_t, \tilde B^{\mathrm{r}}_t, \tilde B^{\mathrm{s}}_t, t, t, t, t)_{t\geq 0}. \tag{5.6}$$

We have noted that the laws of $\{(\varepsilon_n L^{(n)}_{t/\varepsilon_n^2}, \varepsilon_n R^{(n)}_{t/\varepsilon_n^2})_{t\geq 0}\}_{n\in\mathbb{N}}$ are tight. Since $t \mapsto \varepsilon_n^2 T^{(n)}_{t/\varepsilon_n^2}$ increases with slope at most 1, the laws of $\{(\varepsilon_n^2 T^{(n)}_{t/\varepsilon_n^2})_{t\geq 0}\}_{n\in\mathbb{N}}$ are also tight. The same is true for $\{(\varepsilon_n^2 S^{(n)}_{t/\varepsilon_n^2})_{t\geq 0}\}_{n\in\mathbb{N}}$. Therefore, for $n\in\mathbb{N}$, the laws of the 11-tuple, which consists of the 7-tuple in (5.5) joint with $(\varepsilon_n L^{(n)}_{t/\varepsilon_n^2}, \varepsilon_n R^{(n)}_{t/\varepsilon_n^2}, \varepsilon_n^2 T^{(n)}_{t/\varepsilon_n^2}, \varepsilon_n^2 S^{(n)}_{t/\varepsilon_n^2})_{t\geq 0}$, are also tight. By going to a subsequence, we may assume that the 11-tuple converges weakly to some limiting process

$$(\tilde B^{\mathrm{l}}_t, \tilde B^{\mathrm{r}}_t, \tilde B^{\mathrm{s}}_t, t, t, t, t, L_t, R_t, T_t, S_t)_{t\geq 0}. \tag{5.7}$$

By Skorohod's representation theorem (see, e.g., Theorem 6.7 in [4]), we can couple the 11-tuples for $n\in\mathbb{N}$ and the limiting process in (5.7), such that the convergence is almost sure in $\mathcal{C}_{\mathbb{R}^{11}}[0,\infty)$.

Assuming this coupling, we claim that $(L_t, R_t, T_t, S_t)_{t\geq 0}$ solves the equation (2.2), and is therefore determined uniquely in law by Lemma 2.2. Indeed, (2.2)(i)–(iii) follow immediately by taking the limit $n\to\infty$ in (5.4)(i)–(iii). We claim that (2.2)(iv) follows from (5.4)(iv). For each $\delta > 0$, choose a continuous nondecreasing function $\rho_\delta : [0,\infty) \to \mathbb{R}$, such that $\rho_\delta(u) = 0$ for $u \leq \delta$ and $\rho_\delta(u) = 1$ for $u \geq 2\delta$. Then, using (5.4)(iv) and taking the limit $n\to\infty$, we find that

$$\int_0^t \rho_\delta(R_s - L_s)\, dS_s = 0 \tag{5.8}$$

for each $\delta > 0$. Letting $\delta \downarrow 0$, we arrive at (2.2)(iv). $\square$

5.2. *Convergence to the left-right Brownian web.* In this section we prove the convergence, under diffusive scaling, of the collections of all left-most and right-most paths in the arrow configuration $\aleph_\beta$ (and their dual) to the left-right Brownian web (and its dual). As a corollary, we also prove formula (1.22).

Recall the scaling map $S_\varepsilon$ defined in (1.7).

PROPOSITION 5.2 (Convergence of multiple left-right paths). *Let $\beta_n, \varepsilon_n \to 0$ with $\beta_n/\varepsilon_n \to 1$. Let $z_1^{(n)}, \dots, z_k^{(n)}, z_1'^{(n)}, \dots, z_{k'}'^{(n)} \in \mathbb{Z}^2_{\rm even}$ be such that $S_{\varepsilon_n}(z_i^{(n)}) \to z_i$ and $S_{\varepsilon_n}(z_j'^{(n)}) \to z_j'$ for $i = 1, \dots, k$ and $j = 1, \dots, k'$. Let $l_i^{(n)}$ denote the path in $\mathcal{U}^{\rm l}_{\beta_n}$ starting from $z_i$, and let $r_j^{(n)}$ denote the path in $\mathcal{U}^{\rm r}_{\beta_n}$ starting from $z_j'$. Then on the space $\Pi^{k+k'}$,*

$$\mathcal{L}(S_{\varepsilon_n}(l_1^{(n)}, \dots, l_k^{(n)}, r_1^{(n)}, \dots, r_{k'}^{(n)})) \underset{n\to\infty}{\Longrightarrow} \mathcal{L}(l_1, \dots, l_k, r_1, \dots, r_{k'}), \tag{5.9}$$

*where $(l_1, \dots, l_k, r_1, \dots, r_{k'})$ is a collection of left-right coalescing Brownian motions as defined in Section 2.2, starting from $(z_1, \dots, z_k, z_1', \dots, z_{k'}')$.*

PROOF. Recall the inductive construction of $(l_1, \dots, l_k, r_1, \dots, r_{k'})$ from Section 2.2. Note that $(l_1^{(n)}, \dots, l_k^{(n)}, r_1^{(n)}, \dots, r_{k'}^{(n)})$ can be constructed using the same inductive approach. Since the inductive construction pieces together independent evolutions of sets of paths, where each set consists of either a single left-most or right-most path or a pair of left-right paths, the proposition follows easily from Proposition 5.1 and the observation that the stopping times used in the inductive construction are almost surely continuous functionals on $\Pi^{k+k'}$ with respect to the law of independent evolutions of paths in different partition elements. □

Let $\hat{\aleph}_\beta$ denote the arrow configuration dual to $\aleph_\beta$, defined exactly as in (1.17), and let $\hat{\mathcal{U}}_\beta$ denote the set of all $\hat{\aleph}_\beta$-paths. Let $\hat{\mathcal{U}}^{\rm l}_\beta$ (resp. $\hat{\mathcal{U}}^{\rm r}_\beta$) denote the set of $\hat{\aleph}_\beta$-paths dual to $\mathcal{U}^{\rm l}_\beta$ (resp. $\mathcal{U}^{\rm r}_\beta$), that is, the set of all left-most (resp. right-most) paths in $\hat{\mathcal{U}}_\beta$ after rotating the graph of $\hat{\mathcal{U}}_\beta$ by 180°. Let $\hat{\Pi} := \{-\pi : \pi \in \Pi\}$, the image space of $\Pi$ under the rotation map $-$, while preserving the metric. We have the following result:

THEOREM 5.3 (Convergence to the left-right Brownian web and its dual). *Let $\beta_n, \varepsilon_n \to 0$ with $\beta_n/\varepsilon_n \to 1$. Then $S_{\varepsilon_n}(\mathcal{U}^{\rm l}_{\beta_n}, \mathcal{U}^{\rm r}_{\beta_n}, \hat{\mathcal{U}}^{\rm l}_{\beta_n}, \hat{\mathcal{U}}^{\rm r}_{\beta_n})$ are $\mathcal{K}(\Pi)^2 \times \mathcal{K}(\hat{\Pi})^2$-valued random variables, and*

$$\mathcal{L}(S_{\varepsilon_n}(\mathcal{U}^{\rm l}_{\beta_n}, \mathcal{U}^{\rm r}_{\beta_n}, \hat{\mathcal{U}}^{\rm l}_{\beta_n}, \hat{\mathcal{U}}^{\rm r}_{\beta_n})) \underset{n\to\infty}{\Longrightarrow} (\mathcal{W}^{\rm l}, \mathcal{W}^{\rm r}, \hat{\mathcal{W}}^{\rm l}, \hat{\mathcal{W}}^{\rm r}), \tag{5.10}$$

*where $(\mathcal{W}^{\rm l}, \mathcal{W}^{\rm r}, \hat{\mathcal{W}}^{\rm l}, \hat{\mathcal{W}}^{\rm r})$ is the left-right Brownian web and its dual.*

PROOF. It follows from Theorem 6.1 of [8], Theorem 1.2 and Proposition 3.2 that

$$\mathcal{L}(S_{\varepsilon_n}(\mathcal{U}^{\mathrm{l}}_{\beta_n},\hat{\mathcal{U}}^{\mathrm{l}}_{\beta_n})) \underset{n\to\infty}{\Longrightarrow} \mathcal{L}(\mathcal{W}^{\mathrm{l}},\hat{\mathcal{W}}^{\mathrm{l}}) \quad \text{and}$$
$$\mathcal{L}(S_{\varepsilon_n}(\mathcal{U}^{\mathrm{r}}_{\beta_n},\hat{\mathcal{U}}^{\mathrm{r}}_{\beta_n})) \underset{n\to\infty}{\Longrightarrow} \mathcal{L}(\mathcal{W}^{\mathrm{r}},\hat{\mathcal{W}}^{\mathrm{r}}). \tag{5.11}$$

Therefore, $\{S_{\varepsilon_n}(\mathcal{U}^{\mathrm{l}}_{\beta_n},\mathcal{U}^{\mathrm{r}}_{\beta_n},\hat{\mathcal{U}}^{\mathrm{l}}_{\beta_n},\hat{\mathcal{U}}^{\mathrm{r}}_{\beta_n})\}_{n\in\mathbb{N}}$ is a tight family. Let $(X^{\mathrm{l}},X^{\mathrm{r}},\hat{X}^{\mathrm{l}},\hat{X}^{\mathrm{r}})$ be any weak limit point. Then $(X^{\mathrm{l}},\hat{X}^{\mathrm{l}})$ and $(X^{\mathrm{r}},\hat{X}^{\mathrm{r}})$ are distributed as $(\mathcal{W}^{\mathrm{l}},\hat{\mathcal{W}}^{\mathrm{l}})$ and $(\mathcal{W}^{\mathrm{r}},\hat{\mathcal{W}}^{\mathrm{r}})$ respectively. Therefore, $(X^{\mathrm{l}},X^{\mathrm{r}})$ satisfies conditions (i) and (iii) of Theorem 1.5. By Proposition 5.2, $(X^{\mathrm{l}},X^{\mathrm{r}})$ also satisfies condition (ii) of Theorem 1.5, and therefore, $(X^{\mathrm{l}},X^{\mathrm{r}})$ has the same distribution as the standard left-right Brownian web $(\mathcal{W}^{\mathrm{l}},\mathcal{W}^{\mathrm{r}})$. Since $\mathcal{W}^{\mathrm{l}}$ and $\mathcal{W}^{\mathrm{r}}$ determine their duals $\hat{\mathcal{W}}^{\mathrm{l}}$ and $\hat{\mathcal{W}}^{\mathrm{r}}$ almost surely, $(X^{\mathrm{l}},X^{\mathrm{r}},\hat{X}^{\mathrm{l}},\hat{X}^{\mathrm{r}})$ has the same distribution as $(\mathcal{W}^{\mathrm{l}},\mathcal{W}^{\mathrm{r}},\hat{\mathcal{W}}^{\mathrm{l}},\hat{\mathcal{W}}^{\mathrm{r}})$. □

PROOF OF FORMULA (1.22). Since the analogue of (1.22) obviously holds in the discrete setting, (1.22) is a consequence of the convergence in (5.10). □

5.3. *Convergence to the Brownian net.* In this section we prove Theorem 1.1. It suffices to prove (1.8) for $b=1$ and $b=0$. The general case $b>0$ follows the same proof as for $b=1$ if we set $\mathcal{L}(\mathcal{N}_b):=\mathcal{L}(S_{1/b}(\mathcal{N}))$, which automatically gives the scaling relation (1.9). Thus, Theorem 1.1 is implied by the following stronger result.

THEOREM 5.4 (Convergence to the associated Brownian net). *Let $\beta_n,\varepsilon_n\to 0$ with $\beta_n/\varepsilon_n\to b\in\{0,1\}$. Then $S_{\varepsilon_n}(\mathcal{U}_{\beta_n},\mathcal{U}^{\mathrm{l}}_{\beta_n},\mathcal{U}^{\mathrm{r}}_{\beta_n},\hat{\mathcal{U}}^{\mathrm{l}}_{\beta_n},\hat{\mathcal{U}}^{\mathrm{r}}_{\beta_n})$ are $\mathcal{K}(\Pi)^3\times\mathcal{K}(\hat{\Pi})^2$-valued random variables. If $b=1$, then*

$$\mathcal{L}(S_{\varepsilon_n}(\mathcal{U}_{\beta_n},\mathcal{U}^{\mathrm{l}}_{\beta_n},\mathcal{U}^{\mathrm{r}}_{\beta_n},\hat{\mathcal{U}}^{\mathrm{l}}_{\beta_n},\hat{\mathcal{U}}^{\mathrm{r}}_{\beta_n})) \underset{n\to\infty}{\Longrightarrow} \mathcal{L}(\mathcal{N},\mathcal{W}^{\mathrm{l}},\mathcal{W}^{\mathrm{r}},\hat{\mathcal{W}}^{\mathrm{l}},\hat{\mathcal{W}}^{\mathrm{r}}), \tag{5.12}$$

*where $\mathcal{N}$ is the (standard) Brownian net and $(\mathcal{W}^{\mathrm{l}},\mathcal{W}^{\mathrm{r}},\hat{\mathcal{W}}^{\mathrm{l}},\hat{\mathcal{W}}^{\mathrm{r}})$ is its associated left-right Brownian web and its dual. If $b=0$, then*

$$\mathcal{L}(S_{\varepsilon_n}(\mathcal{U}_{\beta_n},\mathcal{U}^{\mathrm{l}}_{\beta_n},\mathcal{U}^{\mathrm{r}}_{\beta_n},\hat{\mathcal{U}}^{\mathrm{l}}_{\beta_n},\hat{\mathcal{U}}^{\mathrm{r}}_{\beta_n})) \underset{n\to\infty}{\Longrightarrow} (\mathcal{W},\mathcal{W},\mathcal{W},\hat{\mathcal{W}},\hat{\mathcal{W}}), \tag{5.13}$$

*where $(\mathcal{W},\hat{\mathcal{W}})$ is the Brownian web and its dual.*

PROOF. We start with the case $b=1$ and then say how our arguments can be adapted to cover also the case $b=0$.

Recall the modulus of continuity $m_K(\cdot)$ of $K\in\mathcal{K}(\Pi)$ from (4.3). Just as in the proof of Lemma 4.6, we see that

$$m_{S_{\varepsilon_n}(\mathcal{U}_{\beta_n})}(\delta)\le m_{S_{\varepsilon_n}(\mathcal{U}^{\mathrm{l}}_{\beta_n}\cup\mathcal{U}^{\mathrm{r}}_{\beta_n})}(\delta), \tag{5.14}$$

hence, the tightness of $\{S_{\varepsilon_n}(\mathcal{U}_{\beta_n})\}_{n\in\mathbb{N}}$ follows from the tightness of $S_{\varepsilon_n}(\mathcal{U}^{\rm l}_{\beta_n})$ and $S_{\varepsilon_n}(\mathcal{U}^{\rm r}_{\beta_n})$ $(n\in\mathbb{N})$. Thus, by going to a subsequence, we may assume that the laws in (5.12) converge to a limit $\mathcal{L}(\mathcal{N}^*,\mathcal{W}^{\rm l},\mathcal{W}^{\rm r},\hat{\mathcal{W}}^{\rm l},\hat{\mathcal{W}}^{\rm r})$, where by Theorem 5.3, $(\mathcal{W}^{\rm l},\mathcal{W}^{\rm r},\hat{\mathcal{W}}^{\rm l},\hat{\mathcal{W}}^{\rm r})$ is the left-right Brownian web and its dual. We need to show that $\mathcal{N}^*$ is the Brownian net associated with $(\mathcal{W}^{\rm l},\mathcal{W}^{\rm r},\hat{\mathcal{W}}^{\rm l},\hat{\mathcal{W}}^{\rm r})$. Our strategy will be to show that $\mathcal{N}_{\rm hop}\subset\mathcal{N}^*\subset\mathcal{N}_{\rm wedge}$, where $\mathcal{N}_{\rm hop}$ and $\mathcal{N}_{\rm wedge}$ are defined as in Section 4. It then follows from the equivalence of the hopping and dual constructions of the Brownian net (Theorems 1.3 and 1.10) that $\mathcal{N}^*=\mathcal{N}$.

Let $\mathcal{D}^{\rm l},\mathcal{D}^{\rm r}\subset\mathbb{R}^2$ be deterministic countable dense sets. For each $z\in\mathcal{D}^{\rm l}$ (resp. $z'\in\mathcal{D}^{\rm r}$), we fix a sequence $z_n\in\mathbb{Z}^2_{\rm even}$ (resp. $z'_n\in\mathbb{Z}^2_{\rm even}$) such that $S_{\varepsilon_n}(z_n)\to z$ (resp. $S_{\varepsilon_n}(z'_n)\to z'$), and we let $\hat{l}^{(n)}_z$ (resp. $\hat{r}^{(n)}_{z'}$) denote the path in $S_{\varepsilon_n}(\hat{\mathcal{U}}^{\rm l}_{\beta_n})$ (resp. $S_{\varepsilon_n}(\hat{\mathcal{U}}^{\rm r}_{\beta_n})$) starting in $S_{\varepsilon_n}(z_n)$ (resp. $S_{\varepsilon_n}(z'_n)$). Let

$$\tau(\hat{\pi}_1,\hat{\pi}_2):=\sup\{t<\hat{\sigma}_{\hat{\pi}_1}\wedge\hat{\sigma}_{\hat{\pi}_2}:\hat{\pi}_1(t)=\hat{\pi}_2(t)\} \tag{5.15}$$

denote the first meeting time of the two dual paths $\hat{\pi}_1,\hat{\pi}_2$. Since, up to their first meeting time, $\hat{l}^{(n)}_z$ and $\hat{r}^{(n)}_{z'}$ are independent random walks, and since random walk paths joint with their first meeting time converge under diffusive scaling to Brownian motions joint with their first meeting time, we have

$$\begin{aligned}&\mathcal{L}(S_{\varepsilon_n}(\mathcal{U}_{\beta_n},\mathcal{U}^{\rm l}_{\beta_n},\mathcal{U}^{\rm r}_{\beta_n},\hat{\mathcal{U}}^{\rm l}_{\beta_n},\hat{\mathcal{U}}^{\rm r}_{\beta_n}),(\tau(\hat{l}^{(n)}_z,\hat{r}^{(n)}_{z'}))_{z\in\mathcal{D}^{\rm l},\ z'\in\mathcal{D}^{\rm r}})\\&\quad\underset{n\to\infty}{\Longrightarrow}\mathcal{L}(\mathcal{N}^*,\mathcal{W}^{\rm l},\mathcal{W}^{\rm r},\hat{\mathcal{W}}^{\rm l},\hat{\mathcal{W}}^{\rm r},(\tau(\hat{l}_z,\hat{r}_{z'}))_{z\in\mathcal{D}^{\rm l},\ z'\in\mathcal{D}^{\rm r}}).\end{aligned} \tag{5.16}$$

By Skorohod's representation theorem, we can construct a coupling such that the convergence in (5.16) is almost sure. Assuming such a coupling, we will show that $\mathcal{N}_{\rm hop}\subset\mathcal{N}^*\subset\mathcal{N}_{\rm wedge}$.

To show that $\mathcal{N}_{\rm hop}\subset\mathcal{N}^*$, it suffices to show that $\mathcal{H}_{\rm cros}(\mathcal{W}^{\rm l}(\mathcal{D}^{\rm l})\cup\mathcal{W}^{\rm r}(\mathcal{D}^{\rm r}))$ is contained in $\mathcal{N}^*$. Any $\pi\in\mathcal{H}_{\rm cros}(\mathcal{W}^{\rm l}(\mathcal{D}^{\rm l})\cup\mathcal{W}^{\rm r}(\mathcal{D}^{\rm r}))$ is constructed by hopping at crossing times between left-most and right-most skeletal paths $\pi_1,\ldots,\pi_m$ as in (1.12). By the a.s. convergence of $S_{\varepsilon_n}(\mathcal{U}^{\rm l}_{\beta_n},\mathcal{U}^{\rm r}_{\beta_n})$ to $(\mathcal{W}^{\rm l},\mathcal{W}^{\rm r})$, there exist $\pi^{(n)}_i\in S_{\varepsilon_n}(\mathcal{U}^{\rm l}_{\beta_n}\cup\mathcal{U}^{\rm r}_{\beta_n})$ such that $\pi^{(n)}_i\to\pi_i$ $(i=1,\ldots,m)$. By the structure of crossing times [Proposition 3.6(a)], the crossing time between $\pi^{(n)}_i$ and $\pi^{(n)}_{i+1}$ converges to the crossing time between $\pi_i$ and $\pi_{i+1}$ for all $i=1,\ldots,m-1$. Therefore, the path $\pi^{(n)}$ that is constructed by hopping at crossing times between $\pi^{(n)}_1,\ldots,\pi^{(n)}_m$ converges to $\pi$. Since $\pi^{(n)}\in S_{\varepsilon_n}(\mathcal{U}_{\beta_n})$ by the nearest-neighbor nature of $\aleph_{\beta_n}$-paths, this proves that $\mathcal{H}_{\rm cros}(\mathcal{W}^{\rm l}(\mathcal{D}^{\rm l})\cup\mathcal{W}^{\rm r}(\mathcal{D}^{\rm r}))\subset\mathcal{N}^*$.

To show that $\mathcal{N}^*\subset\mathcal{N}_{\rm wedge}$, we need to show that a.s. no path $\pi\in\mathcal{N}^*$ enters a wedge $W(\hat{r},\hat{l})$ from outside. If $\pi\in\mathcal{N}^*$ enters a wedge $W(\hat{r},\hat{l})$ from outside, then by Lemma 3.4(b), $\pi$ must enter some skeletal wedge $W(\hat{r}_{z'}, \hat{l}_z)$, with $z \in \mathcal{D}^{\mathrm{l}}$ and $z' \in \mathcal{D}^{\mathrm{r}}$, from outside. By the a.s. convergence of $S_{\varepsilon_n}(\mathcal{U}_{\beta_n})$ to $\mathcal{N}^*$, there exist $\pi^{(n)} \in S_{\varepsilon_n}(\mathcal{U}_{\beta_n})$ such that $\pi^{(n)} \to \pi$. By the a.s. convergence of $\hat{r}^{(n)}_{z'}$ and $\hat{l}^{(n)}_z$ to $\hat{r}_{z'}$ and $\hat{l}_z$ and the convergence of their first meeting time, for $n$ large enough, $\pi^{(n)}$ must enter a discrete wedge from outside, which is impossible.

This concludes the proof for $b = 1$. The proof for $b = 0$ is similar. Note that if in the left-right SDE (1.11), one removes the drift terms $\pm dt$, then solutions $(L, R)$ are just coalescing Brownian motions. Using this fact, it is not hard to generalize Propositions 5.1 and 5.2 in the sense that if $\beta_n/\varepsilon_n \to 0$, then left-most and right-most paths converge to coalescing Brownian motions (with zero drift). Modifying Theorem 5.3 appropriately, we find that

$$\mathcal{L}(S_{\varepsilon_n}(\mathcal{U}^{\mathrm{l}}_{\beta_n}, \mathcal{U}^{\mathrm{r}}_{\beta_n}, \hat{\mathcal{U}}^{\mathrm{l}}_{\beta_n}, \hat{\mathcal{U}}^{\mathrm{r}}_{\beta_n})) \underset{n\to\infty}{\Longrightarrow} (\mathcal{W}, \mathcal{W}, \hat{\mathcal{W}}, \hat{\mathcal{W}}). \tag{5.17}$$

By going to a subsequence if necessary, we may assume that $S_{\varepsilon_n}(\mathcal{U}_{\beta_n})$ converges to some limit $\mathcal{W}^*$. The inclusion $\mathcal{W} \subset \mathcal{W}^*$ is now trivial, while the other inclusion can be obtained by showing that no path in $\mathcal{W}^*$ enters a wedge of $\hat{\mathcal{W}}$ from outside, applying Theorem 1.9. □

**6. Density calculations.** In this section we carry out two density calculations for the Brownian net $\mathcal{N}$, based on the hopping and dual characterizations (Theorem 1.3 and Theorem 1.10), which have been shown in Section 4 to be equivalent. In Section 6.1 we calculate the density of the set of points on $\mathbb{R} \times \{t\}$ that are on the graph of some path in $\mathcal{N}$ starting at time 0, that is, we prove Proposition 1.12. In Section 6.2 we estimate the density of the set of times that are the first meeting times between $l \in \mathcal{W}^{\mathrm{l}}(0,0)$ and some path in $\mathcal{N}_{\mathrm{hop}}$ starting to the left of 0 at time 0. Our calculations show that both sets are a.s. locally finite. The second density calculation gives information on the configuration of meshes on the left of a general left-most path $l$, which will be used in Section 7 to prove that paths in $\mathcal{N}_{\mathrm{mesh}}$ cannot enter the area to the left of $l$. From this, we then readily obtain Theorem 1.7, as well as Propositions 1.4, 1.8 and 1.13.

6.1. *The density of the branching-coalescing point set.* In this section we prove Proposition 1.12. Let $\mathcal{N}$ be the Brownian net, defined by the hopping or dual characterization (Theorem 1.3 and Theorem 1.10). Set

$$\xi_t := \{\pi(t) : \pi \in \mathcal{N}, \sigma_\pi = 0\} \qquad (t > 0). \tag{6.1}$$

Note that $\xi_t = \xi_t^{\overline{\mathbb{R}} \times \{0\}}$, the branching-coalescing point set (defined in Section 1.9) started at time zero from $\overline{\mathbb{R}}$. The exact computation of the density of $\xi_t$ is based on the following two lemmas.

LEMMA 6.1 (Avoidance of intervals). *Almost surely, for each* $s,t,a,b\in\mathbb{R}$ *with* $s<t$ *and* $a<b$, *there exists no* $\pi\in\mathcal{N}(\overline{\mathbb{R}}\times\{s\})$ *with* $\pi(t)\in(a,b)$ *if and only if there exist* $\hat{r}\in\hat{\mathcal{W}}^{\mathrm{r}}(a,t)$ *and* $\hat{l}\in\hat{\mathcal{W}}^{\mathrm{l}}(b,t)$ *such that* $\sup\{u<t:\hat{r}(u)=\hat{l}(u)\}>s$.

PROOF. If $\hat{r},\hat{l}$ with the described properties exist, then by the dual characterization of the Brownian net (Theorem 1.10), no path in $\mathcal{N}$ starting at time $s$ can pass through $(a,b)\times\{t\}$. Conversely, if there exists no $\pi\in\mathcal{N}(\overline{\mathbb{R}}\times\{s\})$ such that $\pi(t)\in(a,b)$, then for each $\varepsilon>0$ and for each $\hat{r}_\varepsilon\in\hat{\mathcal{W}}^{\mathrm{r}}(a+\varepsilon,t)$ and $\hat{l}_\varepsilon\in\hat{\mathcal{W}}^{\mathrm{l}}(b-\varepsilon,t)$, we must have $\tau_\varepsilon:=\sup\{u<t:\hat{r}_\varepsilon(u)=\hat{l}_\varepsilon(u)\}>s$. For if $\tau_\varepsilon\le s$, then by the steering argument used in the proof of Lemma 4.7 (see Figure 7), for each $\delta>0$ we can construct a path in $\mathcal{H}_{\mathrm{cros}}(\mathcal{W}^{\mathrm{l}}\cup\mathcal{W}^{\mathrm{r}})$ starting at time $s+\delta$ in $(\hat{r}_\varepsilon(s+\delta),\hat{l}_\varepsilon(s+\delta))$ and passing through $[a+\varepsilon,b-\varepsilon]\times t$. Letting $\hat{r},\hat{l}$ denote any limits of paths $\hat{r}_{\varepsilon_n},\hat{l}_{\varepsilon'_n}$ along sequences $\varepsilon_n,\varepsilon'_n\downarrow 0$, we see that $\tau:=\sup\{u<t:\hat{r}(u)=\hat{l}(u)\}>s$. In fact, by Lemma 3.4(a), we must have $\tau>s$. $\square$

LEMMA 6.2 (Hitting probability of a pair of left-right SDE). *Let* $L_s$ *and* $R_s$ *be the solution of* (1.11) *with initial condition* $L_0=0$ *and* $R_0=\varepsilon$ *for some* $\varepsilon>0$. *Let* $T_\varepsilon=\inf\{s\ge 0:L_s=R_s\}$. *Then*

$$1-\Psi_\varepsilon(t):=\mathbb{P}[T_\varepsilon<t]=\Phi\left(-\sqrt{2t}-\frac{\varepsilon}{\sqrt{2t}}\right)+e^{-2\varepsilon}\Phi\left(\sqrt{2t}-\frac{\varepsilon}{\sqrt{2t}}\right), \tag{6.2}$$

*where* $\Phi(x)=\int_{-\infty}^{x}\frac{e^{-y^2/2}}{\sqrt{2\pi}}\,dy$.

PROOF. Let $Y_t=B_t+\sqrt{2}t$ with $Y_0=0$, and let $M_t=-\inf_{0\le s\le t}Y_s$. Clearly, $R_t-L_t-\varepsilon$ is equally distributed with $\sqrt{2}Y_t$ before it reaches level $-\varepsilon$. Therefore, $\mathbb{P}[T_\varepsilon<t]=\mathbb{P}[M_t\ge\varepsilon/\sqrt{2}]$. We compute this last probability by first finding the joint density of $B'_t$, a standard Brownian motion, and $M'_t=-\inf_{0\le s\le t}B'_s$. We then apply Girsanov's formula to change the measure from $(B'_s)_{0\le s\le t}$ to that of $(Y_s)_{0\le s\le t}$.

For a standard Brownian motion $B'_t$, it is easy to check by reflection principle that, for $x\ge 0$ and $y\ge -x$,

$$\mathbb{P}[M'_t\ge x,B'_t\ge y]=\mathbb{P}[B'_t\ge 2x+y]=\int_{2x+y}^{\infty}\frac{e^{-z^2/2}}{\sqrt{2\pi}}\,dz. \tag{6.3}$$

Differentiating with respect to $x$ and $y$ gives the joint density

$$\begin{aligned}&\mathbb{P}[M'_t\in dx,B'_t\in dy]\\&\quad=\frac{1}{\sqrt{2\pi t}}\cdot\frac{2(2x+y)}{t}\cdot e^{-(2x+y)^2/(2t)}\,dx\,dy\qquad x\ge 0,y\ge -x.\end{aligned} \tag{6.4}$$

By Girsanov's formula, the measure for $(Y_s)_{0\le s\le t}$ is absolute continuous with respect to the measure for $(B'_s)_{0\le s\le t}$ with density $e^{\sqrt{2}B'_t-t}$. Therefore,

$$
\begin{aligned}
&\mathbb{P}\Big[M_t \ge \frac{\varepsilon}{\sqrt{2}}\Big] \\
&\quad = \int_{\varepsilon/\sqrt{2}}^{\infty}\int_{-x}^{\infty} e^{\sqrt{2}y-t}\frac{1}{\sqrt{2\pi t}}\cdot\frac{2(2x+y)}{t}\cdot e^{-(2x+y)^2/(2t)}\,dy\,dx.
\end{aligned}
\tag{6.5}
$$

Split the integral into two regions: $\mathrm{I}=\int_{-\varepsilon/\sqrt{2}}^{\infty}dy\int_{\varepsilon/\sqrt{2}}^{\infty}dx$; and $\mathrm{II}=\int_{-\infty}^{-\varepsilon/\sqrt{2}}dy\times\int_{-y}^{\infty}dx$. Then we have

$$
\begin{aligned}
\mathrm{I} &= e^{-t}\int_{-\varepsilon/\sqrt{2}}^{\infty}\frac{e^{\sqrt{2}y}}{\sqrt{2\pi t}}\,dy\int_{\varepsilon/\sqrt{2}}^{\infty}\frac{2(2x+y)}{t}\cdot e^{-(2x+y)^2/(2t)}\,dx \\
&= e^{-t}\int_{-\varepsilon/\sqrt{2}}^{\infty}\frac{1}{\sqrt{2\pi t}}\,e^{\sqrt{2}y-(y+\sqrt{2}\varepsilon)^2/(2t)}\,dy \\
&= e^{-2\varepsilon}\int_{-\varepsilon/\sqrt{2}}^{\infty}\frac{1}{\sqrt{2\pi t}}\,e^{-(y+\sqrt{2}\varepsilon-\sqrt{2}t)^2/(2t)}\,dy = e^{-2\varepsilon}\Phi\Big(\sqrt{2t}-\frac{\varepsilon}{\sqrt{2t}}\Big).
\end{aligned}
\tag{6.6}
$$

Similarly,

$$
\begin{aligned}
\mathrm{II} &= e^{-t}\int_{-\infty}^{-\varepsilon/\sqrt{2}}\frac{1}{\sqrt{2\pi t}}\,e^{\sqrt{2}y-y^2/(2t)}\,dy \\
&= \int_{-\infty}^{-\varepsilon/\sqrt{2}}\frac{1}{\sqrt{2\pi t}}\,e^{-(y-\sqrt{2}t)^2/(2t)}\,dy = \Phi\Big(-\sqrt{2t}-\frac{\varepsilon}{\sqrt{2t}}\Big).
\end{aligned}
\tag{6.7}
$$

This concludes the proof. □

PROOF OF PROPOSITION 1.12. It follows from Lemmas 6.1 and 6.2, and the continuity of $\varepsilon\mapsto\Psi_\varepsilon(t)$ that

$$
\mathbb{P}[\xi_t\cap(a,b)\neq\varnothing]=\mathbb{P}[\xi_t\cap[a,b]\neq\varnothing]=\Psi_{b-a}(t)\qquad(t>0)
\tag{6.8}
$$

for deterministic $a<b$. Since the law of $\xi_t$ is clearly translation invariant in space, to prove (1.28), without loss of generality, we may assume $[a,b]=[0,1]$. Let $\mathcal{R}=\{\frac{i}{2^n}:n\in\mathbb{N},0\le i\le 2^n\}$ denote the dyadic rationals. By (6.8), $\mathbb{P}[x\in\xi_t]=0$ for each deterministic $x\in\mathbb{R}$. Since $\mathcal{R}$ is countable, almost surely $\xi_t\cap\mathcal{R}=\varnothing$. Therefore,

$$
|\xi_t\cap[0,1]|=\lim_{n\to\infty}\Big|\Big\{1\le i\le 2^n:\xi_t\cap\Big[\frac{i-1}{2^n},\frac{i}{2^n}\Big]\neq\varnothing\Big\}\Big|\qquad\text{a.s.}
\tag{6.9}
$$

By monotone convergence and translation invariance,

$$
\mathbb{E}[|\xi_t\cap[0,1]|]=\lim_{n\to\infty}2^n\mathbb{P}\Big[\xi_t\cap\Big[0,\frac{1}{2^n}\Big]\neq\varnothing\Big]=\frac{\partial}{\partial\varepsilon}\Psi_\varepsilon(t)|_{\varepsilon=0},
\tag{6.10}
$$

which yields equation (1.28). □

6.2. *The density on the left of a left-most path.* Let $\mathcal{N}$ be the Brownian net, defined by the hopping or dual characterization (Theorem 1.3 and Theorem 1.10), and let $(\mathcal{W}^{\mathrm{l}}, \mathcal{W}^{\mathrm{r}}, \hat{\mathcal{W}}^{\mathrm{l}}, \hat{\mathcal{W}}^{\mathrm{r}})$ be its associated left-right Brownian web and its dual. For each $l \in \mathcal{W}^{\mathrm{l}}$, let

$$C(l) := \{t > \sigma_l : \exists \pi \in \mathcal{N} \text{ s.t. } \sigma_\pi = \sigma_l, \qquad \pi(t) = l(t), \pi(s) < l(s)\ \forall s \in [\sigma_l, t)\} \tag{6.11}$$

be the set of times when some path in $\mathcal{N}$, started at the same time as $l$ and to the left of $l$, first meets $l$. We will prove that, almost surely, $C(l)$ is a locally finite subset of $(\sigma_l, \infty)$ for each $l \in \mathcal{W}^{\mathrm{l}}$. By Lemma 3.4(b), it suffices to verify this property for $l \in \mathcal{W}^{\mathrm{l}}$ with deterministic starting points, in particular, $l$ started at $(0,0)$, which is implied by the following lemma.

PROPOSITION 6.3 (Density on the left of a left-most path). *Let $l$ be the a.s. unique path in $\mathcal{W}^{\mathrm{l}}$ starting at the origin. Then, for each $0 < s < t$,*

$$\mathbb{E}[|C(l) \cap [s,t]|] \le \int_s^t 2\psi(u)^2\, du, \tag{6.12}$$

*where $\psi(t) := \frac{\partial}{\partial\varepsilon}\Psi_\varepsilon(t)|_{\varepsilon=0} = \frac{e^{-t}}{\sqrt{\pi t}} + 2\Phi(\sqrt{2t})$ is the density of the branching-coalescing point set in* (1.28).

PROOF. By a similar argument as in the proof of Proposition 1.12, it suffices to show that

$$\limsup_{\varepsilon \to 0} \frac{1}{\varepsilon}\mathbb{P}[C(l) \cap [t, t+\varepsilon] \ne \varnothing] \le 2\psi(t)^2. \tag{6.13}$$

For $t > 0$, let $\hat{r}_{[t]}$ be the left-most [viewed with respect to the graph of $(\mathcal{W}^{\mathrm{r}}, \hat{\mathcal{W}}^{\mathrm{r}})$] path in $\hat{\mathcal{W}}^{\mathrm{r}}(l(t), t)$ and let $\hat{l}_{[t]}$ be the right-most path in $\hat{\mathcal{W}}^{\mathrm{l}}(l(t), t)$ *that lies on the left of $l$*. Note that, by Lemma 3.3(b), for each deterministic $t > 0$, $\hat{\mathcal{W}}^{\mathrm{l}}(l(t), t)$ almost surely contains two paths, one lying on each side of $l$. Similar arguments as in the proof of Lemma 6.1 show that

$$\begin{aligned} &\mathbb{P}[C(l) \cap [t, t+\varepsilon] \ne \varnothing] \\ &\qquad = \mathbb{P}[\hat{r}_{[t+\varepsilon]}(s) < \hat{l}_{[t]}(s)\ \forall s \in (0,t)]. \end{aligned} \tag{6.14}$$

Set

$$\begin{aligned} L_s &:= l(t+\varepsilon) - l(t-s), && s \in [-\varepsilon, t], \\ \hat{L}_s &:= l(t+\varepsilon) - \hat{l}_{[t]}(t-s), && s \in [0, t], \\ \hat{R}_s &:= l(t+\varepsilon) - \hat{r}_{[t+\varepsilon]}(t-s), && s \in [-\varepsilon, t]. \end{aligned} \tag{6.15}$$

It has been shown in [15] (see also [9]) that paths in $\mathcal{W}$ and $\hat{\mathcal{W}}$ interact by Skorohod reflection. Similar arguments show that if a path $\hat{r} \in \hat{\mathcal{W}}^{\mathrm{r}}$ is started on the left of a path $l \in \mathcal{W}^{\mathrm{l}}$, then $\hat{r}$ is Skorohod reflected off $l$. Therefore, on the time interval $[-\varepsilon, 0]$, the process $(L_s, \hat{R}_s)$ satisfies $L \leq \hat{R}$ and solves the SDE

$$\begin{aligned} dL_s &= dB^{\mathrm{l}}_s - ds, \\ d\hat{R}_s &= dB^{\hat{\mathrm{r}}}_s + ds + d\Delta'_s, \end{aligned} \tag{6.16}$$

where $B^{\mathrm{l}}$ and $B^{\hat{\mathrm{r}}}$ are independent Brownian motions, and $\Delta'_s$ is a reflection term that increases only when $L_s = \hat{R}_s$. Set $\sigma := \inf\{s > 0 : \hat{L}_s = \hat{R}_s\} \wedge t$. Then on the time interval $[0, \sigma]$, the process $(L_s, \hat{L}_s, \hat{R}_s)$ satisfies $L \leq \hat{L} \leq \hat{R}$ and solves the SDE

$$\begin{aligned} dL_s &= dB^{\mathrm{l}}_s - ds, \\ d\hat{L}_s &= dB^{\hat{\mathrm{l}}}_s - ds + d\Delta_s, \\ d\hat{R}_s &= dB^{\hat{\mathrm{r}}}_s + ds, \end{aligned} \tag{6.17}$$

where $B^{\mathrm{l}}, B^{\hat{\mathrm{l}}}, B^{\hat{\mathrm{r}}}$ are independent Brownian motions and $\Delta_s$ increases only when $L_s = \hat{L}_s$. By Lemma 6.4 below,

$$\mathbb{P}[\hat{L}_s < \hat{R}_s \ \forall s \in (0, t)] \leq \int \mathbb{P}[\hat{R}_0 - L_0 \in d\eta] \Psi_\eta(t)^2. \tag{6.18}$$

Set $X_s := \hat{R}_{s-\varepsilon} - L_{s-\varepsilon}$ $(s \in [0, \varepsilon])$. Then $X$ is a Brownian motion with diffusion constant 2 and drift 2, Skorohod reflected at 0, which has the generator $\frac{\partial^2}{\partial \eta^2} + 2\frac{\partial}{\partial \eta}$ with boundary condition $\frac{\partial}{\partial \eta} f(\eta)|_{\eta=0} = 0$. Therefore,

$$\begin{aligned} &\lim_{\varepsilon \to 0} \varepsilon^{-1} \int \mathbb{P}[\hat{R}_0 - L_0 \in d\eta] \Psi_\eta(t)^2 \\ &\quad = \lim_{\varepsilon \to 0} \varepsilon^{-1} \mathbb{E}[\Psi_{X_\varepsilon}(t)^2] \\ &\quad = \left(\frac{\partial^2}{\partial \eta^2} + 2\frac{\partial}{\partial \eta}\right)(\Psi_\eta(t)^2)\Big|_{\eta=0} = 2\psi(t)^2, \end{aligned} \tag{6.19}$$

where we have used that, for fixed $t > 0$, $\eta \mapsto \Psi_\eta(t)^2$ is a bounded twice continuously differentiable function satisfying our boundary condition. □

LEMMA 6.4 (Hitting estimate). *Let $(L, \hat{L}, \hat{R})$ be a solution to the SDE (6.17) started at $(L_0, \hat{L}_0, \hat{R}_0) = (0, 0, \eta)$. Then*

$$\mathbb{P}[\hat{L}_s < \hat{R}_s \ \forall s \in (0, t)] \leq \Psi_\eta(t)^2, \tag{6.20}$$

*where $\Psi_\eta(t)$ is defined in (6.2).*

PROOF. We introduce new coordinates:

$$V_t := \hat{L}_t - L_t, \qquad W_t := \hat{R}_t - L_t. \tag{6.21}$$

The process $(V,W)$ lives in the space $\{(v,w) \in \mathbb{R}^2 : 0 \le v \le w\}$ up to the time $\tau := \inf\{t > 0 : V_t = W_t\}$ and solves the SDE

$$dV_t := dB^{\hat{\mathrm{l}}}_s - dB^{\mathrm{l}}_s + d\Delta_s, \qquad dW_t := dB^{\hat{\mathrm{r}}}_s - dB^{\mathrm{l}}_s + 2\,ds, \tag{6.22}$$

where $\Delta_s$ is a reflection term, increasing only when $V_s = 0$. Changing coordinates once more, we set

$$X_t := W_t - V_t, \qquad Y_t := W_t + V_t. \tag{6.23}$$

Then $(X,Y)$ takes values in $\{(x,y) \in \mathbb{R}^2 : 0 \le x \le y\}$ up to the time $\tau := \inf\{t > 0 : X_t = 0\}$ and solves the SDE

$$dX_s := dB^{\hat{\mathrm{r}}}_s - dB^{\hat{\mathrm{l}}}_s + 2\,ds - d\Delta_s, \qquad dY_s := dB^{\hat{\mathrm{r}}}_s + dB^{\hat{\mathrm{l}}}_s - 2\,dB^{\mathrm{l}}_s + 2\,ds + d\Delta_s, \tag{6.24}$$

where $\Delta_s$ increases only when $X_s = Y_s$. Our strategy will be to compare $(X,Y)$ with a process $(X',Y')$ of the form $X' = U^1 \wedge U^2$ and $Y' = U^1 \vee U^2$, where $U^1, U^2$ are independent processes with generator $\frac{\partial^2}{\partial u^2} + 2\frac{\partial}{\partial u}$. We will show that $X$ hits zero before $X'$. Note that if $U^i_0 = u$, then $\mathbb{P}[U^i_s > 0\ \forall s \in [0,t]] = \Psi_u(t)$, which is defined in (6.2). Therefore,

$$\frac{\partial}{\partial t}\Psi_u(t) = \left(\frac{\partial^2}{\partial u^2} + 2\frac{\partial}{\partial u}\right)\Psi_u(t). \tag{6.25}$$

Moreover, if $(X',Y')$ is started in $(x,y)$, then $\mathbb{P}[X'_s > 0\ \forall s \in [0,t]] = \mathbb{P}[U^1_s > 0\ \forall s \in [0,t]]\mathbb{P}[U^2_s > 0\ \forall s \in [0,t]] = \Psi_x(t)\Psi_y(t)$. With this in mind, we set

$$F(t,x,y) := \Psi_x(t)\Psi_y(t). \tag{6.26}$$

Let $G$ be the operator

$$G := \frac{\partial^2}{\partial x^2} + 2\frac{\partial}{\partial x} + 3\frac{\partial^2}{\partial y^2} + 2\frac{\partial}{\partial y}. \tag{6.27}$$

By Itô's formula,

$$\begin{aligned} &dF(t-s, X_{s\wedge\tau}, Y_{s\wedge\tau}) \\ &\quad = \left(-\frac{\partial}{\partial t} + 1_{\{s<\tau\}}G\right)F(t-s, X_{s\wedge\tau}, Y_{s\wedge\tau})\,ds \\ &\qquad + 1_{\{s<\tau\}}\left(\frac{\partial}{\partial y} - \frac{\partial}{\partial x}\right)F(t-s, X_{s\wedge\tau}, Y_{s\wedge\tau})\,d\Delta_s \end{aligned} \tag{6.28}$$

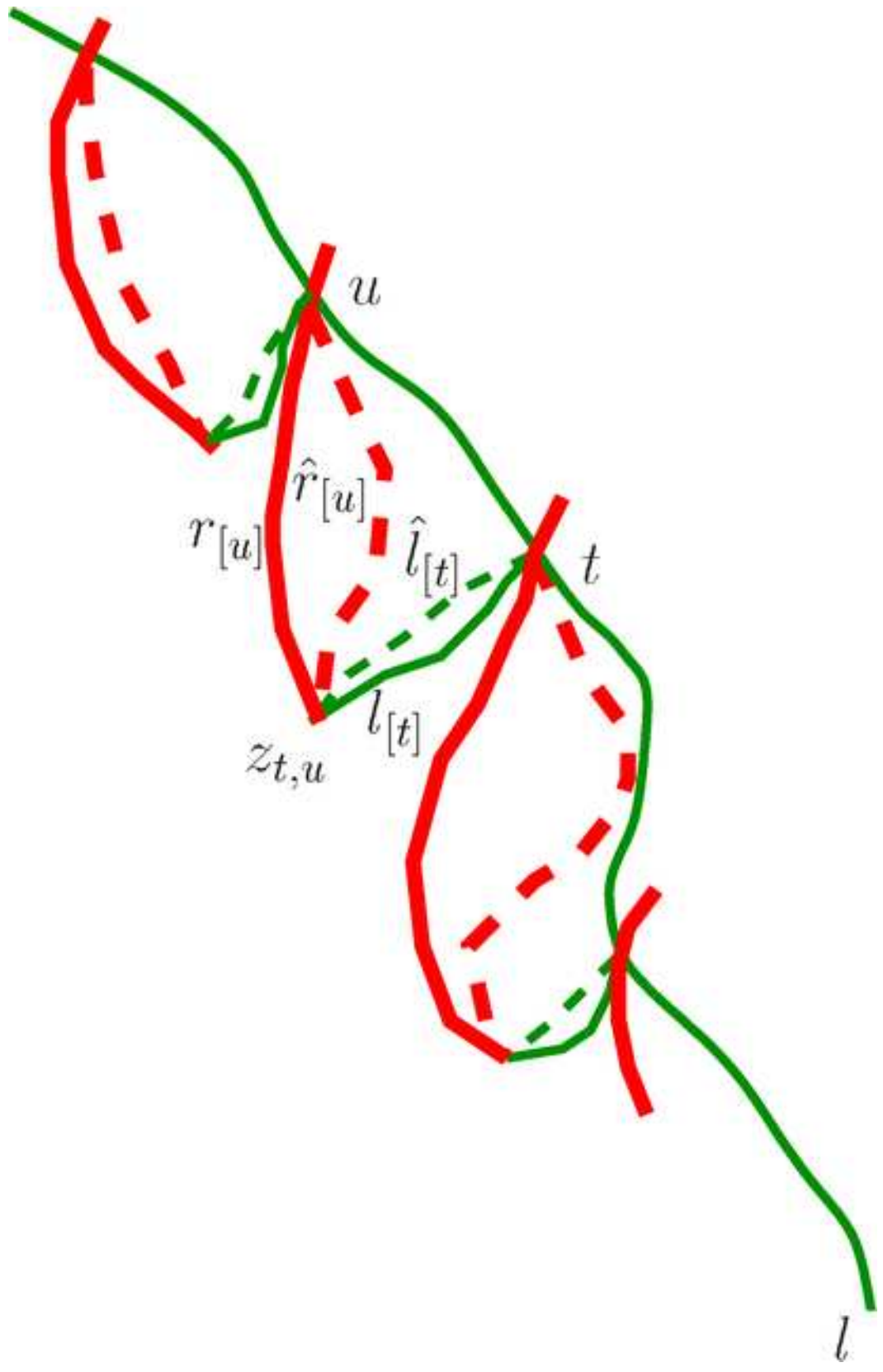

FIG. 8. *Meshes stack up on the left of a left-most path $l \in \mathcal{W}^{\mathrm{l}}$.*

plus martingale terms. It follows from the definition of $\Psi_u(t)$ that $\frac{\partial}{\partial t}\Psi_u(t) \leq 0$ and $\frac{\partial}{\partial u}\Psi_u(t) \geq 0$, and therefore, by (6.25), $\frac{\partial^2}{\partial u^2}\Psi_u(t) \leq 0$. As a result, using (6.25) and (6.26), we see that $(\frac{\partial}{\partial y} - \frac{\partial}{\partial x})F(t,x,y)|_{x=y} = 0$ and

$$\left(-\frac{\partial}{\partial t} + G\right)F(t,x,y) = 2\frac{\partial^2}{\partial y^2}(\Psi_x(t)\Psi_y(t)) \leq 0. \tag{6.29}$$

Inserting this into (6.28), we find that $(F(t-s, X_{s\wedge\tau}, Y_{s\wedge\tau}))_{s\in[0,t\wedge\tau]}$ is a local supermartingale, which implies that

$$\mathbb{P}[\tau > t] = \mathbb{E}[F(t - t\wedge\tau, X_{t\wedge\tau}, Y_{t\wedge\tau})] \leq F(t, X_0, Y_0) = \Psi_\eta(t)^2. \tag{6.30}$$ □

As a corollary to Proposition 6.3, we obtain the following lemma, which describes the configuration of meshes on the left of a left-most path. (See Figure 8.)

LEMMA 6.5 (Meshes on the left of a left-most path). *Almost surely, the set* $C(l)$ *in (6.11) is a locally finite subset of* $(\sigma_l,\infty)$ *for each* $l\in\mathcal{W}^{\mathrm{l}}$. *For each consecutive pair of times* $t,u\in C(l)$ *[i.e.,* $t<u$ *and* $C(l)\cap(t,u)=\varnothing$*], there exists a mesh* $M(r',l')$ *with bottom time* $s\in(\sigma_l,t)$ *and top point* $(l(u),u)$*, such that* $l'<l$ *on* $[s,t)$ *and* $l'=l$ *on* $[t,u]$*. If* $C(l)$ *has a minimal element* $t$*, then there exists a mesh* $M(r',l)$ *with right boundary* $l$*, bottom point* $(l(\sigma_l),\sigma_l)$ *and top point* $(l(t),t)$*.*

PROOF. For any path $\pi$ and $\varepsilon>0$, define a trunctated path by $\pi^{\langle\varepsilon\rangle}:=\{(\pi(t),t):t\in[\sigma_\pi+\varepsilon,\infty]\}$. Let $l_{(0,0)}$ be the a.s. unique left-most path starting in the origin. The proof of Proposition 6.3 applies to $l^{\langle\varepsilon\rangle}_{(0,0)}$ as well; in particular, $C(l^{\langle\varepsilon\rangle}_{(0,0)})$ has the same density as $C(l_{(0,0)})$ for each $\varepsilon>0$. By Lemma 3.4(b), if follows that a.s., $C(l^{\langle\varepsilon\rangle})$ is a locally finite subset of $(\sigma_l+\varepsilon,\infty)$ for each $l\in\mathcal{W}^{\mathrm{l}}$ and $\varepsilon>0$. Since $C(l^{\langle\varepsilon\rangle})\cap(\sigma_l+\delta,\infty)$ decreases to $C(l)\cap(\sigma_l+\delta,\infty)$ as $\varepsilon\downarrow 0$, for each fixed $\delta>0$, it follows that a.s., $C(l)$ is a locally finite subset of $(\sigma_l,\infty)$ for each $l\in\mathcal{W}^{\mathrm{l}}$.

For any $l\in\mathcal{W}^{\mathrm{l}}$ (see Figure 8), consider $t,u\in C(l)\cup\{\sigma_l\}$ such that $t<u$ and $C(l)\cap(t,u)=\varnothing$, that is, either $t,u$ is a consecutive pair of times in $C(l)$, or $t=\sigma_l$ and $u$ is the minimal element of $C(l)$. By an argument similar to the proof of Lemma 6.1, there exist $\hat r_{[u]}\in\hat{\mathcal{W}}^{\mathrm{r}}(l(u),u)$ and $\hat l_{[t]}\in\hat{\mathcal{W}}^{\mathrm{l}}(l(t),t)$ such that $\hat r_{[u]}\le l$ on $[\sigma_l,u]$, $\hat l_{[t]}\le l$ on $[\sigma_l,t]$, and $\tau_{t,u}:=\sup\{s\le t:\hat r_{[u]}(s)=\hat l_{[t]}(s)\}$ satisfies $\tau_{t,u}>\sigma_l$ if $\tau_{t,u}<t$. (Note that possibly $\tau_{t,u}=t$.)

Set $z_{t,u}:=(\hat r_{[u]}(\tau_{t,u}),\tau_{t,u})$. Let $r_{[u]}$ denote the left-most path in $\mathcal{W}^{\mathrm{r}}(z_{t,u})$. Let $l_{[t]}$ denote the right-most path in $\mathcal{W}^{\mathrm{l}}(z_{t,u})$ if $\tau_{t,u}<t$, and let $l_{[t]}$ denote the path in $\mathcal{W}^{\mathrm{l}}(z_{t,u})$ that is the continuation of $l$ if $\tau_{t,u}=t$. Set $u':=\inf\{s>\tau_{t,u}:r_{[u]}(s)=l(s)\}$ and $t':=\inf\{s>\tau_{t,u}:l_{[t]}(s)=l(s)\}$. By Proposition 3.2(c) and (e), $r_{[u]}\le\hat r_{[u]}$ on $[\tau_{t,u},u]$, and therefore, by Propositions 3.6(a), (b), $u'\ge u$. Likewise, since $l_{[t]}\ge\hat l_{[t]}$ on $[\tau_{t,u},t]$, we have $t'\le t$. Now $r_{[u]}$ and $l_{[t]}$ are the left and right boundary of a mesh $M(r_{[u]},l_{[t]})$ with bottom time $\tau_{t,u}$ and top point $(l(u'),u')$, such that $l_{[t]}<l$ on $(\tau_{t,u},t')$ and $l_{[t]}=l$ on $[t',u']$. Since $\mathcal{N}_{\mathrm{hop}}\subset\mathcal{N}_{\mathrm{mesh}}$ (Lemma 4.4) and both $t$ (if $t\sigma_l$) and $u$ are times when a path in $\mathcal{N}_{\mathrm{hop}}$ starting at time $\sigma_l$ first meets $l$ from the left, it follows that $t'=t$ and $u'=u$. (If $t=\sigma_l$, then obviously $\tau_{t,u}=\sigma_l=t=t'$.) To complete the proof, we must show that $\tau_{t,u}<t$ if $t>\sigma_l$. This follows from Lemma 6.6 below. □

LEMMA 6.6 (Top and bottom points of meshes). *Almost surely, no bottom point of one mesh is the top point of another mesh.*

PROOF. Assume that $z\in\mathbb{R}^2$ is the bottom point of a mesh $M(r,l)$ and the top point of another mesh $M(r',l')$. By Propositions 3.2(c) and 3.6(d), any $\hat{r} \in \hat{\mathcal{W}}^{\mathrm{r}}$ starting in $M(r,l)$ must pass through $z$ (and likewise for $\hat{l} \in \mathcal{W}^{\mathrm{l}}$). Therefore, $l', r'$ and $\hat{r}$ are three paths entering $z$ disjointly. This can be ruled out just as in the proof of Theorem 3.11 in [9], where it is argued that a.s. there is no point $z \in \mathbb{R}^2$ where two forward and one backward path in $(\mathcal{W}, \hat{\mathcal{W}})$ enter $z$ disjointly. □

**7. Characterization with meshes.** In this section, we prove Theorem 1.7, as well as Propositions 1.4, 1.8 and 1.13. We fix a left-right Brownian web and its dual $(\mathcal{W}^{\mathrm{l}}, \mathcal{W}^{\mathrm{r}}, \hat{\mathcal{W}}^{\mathrm{l}}, \hat{\mathcal{W}}^{\mathrm{r}})$ and define $\mathcal{N}_{\mathrm{hop}}, \mathcal{N}_{\mathrm{wedge}}$ and $\mathcal{N}_{\mathrm{mesh}}$ as in Section 4. The key technical result is the following lemma, which states that Proposition 1.8 holds for $\mathcal{N}_{\mathrm{mesh}}$.

LEMMA 7.1 (Containment by left-most and right-most paths). *Almost surely, there exist no $\pi \in \mathcal{N}_{\mathrm{mesh}}$ and $l \in \mathcal{W}^{\mathrm{l}}$ such that $l(s) \leq \pi(s)$ and $\pi(t) < l(t)$ for some $\sigma_\pi \vee \sigma_l < s < t$. An analogue statement holds for right-most paths.*

PROOF. Without loss of generality, we may assume that $\sigma_l > \sigma_\pi$; otherwise, consider a left-most path starting at any time in $(\sigma_\pi, s)$ that is the continuation of $l$. By Lemma 6.5, there exists a locally finite collection of meshes on the left of $l$, with bottom times in $[\sigma_l, \infty)$, that block the way of any path in $\mathcal{N}_{\mathrm{mesh}}$ trying to enter the area to the left of $l$. (See Figure 8.) □

PROOF OF THEOREM 1.7 AND PROPOSITION 1.8. We start by proving that $\mathcal{N}_{\mathrm{mesh}} \subset \mathcal{N}_{\mathrm{wedge}}$. Since, by Lemma 4.3, paths in $\mathcal{N}_{\mathrm{mesh}}$ do not enter wedges through their bottom points, it suffices to show that paths in $\mathcal{N}_{\mathrm{mesh}}$ do not cross dual left-most and right-most paths in the wrong direction. By Lemma 3.5, it suffices to show that paths in $\mathcal{N}_{\mathrm{mesh}}$ do not cross forward left-most and right-most paths in the wrong direction. This follows from Lemma 7.1.

Since it has already been proved in Lemmas 4.4, 4.5 and 4.7 that $\mathcal{N}_{\mathrm{mesh}} \supset \mathcal{N}_{\mathrm{hop}} = \mathcal{N}_{\mathrm{wedge}}$, it follows that all these sets are a.s. equal. This proves Theorem 1.7. Lemma 7.1 now translates into Proposition 1.8. □

PROOF OF PROPOSITION 1.4. By Theorem 1.7 and Proposition 1.8, the Brownian net $\mathcal{N}$ associated with a left-right Brownian web $(\mathcal{W}^{\mathrm{l}}, \mathcal{W}^{\mathrm{r}})$ consists exactly of those paths in $\Pi$ that do not enter meshes. It is easy to see that this set is closed under hopping. □

PROOF OF PROPOSITION 1.13. Let $(\mathcal{W}^{\mathrm{l}}, \mathcal{W}^{\mathrm{r}})$ be the left-right Brownian web associated with $\mathcal{N}$. We have to show that, for each $t \in [-\infty, \infty]$ and $\pi \in$ $\Pi_t$ such that $\pi \subset \cup(\mathcal{N} \cap \Pi_t)$, we have $\pi \in \mathcal{N}$. By Theorem 1.7, each mesh of $(\mathcal{W}^{\rm l}, \mathcal{W}^{\rm r})$ with bottom time in $(t, \infty)$ has empty intersection with $\cup(\mathcal{N} \cap \Pi_t)$, and therefore, $\pi$ does not enter any such mesh. Again by Theorem 1.7, it follows that $\pi \in \mathcal{N}$. $\square$

**8. The branching-coalescing point set.** In this section we prove Theorem 1.11. We start with two preparatory lemmas.

LEMMA 8.1 (Hopping paths starting from a closed set). *Let $\mathcal{N}$ be the Brownian net. Let $K \subset \mathbb{R}^2$ be closed. Let $\mathcal{D}^{\rm l}, \mathcal{D}^{\rm r} \subset \mathbb{R}^2$ be deterministic countable dense sets such that, moreover, $\mathcal{D}^{\rm l} \cap K$ is dense in $K$. Then*

$$\mathcal{N}(K) \subset \overline{\Pi(K) \cap \mathcal{H}_{\rm cros}(\mathcal{W}^{\rm l}(\mathcal{D}^{\rm l}) \cup \mathcal{W}^{\rm r}(\mathcal{D}^{\rm r}))}. \tag{8.1}$$

PROOF. By the dual characterization of the Brownian net (Theorem 1.10), it suffices to show that any path $\pi$ starting at some point $z = (x, t) \in K$ that does not enter wedges from outside can be approximated by paths in $\mathcal{H}_{\rm cros}(\mathcal{W}^{\rm l}(\mathcal{D}^{\rm l}) \cup \mathcal{W}^{\rm r}(\mathcal{D}^{\rm r}))$, also starting in $K$. By the compactness of $\mathcal{N}$, it suffices to show that, for each $t < t_1 < \cdots < t_n$, and $0 < \varepsilon < t_1 - t$, there exists a path $\pi^\varepsilon$ started at some time in $(t - \varepsilon, t + \varepsilon)$, such that $|\pi^\varepsilon(t_i) - \pi(t_i)| \le \varepsilon$ for all $i = 1, \ldots, n$. We use the steering argument from the proof of Lemma 4.7 (see Figure 7). We construct a "fish-trap" with left and right boundary $\hat{r}, \hat{l}$ as in Figure 7. Set $S := \sup\{s < t_n : \hat{r}(s) = \hat{l}(s)\}$. For any $a, b \in \mathbb{R}$ with $a < b \le \hat{\sigma}_{\hat{r}} \wedge \hat{\sigma}_{\hat{l}}$, define an open set $V_{(a,b)}$ by

$$V_{(a,b)} := \{(x, s) \in \mathbb{R}^2 : a < s < b,\ S < s, \hat{r}(s) < x < \hat{l}(s)\}. \tag{8.2}$$

We need to show that there exists a path $\pi^\varepsilon$, started at some time in $(t - \varepsilon, t + \varepsilon)$, that stays between $\hat{r}$ and $\hat{l}$. This will follow from the same arguments as in the proof of Lemma 4.7, provided that the set $V_{(t-\varepsilon, t+\varepsilon)} \cap \mathcal{D}^{\rm l}$ is nonempty. Since $\mathcal{D}^{\rm l} \cap K$ is dense in $K$ and $V_{(t-\varepsilon, t+\varepsilon)}$ is open, it suffices to show that $V_{(t-\varepsilon, t+\varepsilon)} \cap K$ is nonempty. Assume that this is not the case. Then, by Lemma 3.4(b), we can find $a < b$ and $\hat{r} \in \hat{\mathcal{W}}^{\rm r}(\mathcal{D}^{\rm r})$, $\hat{l} \in \hat{\mathcal{W}}^{\rm r}(\mathcal{D}^{\rm l})$ starting at times $\hat{\sigma}_{\hat{r}}, \hat{\sigma}_{\hat{l}} > b$, such that $V_{(a,b)} \cap K = \varnothing$ but $z \in \partial V_{(a,b)}$, where we define

$$\partial V_{(a,b)} := \{(x, s) \in \mathbb{R}^2 : a < s < b,\ S \le s,\ x = \hat{r}(s) \text{ or } x = \hat{l}(s)\}. \tag{8.3}$$

We claim that this is impossible. More precisely, we claim that if $K \subset \mathbb{R}^2$ is a deterministic closed set and $\hat{r} \in \hat{\mathcal{W}}^{\rm r}$, $\hat{l} \in \hat{\mathcal{W}}^{\rm r}$ are paths with deterministic starting points, then almost surely, there exist no $a, b \in \mathbb{R}$, with $a < b < \hat{\sigma}_{\hat{r}}, \hat{\sigma}_{\hat{l}}$, such that $V_{(a,b)} \cap K = \varnothing$ and $\partial V_{(a,b)} \cap K \ne \varnothing$. It suffices to prove the statement for deterministic $a, b$. Set

$$Y := \hat{r}(b) + \hat{l}(b),$$

$$Z := \hat{r}(b) - \hat{l}(b),$$
$$(8.4) \qquad R := (\hat{r}(s) - \hat{r}(b))_{s\in[a,b]},$$
$$L := (\hat{l}(s) - \hat{l}(b))_{s\in[a,b]}.$$

We claim that, for any $y \in \mathbb{R}$ and continuous functions $\omega_{\mathrm{r}}, \omega_{\mathrm{l}} : [a, b] \to \mathbb{R}$, the conditional probability

$$(8.5) \qquad \mathbb{P}[V_{(a,b)} \cap K = \varnothing \text{ and } \partial V_{(a,b)} \cap K \neq \varnothing | Y = y,\ R = \omega_{\mathrm{r}},\ L = \omega_{\mathrm{l}}]$$

is zero. Indeed, for given $Y, R$ and $L$, there can be at most one value of $Z$ for which the event $V_{(a,b)} \cap K = \varnothing$ and $\partial V_{(a,b)} \cap K \neq \varnothing$ occurs. Since conditioned on $S < b$, which is necessary for $\partial V_{(a,b)} \neq \varnothing$, the distribution of the random variable $Z$ is absolute continuous with respect to Lebesgue measure, the conditional probability in (8.5) is zero. Integrating over the distributions of $Y, R$ and $L$, we arrive at our result. □

LEMMA 8.2 (Almost sure continuity). *Let $\mathcal{N}$ be the Brownian net, and let $K_n, K \in \mathcal{K}(R^2_{\mathrm{c}})$ be deterministic sets satisfying $K_n \to K$. Assume that $(*, -\infty) \notin K$. Then $\mathcal{N}(K_n) \to \mathcal{N}(K)$ a.s.*

PROOF. Using the compactness of $\mathcal{N}$, by going to a subsequence if necessary, we may assume that $\mathcal{N}(K_n) \to \mathcal{A}$ for some compact subset $\mathcal{A} \subset \mathcal{N}$. Obviously, all paths in $\mathcal{A}$ have starting points in $K$, so $\mathcal{A} \subset \mathcal{N}(K)$. Write $K' := K \cap \mathbb{R}^2$ and $K'' := K \backslash K'$. Since $\mathcal{N}(z)$ is trivial for $z \in K''$, it is easy to see that $\mathcal{A} \supset \mathcal{N}(K'')$. We are left with the task to show $\mathcal{A} \supset \mathcal{N}(K')$. Choose a deterministic countable dense set $\mathcal{D} \subset \mathbb{R}^2$ such that, moreover, $\mathcal{D} \cap K'$ is dense in $K'$. For each $z \in \mathcal{D} \cap K'$, choose $z_n \in K_n$ such that $z_n \to z$. Then $l_{z_n} \to l_z$. If $l_z$ crosses a path $r \in \mathcal{W}^{\mathrm{r}}$, then for $n$ large enough, $l_{z_n}$ also crosses $r$. Therefore, it is not hard to see that

$$(8.6) \qquad \mathcal{A} \supset \Pi(K') \cap \mathcal{H}_{\mathrm{cros}}(\mathcal{W}^{\mathrm{l}}(\mathcal{D}) \cup \mathcal{W}^{\mathrm{r}}(\mathcal{D})).$$

By Lemma 8.1, it follows that $\mathcal{A} \supset \mathcal{N}(K')$. □

REMARK. If $(*, -\infty) \in K$, then by Proposition 1.15(iii), the conclusion of Lemma 8.2 still holds, provided there exist $(x_n, t_n) \in K_n$ such that $t_n \to -\infty$ and $\limsup_{n\to\infty} |x_n|/|t_n| < 1$. Here some of the $(x_n, t_n)$ may be $(*, -\infty)$, with the convention that $|*|/|\infty| := 0$. As the proof of Proposition 1.15(iii) shows, this condition cannot be relaxed very much.

PROOF OF THEOREM 1.11. The continuity of sample paths of $(\xi_t)_{t\geq 0}$ is a direct consequence of the definition of $\xi_t$ and the fact that $\mathcal{N}$ is a $\mathcal{K}(\Pi)$-valued random variable. The fact that $\xi_t$ is a.s. locally finite in $\mathbb{R}$ for deterministic $t > s$ follows from Proposition 1.12.

For $t \geq 0$, the transition probability kernel $P_t$ on $\mathcal{K}(\overline{\mathbb{R}})$ associated with $\xi$ is given by

$$P_t(K, \cdot) := \mathbb{P}[\xi_{s+t}^{K\times\{s\}} \in \cdot], \qquad K \in \mathcal{K}(\overline{\mathbb{R}}). \tag{8.7}$$

Note that the right-hand side of (8.7) does not depend on $s \in \mathbb{R}$ by the translation invariance of the Brownian net. By Lemma 8.2, if $K_n \to K$ and $t_n \to t$, then

$$P_{t_n}(K_n, \cdot) = \mathbb{P}[\xi_0^{K_n\times\{-t_n\}} \in \cdot] \underset{n\to\infty}{\Longrightarrow} \mathbb{P}[\xi_0^{K\times\{-t\}} \in \cdot] = P_t(K, \cdot), \tag{8.8}$$

proving the Feller property of $(P_t)_{t\geq 0}$. We still have to show that $(P_t)_{t\geq 0}$ is a Markov transition probability kernel. This is not completely obvious, but it follows, provided we show that, for any $s < t_0 < t_1$ and compact $K \subset \overline{\mathbb{R}}$,

$$\mathbb{P}[\xi_{t_1}^{K\times\{s\}} \in \cdot \mid (\xi_u^{K\times\{s\}})_{u\in[s,t_0]}] = P_{t_1-t_0}(\xi_{t_0}^{K\times\{s\}}, \cdot) \qquad \text{a.s.} \tag{8.9}$$

Let $\pi|_s^t := \{(\pi(u), u) : u \in [s,t] \cap [\sigma_\pi, \infty]\}$ denote the restriction of a path $\pi \in \Pi$ to the time interval $[s,t]$, and for $\mathcal{A} \subset \Pi$, write $\mathcal{A}|_s^t := \{\pi|_s^t : \pi \in \mathcal{A}\}$. In view of the definition of $\xi_t$, it suffices to show that

$$\mathbb{P}[\mathcal{N}(K\times\{s\})|_{t_0}^{\infty} \in \cdot \mid \mathcal{N}(K\times\{s\})|_s^{t_0}] = \mathbb{P}[\mathcal{N}'(\xi_{t_0}^{K\times\{s\}} \times \{t_0\}) \in \cdot], \tag{8.10}$$

where $\mathcal{N}'$ is an independent copy of $\mathcal{N}$. Let $(\mathcal{W}^{\mathrm{l}}, \mathcal{W}^{\mathrm{r}})$ be the left-right Brownian web associated with $\mathcal{N}$. By the properties of left-right coalescing Brownian motions, $(\mathcal{W}^{\mathrm{l}}, \mathcal{W}^{\mathrm{r}})|_{-\infty}^{t_0}$ and $(\mathcal{W}^{\mathrm{l}}, \mathcal{W}^{\mathrm{r}})|_{t_0}^{\infty}$ are independent, and therefore, by the hopping construction, it follows that $\mathcal{N}|_{-\infty}^{t_0}$ and $\mathcal{N}|_{t_0}^{\infty}$ are independent. In particular, $\xi_{t_0}^{K\times\{s\}}$ and $\mathcal{N}(K\times\{s\})|_s^{t_0}$ are independent of $\mathcal{N}(\overline{\mathbb{R}}\times\{t_0\})$. To show (8.10), it therefore suffices to show that

$$\mathcal{N}(K\times\{s\})|_{t_0}^{\infty} = \mathcal{N}(\xi_{t_0}^{K\times\{s\}} \times \{t_0\}) \qquad \text{a.s.} \tag{8.11}$$

The inclusion $\subset$ is trivial. To prove the converse, we need to show that any path $\pi \in \mathcal{N}(\xi_{t_0}^{K\times\{s\}} \times \{t_0\})$ is the continuation of a path in $\mathcal{N}(K\times\{s\})$; this follows from Lemma 8.3 below.

To prove (1.27), note that $K \in \mathcal{K}'(\overline{\mathbb{R}})$ if and only if $\sup K < \infty$, or $\sup(K \cap \mathbb{R}) = \infty$ and $\infty \in K$, and likewise at $-\infty$. Thus, by symmetry, it suffices to show that, almost surely:

$$\begin{array}{llll} & \text{(i)} & \sup(\xi_s) < \infty & \text{implies} \quad \sup(\xi_t) < \infty \ \forall t \geq s, \\ (8.12) & \text{(ii)} & \sup(\xi_s \cap \mathbb{R}) = \infty & \text{implies} \quad \sup(\xi_t \cap \mathbb{R}) = \infty \ \forall t \geq s, \\ & \text{(iii)} & \infty \in \xi_s & \text{implies} \quad \infty \in \xi_t \ \forall t \geq s. \end{array}$$

Formula (i) follows from the fact that $(\sup(\xi_t))_{t\geq s}$ is the right-most path in $\mathcal{N}(\xi_s \times \{s\})$, which is a Brownian motion with drift $+1$. Formula (ii) is easily proved by considering the right-most paths starting at a sequence of points in $\xi_s \cap \mathbb{R}$ tending to $(\infty, s)$. Last, formula (iii) follows from the fact that $\mathcal{N}(\infty, s)$ contains the trivial path $\pi(t) := \infty$ $(t \geq s)$. □

In the proof of Theorem 1.11 we have used the following lemma, which is of some interest on its own.

LEMMA 8.3 (Hopping at deterministic times). *Let $\mathcal{N}$ be the Brownian net and $t \in \mathbb{R}$. Then almost surely, for each $\pi, \pi' \in \mathcal{N}$ such that $\sigma_\pi \vee \sigma_{\pi'} \leq t$ and $\pi(t) = \pi'(t)$, the path $\pi''$ defined by*

$$\pi'' := \{(\pi(s), s) : s \in [\sigma_\pi, t]\} \cup \{(\pi'(s), s) : s \in [t, \infty]\} \tag{8.13}$$

*satisfies $\pi'' \in \mathcal{N}$.*

PROOF. If $\sigma_\pi = t$, there is nothing to prove, so without loss of generality we may assume that $\sigma_\pi \leq s$ for some deterministic $s < t$. If $\pi'' \notin \mathcal{N}$, then by the dual characterization of the Brownian net, $\pi''$ must enter a wedge from outside, which can only happen if $(\pi(t), t)$ lies on a dual path. But this is not possible since $\pi(t)$ lies in $\xi_t^{\overline{R} \times \{s\}}$, which is locally finite (by Proposition 1.12) and independent of $(\hat{\mathcal{W}}^{\mathrm{l}}, \hat{\mathcal{W}}^{\mathrm{r}})|_{t_0}^{\infty}$, and a.s. no Brownian web path passes through a deterministic point. □

We end this section with a proposition that will be used in the proof of Lemma 9.2, and that is of interest in its own right. Note that the statement below implies that, provided that the initial states converge, systems of branching-coalescing random walks, diffusively rescaled, converge in an appropriate sense to the branching-coalescing point set.

PROPOSITION 8.4 (Convergence of paths started from subsets). *Let $\beta_n$, $\varepsilon_n \to 0$ with $\beta_n / \varepsilon_n \to 1$. Let $K_n \subset \mathbb{Z}^2_{\mathrm{even}}$, $K \in \mathcal{K}(R^2_{\mathrm{c}})$ satisfy $S_{\varepsilon_n}(K_n) \to K$, where $\to$ denotes convergence in $\mathcal{K}(R^2_{\mathrm{c}})$. Assume $(*, -\infty) \notin K$. Then*

$$\mathcal{L}(S_{\varepsilon_n}(\mathcal{U}_{\beta_n}(K_n))) \underset{n \to \infty}{\Longrightarrow} \mathcal{L}(\mathcal{N}(K)). \tag{8.14}$$

PROOF. By going to a subsequence if necessary, we may assume that $\mathcal{L}(S_{\varepsilon_n}(\mathcal{U}_{\beta_n}, \mathcal{U}_{\beta_n}(K_n))) \Rightarrow \mathcal{L}(\mathcal{N}, \mathcal{A})$ for some compact subset $\mathcal{A} \subset \mathcal{N}(K)$. Write $K' := K \cap \mathbb{R}^2$ and $K'' := K \backslash K'$. Since $\mathcal{N}(z)$ is trivial for $z \in K''$, it is easy to see that $\mathcal{A} \supset \mathcal{N}(K'')$. We are left with the task to show $\mathcal{A} \supset \mathcal{N}(K')$. Choose a deterministic countable dense set $\mathcal{D} \subset R^2_{\mathrm{c}}$ such that, moreover, $\mathcal{D} \cap K'$ is dense in $K'$. By the same arguments as those used in the proof of Theorem 5.4 to show that $\mathcal{N}_{\mathrm{hop}} \subset \mathcal{N}^*$, we have

$$\Pi(K') \cap \mathcal{H}_{\mathrm{cros}}(\mathcal{W}^{\mathrm{l}}(\mathcal{D}) \cup \mathcal{W}^{\mathrm{r}}(\mathcal{D})) \subset \mathcal{A}. \tag{8.15}$$

By Lemma 8.1, it follows that $\mathcal{N}(K') \subset \mathcal{A}$. □

**9. The backbone.** In Sections 9.1 and 9.2 we prove Propositions 1.14 and 1.15, respectively.

9.1. *The backbone of branching-coalescing random walks.* Let $\aleph_\beta$ be an arrow configuration. Recall the definition of $\eta_t^A$ from (1.2). Let $\mathbb{Z}_{\rm even} := 2\mathbb{Z}$ and $\mathbb{Z}_{\rm odd} := 2\mathbb{Z}+1$. For any $s \in \mathbb{Z}$ and $A \subset \mathbb{Z}_{\rm even}$ or $A \subset \mathbb{Z}_{\rm odd}$ depending on whether $s$ is even or odd, setting

$$\eta_t := \eta_t^{A\times\{s\}} \qquad (t \in \mathbb{Z},\ t \geq s) \tag{9.1}$$

defines a Markov chain $(\eta_t)_{t\geq s}$ taking values, in turn, in the spaces of subsets of $\mathbb{Z}_{\rm even}$ and $\mathbb{Z}_{\rm odd}$, started at time $s$ in $A$. We call $\eta = (\eta_t)_{t\geq s}$ a *system of branching-coalescing random walks*. We call a probability law $\mu$ on the space of subsets of $\mathbb{Z}_{\rm even}$ an *invariant law* for $\eta$ if $\mathcal{L}(\eta_0) = \mu$ implies $\mathcal{L}(\eta_2) = \mu$, and a *homogeneous invariant law* if $\mu$ is translation invariant and $\mathcal{L}(\eta_0) = \mu$ implies $\mathcal{L}(\eta_1 + 1) = \mu$. Note that we shift $\eta_1$ by one unit in space to stay on $\mathbb{Z}_{\rm even}$.

It is easy to see that $\mathcal{L}(\eta_0^{(*,-\infty)})$ defines a homogeneous invariant law for $\eta$. Our strategy for proving Proposition 1.14 will be as follows. First we prove that the Bernoulli measure $\mu_\rho$ with intensity $\rho = \frac{4\beta}{(1+\beta)^2}$ is a homogeneous invariant law for $\eta$, and that $\mu_\rho$ is reversible in a sense that includes information about the arrow configuration $\aleph_\beta$. Next, we prove Proposition 1.14(iii). From this, we derive that there exists only one nontrivial invariant law for $\eta$, hence, $\mathcal{L}(\eta_0^{(*,-\infty)}) = \mu_\rho$, which proves part (i). Last, part (ii) follows from the reversibility of $\mu_\rho$.

We first need to add additional structure to the branching-coalescing random walks that also keeps track of the arrows in $\aleph_\beta$ that are used by the walks. To this aim, if $(\eta_t)_{t=s,s+1,\dots}$ is defined as in (9.1) with respect to an arrow configuration $\aleph_\beta$, then we define

$$\eta_{t+1/2} := \{\{x,x'\} : x \in \eta_t, ((x,t),(x',t+1)) \in \aleph_\beta\} \qquad (t \in [s,\infty) \cap \mathbb{Z}). \tag{9.2}$$

Note that $\eta_{t+1/2}$ keeps track of which arrows in $\aleph_\beta$ are used by the branching-coalescing random walks between the times $t$ and $t+1$. It is not hard to see that $(\eta_{s+k/2})_{k\in\mathbb{N}_0}$ is a Markov chain.

LEMMA 9.1 (Product invariant law). *The Bernoulli product measure $\mu_\rho$ on $\mathbb{Z}_{\rm even}$ with intensity $\rho = \frac{4\beta}{(1+\beta)^2}$ is a reversible homogeneous invariant law for the Markov chain $(\eta_{s+k/2})_{k\in\mathbb{N}_0}$ defined above, in the sense that, if $\mathcal{L}(\eta_0) = \mu_\rho$, then for all even $t \geq 0$,*

$$\mathcal{L}(\eta_0, \eta_{1/2}, \dots, \eta_{t-1/2}, \eta_t) = \mathcal{L}(\eta_t, \eta_{t-1/2}, \dots, \eta_{1/2}, \eta_0). \tag{9.3}$$

*The same holds for all odd $t \geq 1$, provided that the configurations on the right-hand side of (9.3) are shifted in space by one unit.*

PROOF. It suffices to prove the statement for $t = 1$, that is, we need to prove that if $\mathcal{L}(\eta_0) = \mu_\rho$, then

$$\mathcal{L}(\eta_0, \eta_{1/2}, \eta_1) = \mathcal{L}(\eta_1 + 1, \eta_{1/2} + 1, \eta_0 + 1). \tag{9.4}$$

Indeed, since $(\eta_{t/2})_{t\in\mathbb{N}_0}$ is Markov, $(\eta_0, \ldots, \eta_{s-1/2})$ and $(\eta_{s+1/2}, \ldots, \eta_t)$ are conditionally independent given $\eta_s$ for all $s \in [1, t] \cap \mathbb{Z}$. The identity (9.3) for general even $t \geq 0$, and its analogue for odd $t \geq 0$, then follow easily from (9.4) by induction.

Note that $\eta_{1/2}$ determines $\eta_0$ and $\eta_1$ a.s. Indeed,

$$\begin{aligned} \eta_0 &= \{x \in \mathbb{Z}_{\text{even}} : \exists x' \in \mathbb{Z}_{\text{odd}} \text{ s.t. } \{x, x'\} \in \eta_{1/2}\}, \\ \eta_1 &= \{x' \in \mathbb{Z}_{\text{odd}} : \exists x \in \mathbb{Z}_{\text{even}} \text{ s.t. } \{x, x'\} \in \eta_{1/2}\}. \end{aligned} \tag{9.5}$$

Therefore, (9.4) follows, provided we show that

$$\mathcal{L}(\eta_{1/2}) = \mathcal{L}(\eta_{1/2} + 1). \tag{9.6}$$

We will prove (9.6) by showing that if $\mathcal{L}(\eta_0) = \mu_\rho$ with $\rho = \frac{4\beta}{(1+\beta)^2}$, then $\mathcal{L}(\eta_{1/2})$ is a Bernoulli product measure on the set of all nearest neighbor pairs of integers. Note that, for $x \in \mathbb{Z}_{\text{even}}$, the event $\{x, x \pm 1\} \in \eta_{1/2}$ means that the arrow from $(x, 0)$ to $(x \pm 1, 1)$ is used by a random walker. Since $\mathcal{L}(\eta_0)$ is a product measure, arrows going out of different $x, x' \in \mathbb{Z}_{\text{even}}$ are obviously independent. Thus, it suffices to show that, for $x \in \mathbb{Z}_{\text{even}}$, the events $\{x, x - 1\} \in \eta_{1/2}$ and $\{x, x + 1\} \in \eta_{1/2}$ are independent. Now, for $x \in \mathbb{Z}_{\text{even}}$,

$$\mathbb{P}[\{x, x-1\} \in \eta_{1/2} \text{ and } \{x, x+1\} \in \eta_{1/2}] = \rho\beta, \tag{9.7}$$

while

$$\mathbb{P}[\{x, x-1\} \in \eta_{1/2}] = \mathbb{P}[\{x, x+1\} \in \eta_{1/2}] = \rho\left(\frac{1-\beta}{2} + \beta\right). \tag{9.8}$$

Thus, we obtain the desired independence, provided that $\rho\beta = (\rho\frac{1+\beta}{2})^2$, which has $\rho = \frac{4\beta}{(1+\beta)^2}$ as its unique nonzero solution. □

PROOF OF PROPOSITION 1.14(iii). By going to a subsequence if necessary, we may assume that $\mathcal{U}_\beta(x_n, t_n) \to \mathcal{A}$ for some $\mathcal{A} \subset \mathcal{U}_\beta$. Since all paths in $\mathcal{A}$ start at $(*, -\infty)$, $\mathcal{A} \subset \mathcal{U}_\beta(*, -\infty)$. To prove the other inclusion, it suffices to show that, for each $\pi \in \mathcal{U}_\beta(*, -\infty)$ and $t \in \mathbb{Z}_{\text{even}}$, for $n$ sufficiently large, we can find $\pi' \in \mathcal{U}_\beta(x_n, t_n)$ such that $\pi' = \pi$ on $[t, \infty) \cap \mathbb{Z}$. By hopping, it suffices to show that, for each even $N > 0$ and $t \in \mathbb{Z}_{\text{even}}$, there exists $n_0$ such that, for all $n \geq n_0$,

$$\begin{aligned} &[-N, N] \cap \{\pi(t) : \pi \in \mathcal{U}_\beta(x_n, t_n)\} \\ &\qquad \supset [-N, N] \cap \{\pi(t) : \pi \in \mathcal{U}_\beta(*, -\infty)\}. \end{aligned} \tag{9.9}$$

Let $\hat{l} := \hat{l}_{(-N-1,t)}$ and $\hat{r} := \hat{r}_{(N+1,t)}$ be the dual left-most and right-most paths in $\hat{\aleph}_\beta$ started from $(-N-1,t)$ and $(N+1,t)$, respectively. By the strong law of large numbers, almost surely,

$$\lim_{s\to-\infty} \frac{\hat{l}(s)}{-s} = \beta \quad \text{and} \quad \lim_{s\to-\infty} \frac{\hat{r}(s)}{-s} = -\beta. \tag{9.10}$$

Therefore, by our assumptions on $(x_n, t_n)$, we have $\hat{r}(t_n) < x_n < \hat{l}(t_n)$ for $n$ sufficiently large. Since forward paths and dual paths cannot cross, it follows that eventually $l_{(x_n,t_n)}(t) \le -N$ and $N \le r_{(x_n,t_n)}(t)$. Therefore, any path $\pi \in \mathcal{U}_\beta(*, -\infty)$ passing through $[-N, N] \times \{t\}$ must cross either $l_{(x_n,t_n)}$ or $r_{(x_n,t_n)}$. Since we can hop onto $\pi$ from either $l_{(x_n,t_n)}$ or $r_{(x_n,t_n)}$, formula (9.9) follows. □

Proof of Proposition 1.14(i) and (ii). It is not hard to see that $\mathcal{L}(\eta_0^{(0,-\infty)})$ is the maximal invariant law of $\eta$ with respect to the usual stochastic order. Proposition 1.14(iii) implies that

$$\mathbb{P}[\eta_{2n}^{(0,0)} \in \cdot\,] = \mathbb{P}[\eta_0^{(0,-2n)} \in \cdot\,] \underset{n\to\infty}{\Longrightarrow} \mathbb{P}[\eta_0^{(0,-\infty)} \in \cdot\,]. \tag{9.11}$$

Using monotonicity, it is easy to see from (9.11) that $\mathcal{L}(\eta_0^{(0,-\infty)})$ is the limit law of $\eta_{2n}$ as $n \to \infty$ for any nonempty initial state $\eta_0$. In particular, this implies that $\mathcal{L}(\eta_0^{(0,-\infty)})$ is the unique invariant law of $\eta$ that is concentrated on nonempty states, and therefore, by Lemma 9.1, $\mathcal{L}(\eta_0^{(0,-\infty)}) = \mu_\rho$.

Part (ii) is now a consequence of the reversibility of $\mu_\rho$ as formulated in Lemma 9.1. □

9.2. *The backbone of the branching-coalescing point set.* In this section we prove Proposition 1.15.

Proof of Proposition 1.15(iii). This can be proved by the same arguments as in the proof of Proposition 1.14(iii), except we now need Proposition 1.4 to hop between paths in the net. □

We will derive parts (i) and (ii) of Proposition 1.15 from their discrete counterparts, by means of the following lemma.

Lemma 9.2 (Convergence of the backbone). *If $\beta_n, \varepsilon_n \to 0$ with $\beta_n/\varepsilon_n \to 1$, then*

$$\mathcal{L}(S_{\varepsilon_n}(\mathcal{U}_{\beta_n}(*, -\infty))) \underset{n\to\infty}{\Longrightarrow} \mathcal{L}(\mathcal{N}(*, -\infty)). \tag{9.12}$$

Proof. By going to a subsequence if necessary, using Theorem 1.1, we may assume that

$$\mathcal{L}(S_{\varepsilon_n}(\mathcal{U}_{\beta_n},\mathcal{U}_{\beta_n}(*,-\infty))) \underset{n\to\infty}{\Longrightarrow} \mathcal{L}(\mathcal{N},\mathcal{A}), \tag{9.13}$$

where $\mathcal{N}$ is the Brownian net and $\mathcal{A}\subset\mathcal{N}$. Since all paths in $\mathcal{A}$ start in $(*,-\infty)$, obviously $\mathcal{A}\subset\mathcal{N}(*,-\infty)$. To prove the other inclusion, it suffices to show that (using notation introduced in the proof of Theorem 1.11)

$$\mathcal{N}(*,-\infty)|_t^\infty = \mathcal{A}|_t^\infty, \tag{9.14}$$

for all $t\in\mathbb{R}$. As a first step, we will show that

$$\{\pi(t):\pi\in\mathcal{N}(*,-\infty)\} = \{\pi(t):\pi\in\mathcal{A}\}. \tag{9.15}$$

The inclusion $\supset$ is clear. Taking the limit in Proposition 1.14(i), we see that, for all $t\in\mathbb{R}$, $\{\pi(t):\pi\in\mathcal{A}\}$ is a Poisson point set with intensity 2. On the other hand, taking the limit in Proposition 1.12, we see that $\{\pi(t):\pi\in\mathcal{N}(*,-\infty)\}$ is a translation invariant point set, also with intensity 2. Hence, (9.15) follows.

Since the inclusion $\supset$ in (9.14) is clear, it suffices to show that $\mathcal{N}(*,-\infty)|_t^\infty$ and $\mathcal{A}|_t^\infty$ are equal in law. Let $P$ be the random set in (9.15). By Lemma 8.3, $\mathcal{N}(*,-\infty)|_t^\infty = \mathcal{N}(P\times\{t\})$. By the independence of $\mathcal{N}|_{-\infty}^t$ and $\mathcal{N}|_t^\infty$ (see the proof of Theorem 1.11) and what we have just proved, it follows that $\mathcal{N}(P\times\{t\})$ is equally distributed with $\mathcal{N}(P'\times\{t\})$, where $P'$ is a Poisson point set with intensity 2, independent of $\mathcal{N}$. By Proposition 8.4, the law of $\mathcal{A}|_t^\infty$ is the same as that of $\mathcal{N}(P'\times\{t\})$, and we are done. □

Proof of Proposition 1.15(i) and (ii). The statements follow by a passage to the limit in Propositions 1.14(i) and (ii), using Lemma 9.2. □

## APPENDIX: DEFINITIONS OF PATH SPACE

In this appendix we compare the definition of the path space $\Pi$ and its topology used in the present paper with the definitions used in [7, 8]. Let $\mathcal{P}$ be the space of all functions $\pi:[\sigma_\pi,\infty]\to[-\infty,\infty]$, with $\sigma_\pi\in[-\infty,\infty]$, such that $t\mapsto\Theta_1(\pi(t),t)$ is continuous on $(-\infty,\infty)$. For $\pi_1,\pi_2\in\mathcal{P}$, define $d(\pi_1,\pi_2)$ by (1.5) and define $d'$ in the same way, but with the supremum over all $t\geq\sigma_{\pi_1}\wedge\sigma_{\pi_2}$ replaced by an unrestricted supremum over all $t\in\mathbb{R}$. Call two elements $\pi_1,\pi_2\in\mathcal{P}$ $d$-equivalent (resp. $d'$-equivalent) if $d(\pi_1,\pi_2)=0$ [resp. $d'(\pi_1,\pi_2)=0$], and let $\Pi$ (resp. $\Pi'$) denote the spaces of $d$-equivalence classes (resp. $d'$-equivalence classes) in $\mathcal{P}$. Then $(\Pi,d)$ is in a natural way isomorphic to the set of paths defined in Section 1.2, while $(\Pi',d')$ is the space of paths used in [7, 8]. The difference between these two spaces is small. Indeed, two paths $\pi_1,\pi_2$ are $d$-equivalent if and only if

$$\sigma_{\pi_1}=\sigma_{\pi_2} \quad\text{and}\quad \pi_1(t)=\pi_2(t) \qquad \forall\sigma_\pi\leq t<\infty, \tag{A.1}$$

while they are $d'$-equivalent if and only if

$$\text{(A.2)} \qquad \sigma_{\pi_1} = \sigma_{\pi_2} < \infty \quad \text{and} \quad \pi_1(t) = \pi_2(t) \qquad \forall \sigma_\pi \le t < \infty.$$

Thus, the only difference between $\Pi$ and $\Pi'$ is that, while the former has only one path with starting time $\infty$, the latter has a one-parameter family $(\pi^{(r)})_{r\in[-\infty,\infty]}$ of such paths, given by

$$\text{(A.3)} \qquad \sigma_{\pi^{(r)}} := \infty, \qquad \pi^{(r)}(\infty) := r \qquad (r \in [-\infty,\infty]).$$

A sequence of paths $\pi_n$ converges in $d'$ to the limit $\pi^{(r)}$ if and only if $\sigma_{\pi_n} \to \infty$ and $\pi_n(\sigma_{\pi_n}) \to r$. Both the spaces $(\Pi, d)$ and $(\Pi', d')$ are complete and separable, and the former is the continuous image of the latter under a map that identifies the family of paths $(\pi^{(r)})_{r\in[-\infty,\infty]}$ with a single path.

Of course, it is more natural to identify all paths starting at infinity. In fact, it seems that the authors of [8] used the metric in (1.5) in earlier versions of their manuscript, but then by accident dropped the restriction that $t \ge \sigma_{\pi_1} \wedge \sigma_{\pi_2}$ in the supremum [C. M. Newman, personal communication].

**Acknowledgments.** We thank Charles Newman, Krishnamurthi Ravishankar and Emmanuel Schertzer for answering questions about [8] and [12]. We thank Peter Mörters for answering questions about special points of Brownian motion and the referee for helpful suggestions leading to a better exposition.

R. Sun thanks ÚTIA and J. M. Swart thanks EURANDOM for their hospitality during short visits, which were supported by RDSES short visit grants from the European Science Foundation.

R. Sun was a postdoc at EURANDOM from Oct 2004 to Oct 2006, when this work was completed.

TU Berlin
MA 7-5, Fakultät II—Institut für Mathematik
Strasse des 17. Juni 136
10623 Berlin
Germany
E-mail: sun@math.tu-berlin.de

ÚTIA AV ČR
Pod vodárenskou věží 4
18208 Praha 8
Czech Republic
E-mail: swart@utia.cas.cz
URL: http://staff.utia.cas.cz/swart/